\documentclass[leqno, myheadings, twoside]{amsart}
\usepackage{amsmath, amsthm, amssymb, amscd, amsxtra,graphicx}
\usepackage{latexsym, amsfonts}
\usepackage{url}
\usepackage{texdraw}
\usepackage{epsfig}
\def\1-1{(1.1)}

\makeatletter \@addtoreset{equation}{section}

\makeatletter \renewcommand{\@biblabel}[1]{#1.}

\theoremstyle{remark}

\begin{document}
\title [Asymptotic formula for Steklov eigenvalues]
{The Weyl-type asymptotic formula for biharmonic Steklov eigenvalues
with Neumann boundary condition in Riemannian manifolds}
\author{Genqian Liu}

\subjclass{35P20, 58J50, 58C40, 65L15\\   {\it Key words and
phrases}.
 biharmonic equation, Steklov eigenvalue, asymptotic formula, Riemannian manifold}

\maketitle Department of Mathematics, Beijing Institute of
Technology,
 Beijing, the People's Republic of China.
 \ \
E-mail address:  liugqz@bit.edu.cn

\vskip 0.46 true cm

\vskip 15 true cm

\begin{abstract} \  In this paper, by a new method we establish the Weyl-type asymptotic formula
 for the counting function of biharmonic Steklov eigenvalues with
    Neumann boundary condition
 in a bounded domain of an $n$-dimensional
   Riemannian manifold.
   \end{abstract}

\vskip 1.39 true cm

\section{ Introduction}

\vskip 0.45 true cm

  Let $(\mathcal{M},g)$ be an oriented Riemannian manifold of
dimension $n$ with a positive definite metric tensor $g$, and let
$D\subset \mathcal{M}$ be a bounded domain with $C^2$-smooth
boundary $\partial D$. Assume $\varrho$ is a non-negative bounded
function defined on $\partial D$.
  Consider the following
 biharmonic Steklov eigenvalue problem with Neumann boundary
 condition:
         \begin{eqnarray} \label{1-1}  \left\{\begin{array}{ll}\triangle^2_g u=0
      \quad \;\; &\mbox{in}\;\; D,\\
  \frac{\partial u}{\partial \nu}=0 \;\; & \mbox{on}\;\; \partial D,\\
   \frac{\partial (\triangle_g u)}{\partial \nu}-\lambda^3 \varrho^3 u=0
   \;\;& \mbox{on}\;\;  \partial D, \end{array} \right. \end{eqnarray}
   where $\nu$ denotes the inward unit
   normal vector to $\partial D$,
   and $\triangle_g$ is the Laplace-Beltrami operator
   given in local coordinates by
   \begin{eqnarray*} \triangle_g =\frac{1}{\sqrt{|g|}}\sum_{i,j=1}^n
   \frac{\partial}{\partial x_i} \left( \sqrt{|g|}\,g^{ij}
   \frac{\partial}{\partial x_j}\right).\end{eqnarray*}
Here $| g | : =  det(g_{ij})$ is the determinant of the metric
tensor, and $g^{ij}$ are the components of the inverse of the metric
tensor $g$.

The problem (\ref{1-1}) was first discussed in 1968 by J. R. Kuttler
and V. G. Sigillito (see \cite{KS1}) since it describes the
deformation $u$ of the linear elastic supported plate $D$ under the
action of the transversal exterior force $f(x) = 0$, $x\in D$ (for
example, when the weight of the body $D$ is the
 only body force) with Neumann boundary condition $\frac{\partial u}{\partial
\nu}\big|_D =0$ (see, \cite{Vi} or p.$\,$32 of \cite{TG}).

\vskip 0.12 true cm

 It is well-known that the problem (\ref{1-1}) has nontrivial solutions $u$ only for a discrete set
 of $\lambda^3 = \lambda_k^3$, which are called biharmonic Steklov
  eigenvalues with Neumann boundary condition.
  Let us enumerate the eigenvalues in increasing order:
\begin{eqnarray*}  & 0=\lambda_0^3<\lambda_1^3 \le \lambda_2^3 \le  \cdots \le
   \lambda_k^3\le\cdots,
   \end{eqnarray*}
   where each eigenvalue is counted as many times as its multiplicity.
    The corresponding eigenfunctions $u_{01}, u_{02}, u_1, u_2, \cdots, u_k, \cdots$
   form a complete orthonormal basis in $L^2_{\varrho} (\partial D)$ (see, Proposition 3.5), where
   $u_{01}(x)=1$, $u_{02}(x)= \int_{D} F(x,y)dR_y$ on $x\in \partial D$,
   and $F(x,y)$ is Green's function in $D$ with
   Neumann boundary condition.
It is clear that $\lambda_k^3$ can be characterized variationally as
\begin{eqnarray*}
 \lambda_0=0,\quad \, \lambda_1^3=\frac{\int_D |\triangle_g u_1|^2 dR} {
\int_{\partial D} \varrho^3 u_1^2 ds} =
  \inf_{\underset{\int_{\partial D} \varrho^3 v\, ds =0,\;
   \int_{\partial D}\varrho^3 u_{{}_{02}}v\, ds =0}
  {v\in Lip(\bar D)\cap H^2 (D)}} \,
  \frac{\int_D |\triangle_g v|^2 dR} {
\int_{\partial D} \varrho^3  v^2 ds},\;\;\qquad \qquad \\
 \lambda_k^3=\frac{\int_D |\triangle_g u_k|^2 dR} {
\int_{\partial D} \varrho^3  u_k^2 ds} =
  \max_{\underset{\int_{\partial D}
 \varrho^3 v\,ds=0,\; \int_{\partial D} \varrho^3 u_{{}_{02}} v \,ds
 =0\},\;  codim(\mathcal {F})=k+1}
  {\mathcal{F}\subset \{v\big|v\in Lip(\bar D)\cap H^2 (D),}}
   \inf_{v\in \mathcal{F}}
\frac{\int_D |\triangle_g v|^2 dR} { \int_{\partial D} \varrho^3 v^2
ds}, \quad  k=2,3, \cdots
\end{eqnarray*} where
 $H^m(D)$ is the Sobolev space, $Lip(\bar D)$ is the set of Lipschitz functions on
 $\bar D$,
 and where $dR$ and $ds$ are the Riemannian
elements of volume and area on $D$ and $\partial D$, respectively.

 The problem (\ref{1-1}) is also important in biharmonic analysis because the set of the eigenvalues
 for the biharmonic Steklov problem is the same as the set
 of eigenvalues of the well-known ``Dirichlet to normal derivative of
 Laplacian'' map for biharmonic equation
 (This map associates each function $u$ defined on the boundary
  $\partial D$ to the normal derivative $\frac{\partial (\triangle_g u)}{\partial \nu}$
  of $\triangle_g u$,
   where the biharmonic function $u$ in $D$ is uniquely
 determined by $u\big|_{\partial D}$  and
  $(\partial u/\partial \nu)\big|_{\partial D}=0$).

In the general case the eigenvalues $\lambda_k^3$ can not be
evaluated explicitly. In particular, for large $k$ it is difficult
to calculate them numerically. In view of the important
applications, one is interested in finding the asymptotic formula
for $\lambda_k^3$ as $k\to \infty$. Let us introduce the counting
function $A(\tau)$ defined as the number of eigenvalues
$\lambda_k^3$ less than or equal to a given $\tau^3$. Then our
asymptotic problem is reformulated as the study of the asymptotic
behavior of
 $A(\tau)$ as $\tau\to +\infty$.

In order to better understand our problem (\ref{1-1}) and its
asymptotic behavior, let us mention the Steklov eigenvalue problem
for harmonic equation
\begin{eqnarray} \label{1-3}  \left\{\begin{array}{ll} \triangle_g v=0&
\quad \, \mbox{in}\;\; D,\\
\frac{\partial v}{\partial \nu}+ \iota \varrho v =0 &\quad \,
\mbox{on}\;\; \partial D,
\end{array} \right.\end{eqnarray}
 where $\iota$ is a real number. This problem was  introduced by
 M. W. Steklov for bounded domains in the plane in \cite{St}. His motivation
 came from physics. The function $v$ represents
 the steady state temperature on $D$ such that
 the flux on the boundary is proportional
 to the temperature.
 For the harmonic Steklov eigenvalue problem (\ref{1-3}), in a special case
 in two dimensions,
 {\AA}. Pleijel \cite{Pl} outlined an investigation of the
 asymptotic behavior of both eigenvalues $\iota_k$ and the eigenfunctions $v_k$.
 In 1955, L. Sandgren \cite{Sa} established the asymptotic formula of the counting function $B(\tau)=
 \#\{\iota_{{}_k}\big|\iota_{{}_k} \le \tau\}$:
   \begin{eqnarray} \label{1-4}   B(\tau)\sim
   \frac{\omega_{n-1}\tau^{n-1}}{(2\pi)^{n-1}}
 \int_{\partial D} \varrho^{n-1} ds \quad \;\mbox{as}\;\;
 \tau\to +\infty,\end{eqnarray}
 i.e.,
 \begin{eqnarray*}  \lim_{\tau \to +\infty} \; \frac{B(\tau)}{\tau^{n-1}}=
 \frac{\omega_{n-1}}{(2\pi)^{n-1}}
 \int_{\partial D} \varrho^{n-1} ds,\end{eqnarray*}
 where $\omega_{n-1}$ is the volume of the $(n-1)$-dimensional unit
 ball, $ds$ is the Riemannian element of area on $\partial
 D$.
This asymptotic behaviors is motivated by the similar one for the
 eigenvalues of the Dirichlet-Laplacian. The classical result for the
 Dirichlet-Laplacian on
 smooth domain $D$ is
Weyl's formula (see \cite{We1}, \cite{We2} or \cite{CH}):
\begin{eqnarray}  \label{1-5}   N_D(\tau,D)\sim
 \frac{\omega_n}{(2\pi)^{n}} \big(\mbox{vol}(D)\big) \tau^{n/2}
\quad \mbox{as}\;\; \tau\to +\infty,\end{eqnarray}
 where $N(\tau, D)=\# \{\mu_k \le \tau\}$ and
 $\mu_k$ is the $k$-th Dirichlet eigenvalue for $D$.

The study of asymptotic behavior
 for the biharmonic Steklov eigenvalues with Neumann boundary
condition is much more difficult than that for the harmonic
 Steklov eigenvalues. It has been a tempting and challenging problem
 in the past 40 years.
  The main stumbling
   block that lies in its way is the
    estimates for the different kinds of Steklov
  eigenvalues  corresponding to the different kinds of boundary conditions.
For the simpler biharmonic Steklov eigenvalue problem with Dirichlet
boundary condition, the
 author established the leading term asymptotic formula of the eigenvalues (see, \cite{Liu}).
\vskip 0.24 true cm

In this paper, for the biharmonic Steklov eigenvalues
  with Neumann boundary condition, by a new method we establish the Weyl-type
asymptotic formula of the counting function. The main results are
the following:

\vskip 0.25 true cm

 \noindent  {\bf Theorem 1.1.} \ \  {\it
Let $(\mathcal{M},g)$ be an $n$-dimensional oriented Riemannian
manifold, and let $D\subset M$ be a bounded domain with
$C^{1}$-smooth boundary $\partial D$.  Then
\begin{eqnarray} \label{1-6} A(\tau)=
\frac{\omega_{n-1}\tau^{n-1}}{(\sqrt[3]{16}\pi)^{n-1}}
 \int_{\partial D} \varrho^{n-1} ds + o( \tau^{n-1})\quad \, \; \mbox{as}\;\; \tau\to +\infty,\end{eqnarray}
where $A(\tau)$ is defined as before}.

\vskip 0.25 true cm

 \noindent  {\bf Corollary 1.2.} \  {\it
 Under hypothesis Theorem 1.1, if  $\varrho\equiv 1$ on
$\partial D$ for problem (\ref{1-1}),
 then \begin{eqnarray} \label{1-7}  \lambda_k \sim
 \sqrt[3]{16}\pi \left(\frac{k+2}{\omega_{n-1} (\mbox{vol}(\partial D))}\right)^{1/(n-1)}
 \quad \, \mbox{as}\;\; k\to +\infty.\end{eqnarray}}

\vskip 0.25 true cm

 However, when the boundary of a bounded domain is smooth, we have the
following Weyl-type asymptotic formula with a better remainder
estimate:

\vskip 0.25 true cm

 \noindent  {\bf Theorem 1.3.} \ \  {\it
Let $(\mathcal{M},g)$ be a smooth, $n$-dimensional oriented
 Riemannian manifold, and let $D\subset M$ be a bounded domain with
smooth boundary $\partial D$.  Then
\begin{eqnarray} \label{1;-6;} A(\tau)\sim
\frac{\omega_{n-1}\tau^{n-1}}{(\sqrt[3]{16}\pi)^{n-1}}
 \int_{\partial D} \varrho^{n-1} ds + O( \tau^{n-2})\quad \, \; \mbox{as}\;\; \tau\to +\infty,\end{eqnarray}
where $A(\tau)$ is defined as before}.

\vskip 0.28 true cm

 The proofs of our main results uses four   key techniques:
The first technique is the compact trace lemmas for the domain which
is the union
 of a finite number of Lipschitz images of cubes.
 The second technique is to
 give the explicit formula for the different kinds of
 biharmonic Steklov eigenvalues and eigenfunctions in a cube of ${\Bbb R}^n$
 (by the method of separation variables we seek the product form of eigenfunctions, one of
 factors is the Dirichlet eigenfuction or Neumann eigenfunction, see Section
 4).   Then we can  use the well-known variational methods,
which H. Weyl \cite{We4}  and R. Courant  and D. Hilbert \cite{CH}
 have employed in the case of the membrane to give the asymptotic formulas for
 the two kinds of the  Steklov eigenvalues in the cube.
  The third technique is
 put the biharmonic Steklov problem into an
 abstract Hilbert space theory. That is,
  we first make a division of $\bar D$ into subdomains. From this division we construct
two Hilbert spaces ${\mathcal{K}}^0$ and ${\mathcal{K}}^d$ and
isometric mappings of ${\mathcal{K}}^0$ into $\mathcal{K}$ and
$\mathcal{K}$ into ${\mathcal{K}}^d$. Those of subdomains situated
at the boundary $\Gamma_\varrho$ we can map on cylinders of type
treated in Section 4. In a sufficiently fine division of these, the
variant of the $g^{ik}$ and $\varrho$ will be small and they can be
replaced by constants. Then, we can estimate the asymptotic behavior
of the eigenvalues.
 By means of the results of Section 5, we get the asymptotic formula with leading asymptotic for $A(\tau)$.
  Finally, applying Theorem 1.1 and a standard technique based on
  the asymptotic behavior of spectral function of pseudodifferential
  operator (see, p.$\;$162 of \cite{Sh},  \cite{Ho} or
  \cite{Ta2}), we obtain the desired result of Theorem 1.3 (with a better remainder estimate).

\vskip 1.39 true cm

\section{Compact trace Lemmas}

\vskip 0.45 true cm

An $n$-dimensional cube in ${\Bbb R}^n$ is the set $\{x\in {\Bbb
R}^n\big| 0\le x_k \le a, \, k=1, \cdots, n\}$.
 A set $D\subset {\Bbb R}^n$ is said to be a Lipschitz image
 of a set $\Omega\subset {\Bbb R}^n$ (see \cite{Sa}) if there is a one-to-one
 map from $\Omega$ to $D$ defined by
    \begin{eqnarray} \label{2-1} x=\Psi(x'), \; \; x'\in \Omega \end{eqnarray}
      satisfying a Lipschitz condition
\begin{eqnarray}  \label{2-2}  c^{-1} |x'-y'|\le |\Psi(x')-\Psi(y')|\le c|x'-y'|\end{eqnarray}
  for some constant $c$ and all $x'$ and $y'$ in $\Omega$.
  ($|x|=(x_1^2 +\cdots +x_n^2)^{1/2}$, $\, |x'|=({x'_1}^2+\cdots
  +{x'_n}^2 )^{1/2}$).

  A set $D\subset {\Bbb R}^n$  is said to be star-shaped with respect to a point $x^0$ if
  $x\in D$  implies that the closed segment $\{(1-t)x^0 +tx\big|0\le t\le 1\}$
  is contained in $D$. Now assume
  that $D$ is a bounded domain in ${\Bbb R}^n$ and that the
  closed domain $\bar D$ is star-shaped with
  regard to all points in an open neighborhood of a point $x^0\in D$.
   We can assume $x^0=(0, \cdots, 0)$. In this section,
   $\|x\|$ denotes an arbitrary norm in ${\Bbb R}^n$ with the usual
   properties of a norm, that is
   \begin{eqnarray*}  a) \quad  \|x\|=0\Leftrightarrow x=0, \,\quad \;\;  b)\quad
    \|x+y\|\le \|x\|+\|y\| \,\quad \; \; \mbox{and} \,\quad \; c) \quad
   \|ct\|=|t|\|x\|,\end{eqnarray*}
   where $t$ is a real number. Then evidently (see \cite{Sa}) there is a $\delta>0$
   such that $\bar D$ is star-shaped with respect to all points
 in $B_\delta =\{x\big|\|x\|<\delta\}$.
 Since $B_\delta$ is open, it is clearly that $x\in B_\delta$
 and $y\in \bar D$ implies that all the inner points of the
 segment $\{(1-t)x +t y\big|0\le t\le 1\}$ belong to $D$.

\vskip 0.16 true cm

 \noindent  {\bf Lemma 2.1 (Sandgren, p$\;$21 of \cite{Sa}).} \ \  {\it
    If a bounded domain $\bar D\subset {\Bbb R}^n$ is star-shaped with respect
    to all points in the open cube $\sum_\delta= \{x\in {\Bbb R}^n\big| \max_{1\le i\le n}
    |x_i| <\delta\}$, then $\bar D$ is a
    Lipschitz image of the cube $\bar D'$ (with the side-length $2a$) given by a transformation
    \begin{eqnarray} \label{2-3} x=\Psi(x'), \;\; x\in \bar D \;\; \mbox{and}\;\; x'\in \bar
    D'\end{eqnarray}
     satisfying the Lipschitz condition
     \begin{eqnarray}  \label{2-4} c^{-1} \|x'-y'\|\le \|\Psi(x')-\Psi(y')\|\le
     c\|x'-y'\|,\end{eqnarray}
     where $c= \max (3a/\delta, 3b^2/ \delta a)$ and $\,b=\max_{x\in \bar D} \|x\|$.}

\vskip 0.16 true cm

 Let $f$ be a real-valued function defined
in an open set $D$ in ${\Bbb R}^{n}$ ($n\ge 1$). For $y\in D$ we
call $f$ {\it real analytic at $y$} if there exist $a_\beta \in
{\Bbb R}^1$ and a neighborhood $U$ of $y$ (all depending on $y$)
such that $$ f(x)= \sum_\beta a_\beta (x-y)^\beta$$ for all $x$ in
$U$. We say $f$ is {\it real analytic in $D$}, if $f$ is real
analytic at each $y\in D$.

 From here up to Section 5, let $M$ be an
 $n$-dimensional Riemannian manifold with real analytic metric tensor $g$.
  We say that $\bar D$ is a Lipschitz image of a cube
  if it is contained in some coordinate neighborhood $U$ and its
  image $\bar D_1$ in ${\Bbb R}^n$  given by the coordinates of $U$
  is a Lipschitz image (see, previous definition) of a closed cube
  in ${\Bbb R}^n$.

A subset $\Gamma$ of $(\mathcal{M},g)$ is said to be an
$(n-1)$-dimensional smooth surface if  $\Gamma$ is nonempty and if
for every point $x$ in $\Gamma$, there is a smooth  diffeomorphism
of the open unit ball $B(0,1)$ in ${\Bbb R}^n$ onto an open
neighborhood $U$ of $x$ such that $B(0,1)\cap \{x\in {\Bbb
R}^n\big|x_n=0\}$ maps onto $U\cap \Gamma$.

Let $D$ together with its boundary  be transformed pointwise into
the domain $D'$ together with its boundary by equations of the form
\begin{eqnarray} \label {02-001} x'_i= x_i+ f_i(x_1, \cdots,
x_{n}), \quad \; i=1,2, \cdots, n. \end{eqnarray}
 where the functions $f_i$ and their first order derivatives are
 Lipschitz continuous throughout the domain, and they are less in absolute
 value than a small positive number $\epsilon$.
 Then we say that the domain $D$ is approximated by the domain $D'$ with the
 degree of accuracy $\epsilon$.

 It is
well-known (see, for example, p.$\;$133 of \cite{HA} or p.$\;$24 of
\cite{Sa})
 that every element $u$ in $Lip(\bar D)$ has partial derivatives  $\partial u/ \partial x_k$,
 $k=1, \cdots, n$, which are defined a.e. in $D$ and belong to $L^\infty (D)$. In particular,
 $Lip(\bar D)\subset H^1(D)$.

\vskip 0.12 true cm

A subset $\mathfrak{F}$ of $L^2(\partial D)$ is said to be
precompact if any infinite sequence $\{u_k\}$ of elements of
$\mathfrak{F}$ contains a Cauchy subsequence $\{u_{k'}\}$, i.e., one
for which
\begin{eqnarray*} \int_{\partial D} (u_{k'} -u_{l'})^2 ds\to 0, \quad
\mbox{when}\;\; k', l'\to \infty.\end{eqnarray*}

\vskip 0.25 true cm

 \noindent  {\bf Lemma 2.2.} \ \  {\it
    Let $\bar D\subset (M, g)$ be a Lipschitz image of a cube and let
    $\varrho$ be a non-negative function in $L^\infty (\partial D)$
    such that $\int_{\partial D} \varrho^3 ds > 0$.
     Assume $\mathfrak{M}$ is a set of functions $u$ in
 $\tilde N(D)=\{u\big|u\in Lip(\bar D)\cap H^2(D),\;
 \frac{\partial u}{\partial \nu}=0 \;\;\mbox{on}\;\;
 \partial D\}$ for which
  \begin{eqnarray} \label{-2-1}\int_D |\triangle u|^2 dR + \big(\int_{\partial D}
   \varrho^3 u\, ds \big)^2  \end{eqnarray}
  is uniformly bounded. Then the set $\{u\big|_{\partial D}: u\in \mathfrak{M}\}$
  is precompact in $L^2(\partial D)$. }

\vskip 0.28 true cm

\noindent  {\bf Proof.} \
 It follows from Theorem 2.3 of \cite{Sa} that
 there exists a constant $C_1 >0$, only depending on $D$ and
 $\varrho $, such that for all  $u\in Lip(\bar D)$,
  \begin{eqnarray} \label{2.001} \int_D u^2dR \le
  C_1 \left[ \int_D |\nabla_g u|^2 dR + \big(\int_{\partial D} \varrho^3  u \, ds\big)^2\right],\end{eqnarray}
where \begin{eqnarray*} \int_D |\nabla_g u|^2 dR = \int_{\phi(D)}
g^{ik} (x)
 \frac{\partial u}{\partial x_i}\,\frac{\partial u}{\partial x_k}\, \sqrt{|g|}dx,\end{eqnarray*}
  and $\phi(D)$ is the coordinate image of $D$.
  Put
 \begin{eqnarray} \qquad \; \label{2.13} \Lambda_1(D)= \inf_{\underset
 {\int_D |\nabla_g u|^2dR +
 (\int_{\partial D} \varrho^3 u \, ds)^2=1} {
 u\in \tilde N (D)} }
\;\frac{\int_D|\triangle_g u|^2 dR+\big(\int_{\partial D} \varrho^3
u \,
 ds\big)^2}{\int_D |\nabla_g u|^2dR
 +\big(\int_{\partial D} \varrho^3 u \, ds\big)^2}.\end{eqnarray}
In order to prove the existence of a minimizer to (\ref{2.13}),
  consider a minimizing sequence $u_m$ in the set $\tilde N(D)$,
 i.e.,
    \begin{eqnarray*}  \int_{D}|\triangle_g u_m|^2 dR
   +\big(\int_{\partial D} \varrho^3 u_m \,
  ds\big)^2  \to
    \Lambda_{1} (D)\quad \mbox{as}\;\, m\to +\infty\end{eqnarray*}
 with $\int_D |\nabla_g u_m|^2dx
 +\big(\int_{\partial D} \varrho^3  u_m \, ds\big)^2
 =1$. Thus, there is a constant $C_2>0$ such that
  \begin{eqnarray} \label{2-021} &&\|\triangle_g u_m\|_{L^2(D)}
   \le C_2, \quad  \int_D |\nabla_g u_m|^2dx\le C_2, \quad
 \big(\int_{\partial D} \varrho^3  u_m \, ds\big)^2
 \le C_2
   \end{eqnarray}
 for all $m\ge 1$. It follows from
 the {\it a priori} estimate for elliptic
equations  (see, for example, Proposition 7.2 of p.345 in \cite
{Ta1})  that there exists a constant $C_3>0$ depending only on
 $n, D$ such that
  \begin{eqnarray} \label {2-22} \| u_m\|_{H^2(D)}
  \le C_4(\|\triangle_g  u_m\|_{L^2(D)}
  +\|u_m\|_{H^1(D)}).\end{eqnarray}
From this, (\ref{2.001}) and (\ref{2-021}), we have that
\begin{eqnarray*} \|u_m\|_{H^2(D)}\le C_4 \;\;\mbox{for all}\;\; m.\end{eqnarray*}
  By the Banach-Alaoglu
theorem we can then extract a subsequence, which we still call
$\{u_m\}$, converging  weakly in $H^2(D)$ to a limit $u$, and
strongly converging to $u$ in $L^2(D)$.  Since the functional
$\int_{D} |\triangle_g u|^2 dR$ is lower semicontinuous in the weak
$H^2(D)$ topology, we have
$$\int_{D} |\triangle_g u|^2 dR\le
 \underset {m\to \infty} {\underline{\lim}} \int_{D} |\triangle_g u_m|^2 dR,$$
 Since $u_m \to u$ weakly in
$H^2(D)$, we get that $u_m\to u$ strongly in $H^r (D_l)$ for any
$0<r<2$. Note that $\frac{\partial u_m}{\partial \nu}\big|_{\partial
D}=0$.  It follows that $\frac{\partial u}{\partial
\nu}\big|_{\partial D}=0$.
  Therefore $u\in \mathfrak {M}$ is a minimizer.

 We claim that $\Lambda_1(D)>0$. Suppose by contradiction that
 $\Lambda_1(D)=0$. Then
 \begin{eqnarray} \label{2,1} \left\{ \begin{array}{ll}  \triangle_g u=0 \;\;
  & \mbox{in}\;\; D,\\
  \frac{\partial u}{\partial \nu}=0  \;\; &\mbox{on}\;\; \partial D
    \end{array}\right. \end{eqnarray}
    and \begin{eqnarray} \label{2,2} \int_{\partial D} \varrho^3 u\,
    ds=0.\end{eqnarray}
    The boundary value problem (\ref{2,1}) implies that $u\equiv constant$ in $D$.
    By $\int_{\partial D} \varrho^3 \,
  ds >0$ and (\ref{2,2}), we then get $u=0$ in $D$. This contradicts the fact that
 $\int_D |\nabla_g u|^2 dR +\big(\int_{\partial D} \varrho^3 u\,
 ds\big)^2 =1$, and the claim is proved.

 By (\ref{2.13}), we obtain that for
 every $u\in \tilde N (D)$,
    \begin{eqnarray} \qquad \quad \;\; \label {2-16} \int_D |\nabla_g u|^2 dR +
  \big(\int_{\partial D} \varrho^3  u \, ds\big)^2
    \le
  \frac{1}{\Lambda_1(D)}\left( \int_D|\triangle_g u|^2 dR
 +\big(\int_{\partial D} \varrho^3
 u \, ds\big)^2\right).\end{eqnarray}
 It follows from (\ref{2.001}) and (\ref{2-16}) that $\mathfrak{M}$
 is a bounded set in $H^1(D)$. Since $\bar D\subset (M, g)$ is a
 Lipschitz image of a cube, it follows  from \cite{AGMT} (see also,
   Chs V, VI of \cite{EE}) that
the set $\{u\big|_{\partial D}: u\in \mathfrak{M}\}$
  is precompact in $L^2(\partial D)$.
     \quad $\square$

\vskip 0.25 true cm

 \noindent  {\bf Lemma 2.3.} \ \  {\it Let $(\mathcal{M},g)$ be a real
 analytic Riemannian manifold, and let  $\bar D\subset (M, g)$ be a
  Lipschitz image of a cube.
 Assume that $\Gamma_1$ is a portion of $\partial D$ and
  that $\Gamma_0$ is
     an $(n-1)$-dimensional $C^{2,\varepsilon}$-smooth surface in $\partial D$
 satisfying $\Gamma_0\subset \subset \partial D-\bar \Gamma_1$.
   Suppose that $\tilde\varrho$ is a non-negative function defined on
   $\partial D$ and assume that $\mathfrak{E}$
   (respectively, $\mathfrak{S}$) is a set of functions $u$ in
 $K^d(D)=\{u\big| u\in Lip(\bar D)\cap H^2(D),\, \frac{\partial u}
 {\partial\nu}=0 \;\;\mbox{on}\;\;\Gamma_1,\,
 \, u=\frac{\partial u}{\partial \nu}=0 \,\, \mbox{on}\;\; \Gamma_0\}$
 (respectively, $K(D)= \{u\big|u\in Lip(\bar D)\cap H^2(D), \;\, \frac{\partial u}{\partial \nu}=0 \,\;
  \mbox{on}\;\, \Gamma_1\cup \Gamma_0, \;\, \mbox{and}\,\;u=0 \,\;\mbox{on}\,\;
  \partial D-\Gamma_1\}$) for which
  \begin{eqnarray} \label{-2-2}\int_D |\triangle_g u|^2 dR
  + \left(\int_{\partial D} {\tilde \varrho}^3 u\, ds \right)^2\end{eqnarray}
  is uniformly bounded. Then the set $\{{u}\big|_{\partial D} :
  u\in \mathfrak{E}\}$ (respectively, $\{{u}\big|_{\partial D} :
  u\in \mathfrak{S}\}$) is precompact in $L^2(\partial D)$. }

\vskip 0.28 true cm

\noindent  {\bf Proof.} \   We only prove the case in $K^d(D)$
because the method is similar for the case in $K(D)$.
     It follows from Theorem 2.3 of \cite{Sa} that there exists a constant
 $C>0$, depending only on $D$ and $\tilde \varrho$, such that for every  $u\in Lip (\bar D)$,
  \begin{eqnarray} \label{2.01} \int_D u^2dR &\le&
  C \left[ \int_D |\nabla_g u|^2 dR + \left(\int_{\partial D}
   \tilde\varrho^3  u \, ds\right)^2\right].\end{eqnarray}

Let
 \begin{eqnarray}\label{2-19} \qquad\qquad \; \Lambda_{\tilde \varrho} (D)=
  \inf_{ \underset{\int_{D} |\nabla_g v|^2 dR+
  \left(\int_{\partial D} {\tilde \varrho}^3 v\, ds \right)^2 =1}{v\in
  k^d(D)}}\,
   \; \frac{\int_{D}|\triangle_g v|^2 dR+  \left(\int_{\partial D} {\tilde \varrho}^3 v\, ds \right)^2
 }{\int_{D} |\nabla_g v|^2 dR + \left(\int_{\partial D} {\tilde \varrho}^3 v\, ds \right)^2}.\end{eqnarray}
  In order to prove the existence of a minimizer to (\ref{2-19}),
  consider a minimizing sequence $v_m$ in $K^d(D)$,
 i.e.,
    \begin{eqnarray*}  \int_{D}|\triangle_g v_m|^2 dR  +\left(\int_{\partial D} {\tilde \varrho}^3 v_m\,
    ds \right)^2 \to
    \Lambda_{\tilde \varrho} (D)\quad \mbox{as}\;\, m\to +\infty\end{eqnarray*}
 with $\int_{D}  |\nabla_g v_m|^2 dR  +\left(\int_{\partial D} {\tilde \varrho}^3 v_m\, ds \right)^2  =1$.
Thus, there is a constant $\tilde C>0$ such that
  \begin{eqnarray} \label {2-21} &&\|\triangle_g v_m\|^2_{L^2(D)}
  +\left(\int_{\partial D} {\tilde \varrho}^3 v_m\, ds \right)^2 \le \tilde C,\\
   &&\|\nabla_g v_m\|^2_{L^2(D)}+\left(\int_{\partial D} {\tilde \varrho}^3 v_m\,
   ds \right)^2 \le \tilde C \quad \;
  \mbox{for all}\;\; m.\nonumber \end{eqnarray}
 Let $\{D_l\}$  be a sequence of Lipschitz domains such that
   $D_1\subset D_2\subset \cdots \subset D_l\subset \cdots \subset
   \subset (D\cup \Gamma_1\cup \Gamma_0)$, $\;\cup_{l=1}^\infty D_l= D$,
   and $\Gamma_1\cup \Gamma_0 \subset \partial D_l$ for
   all $l$.
It follows from
 the {\it a priori} estimate for elliptic
equations  (see, for example, the proof of Proposition 7.2 of \cite
{Ta1})  that there exists a constant $C'_l>0$ depending only on
 $n, D_l, D, \Gamma_1$ and $\Gamma_0$
  such that
  \begin{eqnarray} \label {2-22} \| v_m\|_{H^2(D_l)}
  \le C'_l(\|\triangle_g  v_m\|_{L^2(D)}
  +\|v_m\|_{H^1(D)}).\end{eqnarray}
From this, (\ref{2.01}) and (\ref{2-21}), we have that
\begin{eqnarray*} \|v_m\|_{H^2(D_l)}\le C''_l \;\;\mbox{for all}\;\; m,\end{eqnarray*}
where $C''_l$ is a constant.  For each $l$, by the Banach-Alaoglu
theorem we can extract a subsequence $\{v_{l,m}\}_{m=1}^\infty$ of
$\{v_m\}$, which converges weakly in $H^2(D_l)$ to a limit $u$, and
strongly converges to $u$ in $L^2(D_l)$. We may assume that
$\{v_{l+1,m}\}$ is a subsequence of $\{v_{l,m}\}$ for every $l$.
Then, the diagonal sequence $\{v_{l,l}\}$ converges weakly in $H^2$
to $u$, and strongly converges to $u$ in $L^2$, in every compact
subset $E$ of $D$.  It is obvious that $\|\nabla_g u\|^2_{L^2
(D)}+\left(\int_{\partial D} {\tilde \varrho}^3 u\, ds \right)^2
   =1$. Since the functional $\int_{D_l} |\triangle_g u|^2 dR$ is lower
semicontinuous in the weak $H^2(D_l)$ topology, we have
$$\int_{D_l} |\triangle_g u|^2 dR\le
 \underset {k\to \infty} {\underline{\lim}} \int_{D_l} |\triangle_g v_{k,k}|^2 dR,$$
 so that
 \begin{eqnarray*} \int_D |\triangle_g u|^2 dR &=&\lim_{l\to \infty} \int_{D_l}
 |\triangle_g u|^2 dR \le \lim_{l\to \infty} \left(\underset {k\to
 \infty}{\underline{\lim}}
 \int_{D_l} |\triangle_g v_{k,k}|^2 dR\right)\\
 &\le& \lim_{l\to \infty} \left(\underset {k\to
 \infty}{\underline{\lim}}
 \int_{D} |\triangle_g v_{k,k}|^2 dR\right)
 =\underset {k\to
 \infty}{\underline{\lim}}
 \int_{D} |\triangle_g v_{k,k}|^2 dR.\end{eqnarray*}
In addition, for each fixed $l$, since $v_{k,k} \to u$ weakly in
$H^2(D_l)$, we get that $v_{k,k}\to u$ strongly in $H^r (D_l)$ for
any $0<r<2$. Note that $\frac{\partial v_{k,k}}{\partial
\nu}\big|_{\Gamma_1}=0$ and
 $v_{k,k}\big|_{\Gamma_0}=\frac{\partial v_{k,k}}{\partial
\nu}\big|_{\Gamma_0}=0$. It follows that $\frac{\partial u}{\partial
\nu}\big|_{\Gamma_1}=0$ and $u\big|_{\Gamma_0}=\frac{\partial
u}{\partial \nu}\big|_{\Gamma_0}=0$.
  Therefore $u\in K^d(D)$ is a minimizer.

 We claim that  $\Lambda_{\tilde \varrho}(D)>0$.
Suppose by contradiction that
\begin{eqnarray*}\Lambda_{\tilde \varrho}(D)= \frac{\int_{D}|\triangle_g
u|^2 dR+ \left(\int_{\partial D} {\tilde \varrho}^3 u\, ds
\right)^2}{\int_{D} |\nabla_g u|^2 dR +\left(\int_{\partial D}
 {\tilde \varrho}^3 u\, ds \right)^2}= 0.\end{eqnarray*} It follows that
$\triangle_g u=0$ in $D$. Since the coefficients of the Laplacian
are real analytic in $D$,
  and since $\Gamma_0$ is a $C^{2,\varepsilon}$-smooth surface,
  we find with the aid of the
  regularity for elliptic equations (see,
  Theorem A of \cite{MN}, \cite{Mo} or \cite{ADN}) that
  $u$ is $C^{2,\varepsilon}$-smooth
  up to the partial boundary $\Gamma_0$.
 Note that $u=\frac{\partial
  u}{\partial \nu}=0$ on $\Gamma_0$.
 Applying Holmgren's uniqueness theorem (see, Corollary 5 of p.$\,$39 in \cite{Ra})
  for the real analytic
 elliptic equation $\triangle_g u=0$ in $D$, we get that $u\equiv 0$ in
 $D$. This
     contradicts the fact $\int_{D}|\nabla_g u|^2 dR+ \left(\int_{\partial D}
\varrho^3 u\, ds \right)^2=1$, and the claim is proved. Therefore we
have that
 \begin{eqnarray}\label{2-23} &&
\int_{D}|\nabla_g u|^2 dR+  \left(\int_{\partial D} {\tilde
\varrho}^3 u\, ds \right)^2\\&& \quad \quad \le
\frac{1}{\Lambda_{\tilde
 \varrho}} \left[\int_{D} |\triangle_g u|^2 dR+
  \left(\int_{\partial D} {\tilde \varrho}^3 u\, ds \right)^2\right]
    \quad \;\mbox{for}\;\; u\in K^d(D).\nonumber\end{eqnarray}
According to the assumption,
 there is a constant $C''$ such that  \begin{eqnarray*}
   \label{2-24}\|\triangle_g u
\|_{L^2(D)}^2 + \left(\int_{\partial D} {\tilde \varrho}^3 u\, ds
\right)^2 \le C''\quad \; \mbox{for all}\;\; u \in
 \mathfrak{E},\end{eqnarray*}
 and hence \begin{eqnarray}
 \label{2-25}\|\nabla_g u
\|_{L^2(D)}^2 + \left(\int_{\partial D} {\tilde \varrho}^3 u\, ds
\right)^2 \le C''\quad \; \mbox{for all}\;\; u \in
 \mathfrak{E}.\end{eqnarray}
   Combining this and (\ref{2.01}), we have
   $$\|u\|_{H^1(D)}\le C, \quad\mbox{for all}\;\; u\in \mathfrak{E},$$
   which implies that $\{u\big|_{\partial D} : u\in {\mathfrak{E}}\}$
     is precompact in $L^2(\partial D)$. \ \ $\square$

\vskip 0.25 true cm

 \noindent  {\bf Corollary 2.4.} \ \  {\it Lemma 2.3 is still true if we exchange
 (\ref{-2-2}) for \begin{eqnarray}\label {0-2-1}\int_D |\Delta_g u|^2
dR.\end{eqnarray} }

\vskip 0.2 true cm

\noindent  {\bf Proof.} \  \  Let ${\tilde \varrho}\in L^\infty
(\partial D)$ be the characteristic functions of $\Gamma_0$
(respectively, $\Gamma-\Gamma_1$) in the case of $K^d(D)$
(respectively, $K(D)$). Then $\int_\Gamma {\tilde \varrho}^3 u\, ds
=0$ for all $u$ in $\mathfrak{M}$, and hence according to
(\ref{0-2-1}), the condition (\ref{-2-2}) is satisfies such that
Lemma 2.3 can be applied. \ \ $\square$

\vskip 0.2 true cm

 Let $\bar D_1,\bar D_2\subset (M, g)$ be
 Lipschitz images of cubes with boundaries
$\partial D_1$ and $\partial D_2$.
 We say an open domain
$D$ with boundary $\partial D$ and closure $\bar D= D\cup \partial D$ is composed of
$\bar D_1$ and $\bar D_2$ if

a) \ \  $\partial D_1 \cap \partial D_2$ has positive measure

b) \ \  $\bar D=\bar D_1 \cup \bar D_2$

c) \ \  $\partial D\subset \partial D_1\cup \partial D_2$

d) \ \  $\partial D$  has positive measure.

The requirement d) excludes for instance the possibility that $M$ is
a sphere and $D=M$ the union of two hemispheres  $\bar D_1$ and
$\bar D_2$ (see, p.27 of \cite{Sa}).

By a finite number of domains, each of which is a Lipschitz image of
a cube, we can obtain more domains according to the above method.
Denoted by $\mathcal{F}$ all such domains. Completely similar to the
proofs of Lemmas 2.2, 2.3, we find that the compact trace Lemmas
2.2, 2.3 are also true for each domain in class $\mathcal{F}$.

\vskip 1.39 true cm

\section{Some completely continuous transformations and their eigenvalues}

\vskip 0.45 true cm

  Let $(\mathcal{M}, g)$ be an $n$-dimensional real analytic Riemannian
 manifold and let $D\subset M$ be a bounded
 domain with boundary $\Gamma$. Suppose that $D$ is of the type defined in Section 2 (i.e.,
 $D\in \mathcal{F}$) so that
 the compact trace lemmas 2.2, 2.3 are true. Let $\varrho$
 be a non-negative bounded function defined on
  $\Gamma$ or only on a portion $\Gamma_\varrho$ of $\Gamma$
  (measure $\Gamma_\varrho=\int_{\Gamma_\varrho}
  ds >0$) and assume that
  $\int_{\Gamma_\varrho}\varrho^3 ds>0$.
  In the case $\Gamma_\varrho \ne \Gamma$ we let
   $\Gamma_0$ be a  $C^{2,\varepsilon}$-smooth $(n-1)$-dimensional surface in $\Gamma-\bar\Gamma_{\varrho}$.

If $\Gamma_\varrho \ne \Gamma$ (measure $\Gamma-\Gamma_\varrho >0$),
we denote
\begin{eqnarray*} &K(D)=\{u\big|u\in Lip(\bar D)\cap H^2(D),\;\;
\frac{\partial u}{\partial \nu}=0 \;\; \mbox{on}\;\;
\Gamma_\varrho \cup \Gamma_0, \, \, \mbox{and}\;\; u=0\;\;\mbox{on}\;\; \Gamma-\Gamma_\varrho\}, \\
 &K^d(D)=\{u\big|u\in Lip(\bar D)\cap H^2(D),\,\, \frac{\partial u}{\partial \nu}=0\;\;
 \mbox{on}\;\;  \Gamma_\varrho, \,\,
 \mbox{and}\,\, u=\frac{\partial u}{\partial \nu}=0 \;\; \mbox{on}\;\; \Gamma_0\}.\end{eqnarray*}

If $\Gamma_\varrho =\Gamma$, we denote
\begin{eqnarray*} N(D)=\{u\big|u\in Lip (\bar D)\cap  H^2(D),\;
\frac{\partial u}{\partial \nu}=0 \;\;\mbox{on}\;\; \Gamma,
\;\int_\Gamma \varrho^3 uu_{01}\, ds =0 \;\, \mbox{and}\;\,
\int_\Gamma \varrho^3 uu_{{}_{02}}\, ds=0\},\end{eqnarray*} where
$u_{01}=1$ and $u_{02}$ are two eigenfunctions corresponding to the
Steklov eigenvalue $\lambda=0$ (see, Section 1).
 We shall also use the notation
\begin{eqnarray*}    \langle u,v\rangle^{\star} =\int_D (\triangle_g u)(\triangle_g v) dR,
 \quad \, u, v \in K(D) \;\; \mbox{or}\;\; L^d(D) \;\;\mbox{or}\;\;
 N(D).\end{eqnarray*}
 The bilinear functional $\langle u,v\rangle^{\star}$ can be used
as an inner product in each of the spaces $K(D)$, $K^d(D)$ and
$N(D)$.
 In fact, $\langle u,v\rangle^{\star}$ is a positive, symmetric, bilinear functional.
  In addition, if $\langle u,u\rangle^{\star}=0$, then $\triangle_g u=0$ in $D$.
  For $u\in K(D)$ or $N(D)$, by applying
  Green's formula, we have  \begin{eqnarray*} 0=\int_D u(\triangle_g u) dR =-\int_D |\nabla_g u|^2 dR
  -\int_{\partial D} u\, \frac{\partial u}{\partial \nu}\, ds
  = -\int_D |\nabla_g u|^2 dR, \end{eqnarray*}
 which implies $u\equiv constant$ in $D$.
  In the case $u\in K(D)$, in view of  $u=0$ on $\Gamma-\Gamma_\varrho$, we get that $u\equiv 0$ in $D$;
In the case $u\in N(D)$, since $\int_D \varrho^3 u\, ds =0$ with $
\int_{\Gamma}\varrho^3 \, ds >0$, we obtain $u\equiv 0$ in $D$;
 In the
case  $u\in K^d (D)$, since $\Gamma_0$ is a
$C^{2,\varepsilon}$-smooth surface and  $u=\frac{\partial
 u}{\partial \nu}=0$ on $\Gamma_0$, we find by
 $\Delta u=0$ in $D$ and
 Holmgren's uniqueness theorem
 (see, Corollary 5 of p.$\,$39 in \cite{Ra})  that
  $u\equiv 0$ in $D$.
  Closing $K(D)$, $K^d(D)$ and $N(D)$ with respect to  the norm
  $\|u\|^\star=\sqrt{\langle u, u\rangle^{\star}}$,
  we get the Hilbert spaces $(\mathcal{K}, \|\cdot\|^\star)$,
   $({\mathcal{K}}^d, \|\cdot\|^\star)$
  and $(\mathcal{N}, \|\cdot\|^\star)$, respectively.

  Next, we consider two linear
functionals
\begin{eqnarray*} [u,v]=\int_{\Gamma_\varrho} \varrho^3u v  \, ds
\end{eqnarray*}
and \begin{eqnarray} \label{3.0}\langle u,v\rangle
 =\langle u,v\rangle^{\star}  +[u,v], \end{eqnarray} where $u,v\in K(D)$ or
 $u,v\in K^d(D)$ or $u,v\in N(D)$.  It is clear that $\langle
u,v\rangle$ is an inner product in each of the spaces $K(D)$, $K^d
(D)$ and $N(D)$.

\vskip 0.25 true cm

 \noindent  {\bf Lemma 3.1.} \ \  {\it The norm
 \begin{eqnarray*}  \|u\|^\star=\sqrt{\langle u,u\rangle^{\star} }\end{eqnarray*}
 and \begin{eqnarray*} \|u\|=\sqrt{\langle u,u\rangle}\end{eqnarray*}
are equivalent in $K(D)$, $K^d(D)$ and $N(D)$.}

\vskip 0.2 true cm

 \noindent  {\bf Proof.} \ \ We have to show that there exist
 positive constants $C_1$ and $C_2$  such that
 $$ C_1\|u\|^\star\le \|u\| \le C_2\|u\|^\star\quad \, \mbox{for all}\;\; u\;\;
 \mbox{in}\;\; K(D)\;\; \mbox{or}\;\; K^d(D)\;\;\mbox{or}\;\;
 N(D).$$
 Obviously, $\|u\|^\star\le \|u\|$ for all $u$
 in each of three spaces.
  Let us first consider the case $u\in N(D)$.
It suffices to show that $\|u\|$ is bounded when $u$ belongs to the
set
\begin{eqnarray*} \mathfrak{M}=\{u\big|u\in N(D), \;  \|u\|^\star\le 1\}.\end{eqnarray*}
   It follows from Lemma 2.2 that
${\mathfrak{M}}_\Gamma:=\{u\big|_{\Gamma}:\, u\in \mathfrak{M}\}$ is
precompact in $L^2(\Gamma)$. This implies that there exists a
constant $C>0$ such that $\int_\Gamma  u^2 ds \le C$ for all $u\in
\mathfrak{M}$. Therefore, $[u,u]=\int_{\Gamma} \varrho^3 u^2 ds$ is
bounded in $\mathfrak{M}$, and so is $\|u\|^2=\langle
u,u\rangle^{\star} +[u,u]$.

Next, we consider the case for $u\in K (D)$ and still denote
\begin{eqnarray*} \mathfrak{M}=\{u\big|u\in K(D), \;  \|u\|^\star\le
1\}.\end{eqnarray*}
 By taking
 \begin{eqnarray*} \tilde \varrho =\left\{\begin{array}{ll} 1 \quad \;
 \mbox{for}\;\; x\in \Gamma-\Gamma_\varrho,\\
 0\quad \; \mbox{for}\;\; x\in
 \Gamma_\varrho,\end{array}\right.\end{eqnarray*}
  we have \begin{eqnarray} \int_D |\triangle u|^2 dR + \big(\int_\Gamma {\tilde \varrho}^3 u\, ds
\big)^2\le  C\;\; \mbox{for every}\;\; u\in \mathfrak{M}.
\end{eqnarray} It follows from
 Lemma 2.3 that
  ${\mathfrak{M}}_\Gamma$ is precompact in $L^2(\Gamma)$.
  In particular, $\int_{\Gamma_\varrho} u^2 \, ds$ is bounded on
  $\mathfrak{M}$ and hence also $[u,u]=\int_{\Gamma_\varrho} \varrho^3 u^2 \, ds
  $ and $\|u\|=\|u\|^\star +[u,u]$.

 Similarly, by taking $\tilde \varrho $ to be the characteristic function of $\Gamma_0$
 and by applying Lemma 2.3,  we can prove the
 corresponding results for the space $K^d(D)$.
  \ \  $\square$

\vskip 0.22 true cm

From Lemma 3.1, it follows that
\begin{eqnarray*}  \qquad  |[u,u]|=\bigg|\int_{\Gamma_\varrho}
\varrho^3
  u^2  ds \bigg|\le C \langle u,u\rangle^{\star} \; \; \mbox{for
  all}\;\; u\;\; \mbox{in}\;\; K(D) \;\;\mbox{or}\;\; K^d(D) \;\;
\mbox{or}\;\; N(D).\end{eqnarray*} Therefore, $[u,v]$ is a bounded,
symmetric,
  bilinear functional in
$(K(D)$,$\langle\cdot$,$\cdot\rangle^{\star})$,
$(K^d(D)$,$\langle\cdot$,$\cdot\rangle^{\star})$ and $(N(D)$,
$\langle \cdot$,$\cdot\rangle^{\star})$. Since it is densely defined
in $(\mathcal{K}$,$\langle\cdot,\cdot\rangle^{\star})$,
$({\mathcal{K}}^d$, $\langle\cdot,\cdot\rangle^{\star})$ and
$\,(\mathcal{N}, \langle\cdot,\cdot\rangle^{\star})$, respectively,\
it can immediately be extended to
$\,(\mathcal{K}$,$\langle\cdot,\cdot\rangle^{\star})$,
$({\mathcal{K}}^d$,$\langle\cdot,\cdot\rangle^{\star})$ and
$(\mathcal{N}, \langle\cdot,\cdot\rangle^{\star})$. We still use
$[u,v]$ to express the extended functional. Then there is a bounded
linear transformation $G_{\mathcal{K}}^{(\star)}$ of
$(\mathcal{K},\langle\cdot,\cdot\rangle^{\star})$ into
$(\mathcal{K},\langle\cdot,\cdot\rangle^{\star})$
 (respectively, $G_{{\mathcal{K}}^d}^{(\star)}$ of
 $({\mathcal{K}}^d,\langle\cdot,\cdot\rangle^{\star})$
into $({\mathcal{K}}^d,\langle\cdot,\cdot\rangle^{\star})$,
    $G_{\mathcal{N}}^{(\star)}$ of
$(\mathcal{N},\langle\cdot,\cdot\rangle^{\star})$ into
$(\mathcal{N},\langle\cdot,\cdot\rangle^{\star})$) such that
   \begin{eqnarray} \label {3-1} [u,v]=\langle G_{\mathcal{K}}^{(\star)}
    u, v\rangle^{\star} \quad \,\;
   \mbox{for all}\;\; u \;\; \mbox{and}\;\; v\;\;\mbox{in}\;\;
   \mathcal{K} \end{eqnarray}
  (respectively, \begin{eqnarray} \label {3.1} [u,v]=
  \langle G_{{\mathcal{K}}^d}^{(\star)} u, v\rangle^{\star} \quad \,\;
   \mbox{for all}\;\; u \;\; \mbox{and}\;\; v\;\;\mbox{in}\;\;
   {\mathcal{K}}^d, \end{eqnarray}
    \begin{eqnarray} \label{3.2} [u,v]= \langle G_{\mathcal{N}}^{(\star)} u, v\rangle^{\star} \quad \;\,
   \mbox{for all}\;\; u \;\; \mbox{and}\;\; v\;\;\mbox{in}\;\;
   \mathcal{N}).\end{eqnarray}

 \vskip 0.25 true cm

 \noindent  {\bf Lemma 3.2.} \ \  {\it The transformations
 $G_{\mathcal{K}}^{(\star)}$, $G_{{\mathcal{K}}^d}^{(\star)}$
 and $G_{\mathcal{N}}^{(\star)}$ are self-adjoint and compact.}

\vskip 0.2 true cm

 \noindent  {\bf Proof.} \  Since $[u,v]$ is symmetric, we immediately get
 that the transformation $G_{\mathcal{K}}^{(\star)}$, $G_{{\mathcal{K}}^d}^{(\star)}$ and
$G_{\mathcal{N}}^{(\star)}$ are all self-adjoint. For the
compactness,
 we only discuss the case for the transformation
$G_{\mathcal{K}}^{(\star)}$. It suffices to
  show (see, p.$\,$204 of \cite{RN}): From every sequence $\{u_m\}$ in $K(D)$ which
is bounded
\begin{eqnarray} \label{3-4} \|u_m\|^\star \le constant, \; \, m=1,2, 3, \cdots, \end{eqnarray}
   we can pick out a subsequence $\{u_{m'}\}$ such that
   \begin{eqnarray} \label{3-5} \langle G_{\mathcal{K}}^{(\star)} (u_{m'}-u_{l'}),
   (u_{m'}-u_{l'})\rangle^{\star}
   \to 0 \quad \, \mbox{when}\;\; m', l'\to \infty.\end{eqnarray}

   Applying Lemmas 2.3, 3.1 with the aid of (\ref{3-4}), we find that the
   sequence $\{u_m\big|_{\Gamma_\varrho}\}$ is
   precompact in $L^2(\Gamma_\varrho)$, so that
  there is a subsequence $\{u_{m'}\}$ such that
  \begin{eqnarray*} \int_{\Gamma_\varrho} \left(u_{m'}
  -u_{l'}\right)^2 ds \to 0 \quad \,
  \mbox{as} \;\;
    m', l'\to \infty. \end{eqnarray*}
   Therefore
 \begin{eqnarray*} [u_{m'}- u_{l'}, u_{m'} -u_{l'}] =
 \int_{\Gamma_\varrho} \varrho^3 \left(u_{m'}
  -u_{l'}\right)^2 ds \to 0 \quad \,
  \mbox{as} \;\;
    m', l'\to \infty, \end{eqnarray*}
which implies (\ref{3-5}). This proves the compactness of
$G_{\mathcal{K}}^{(\star)}$.
   \ \  $\square$

\vskip  0.25 true cm

Except for the transformations $G^{(\star)}_{\mathcal{K}}$,
$G^{(\star)}_{{\mathcal {K}}^d}$ and $G^{(\star)}_{\mathcal{N}}$,
 we need introduce corresponding transformations
 $G_{\mathcal{K}}$,
  $G_{{\mathcal{K}}^d}$
   and  $G_{\mathcal{N}}$
  by the inner product $\langle \cdot, \cdot \rangle$.
  Since
\begin{eqnarray} 0\le [u,v]\le \langle u,v\rangle \quad \;\; \mbox{for all } u
\;\;\mbox{in}\;\; K(D)\;\;\mbox{or}\;\; K^d(D) \;\; \mbox{or}\;\;
 N(D),\end{eqnarray}  there is a bounded linear self-adjoint
transformation $G_{\mathcal{K}}$ of $(\mathcal{K}, \langle \cdot,
\cdot\rangle)$
 (respectively,  $G_{{\mathcal{K}}^d}$ of $({\mathcal{K}}^d, \langle \cdot,
\cdot\rangle)$,  $\,G_{\mathcal{N}}$ of $(\mathcal{N}, \langle
\cdot, \cdot \rangle)$)
 such that
\begin{eqnarray} \label{3-6} [u,v]=\langle G_{\mathcal{K}} u,v\rangle
\quad \;\; \mbox{for all} \;\; u\;\; \mbox{and}\;\; v\;\;
\mbox{in}\;\; \mathcal{K}\end{eqnarray}
 (respectively,
\begin{eqnarray} \label{3-7} [u,v]=\langle G_{{\mathcal{K}}^d}u,v\rangle
\quad \;\; \mbox{for all} \;\; u\;\; \mbox{and}\;\; v\;\;
\mbox{in}\;\; {\mathcal{K}}^d, \end{eqnarray}   \begin{eqnarray}
\label{3-8} [u,v]=\langle G_{\mathcal{N}}u,v\rangle
 \quad \;\; \mbox{for all} \;\; u\;\;
\mbox{and}\;\; v\;\; \mbox{in}\;\; \mathcal{N}).\end{eqnarray}

\vskip 0.25 true cm

 \noindent  {\bf Lemma 3.3.} \ \  {\it The transformations
 $G_{\mathcal{K}}$, $G_{{\mathcal{K}}^d}$
 and  $G_{\mathcal{N}}$ are positive and compact.}

\vskip 0.2 true cm

\noindent  {\bf Proof.} \ From $[u,u]\ge 0$  for any $u\in
{\mathcal{K}} $ or ${\mathcal{K}}^d$ or
 ${\mathcal{N}}$, we immediately know that
 $G_{\mathcal{K}}$, $G_{{\mathcal{K}}^d}$
 and  $G_{\mathcal{N}}$ are positive. The proof of the compactness is completely
 similar to that of Lemma 3.2. \ \   $\square$

 It follows from Lemma 3.3 that
 $G_{\mathcal{K}}$ (respectively, $G_{{\mathcal{K}}^d}$,  $G_{\mathcal{N}}$)
  has only non-negative eigenvalues and that the positive
 eigenvalues form an enumerable sequence  $\{\mu_{\mathcal{K}}\}$
 (respectively, $\{\mu_{{\mathcal{K}^d}}\}$,
   $\{\mu_{\mathcal{N}}\}$) with $0$ as the
 only limit point.

\vskip 0.25 true cm

 \noindent  {\bf Theorem 3.4.} \ \  {\it The transformations
 $G_{\mathcal{K}}^{(\star)}$ and
  $G_{\mathcal{K}}$ (respectively,  $G_{{\mathcal{K}}^d}^{(\star)}$ and
  $G_{{\mathcal{K}}^d}$,  $\,G_{\mathcal{N}}^{(\star)}$ and
  $G_{\mathcal{N}}$) have the same eigenfunctions. If $\mu_{\mathcal{K}}^{\star}$
 and $\mu_{\mathcal{K}}$
  (respectively,
  $\mu_{{\mathcal{K}}^d}^{\star}$
 and $\mu_{{\mathcal{K}}^d}$,
    $\, \mu_{\mathcal{N}}^{\star}$
 and $\mu_{\mathcal{N}}$)
 are eigenvalues corresponding to
 the same eigenfunction we have
 \begin{eqnarray} \label{3.10}\mu_{\mathcal{K}}=
 \frac{\mu_{\mathcal{K}}^{\star}}{1+\mu_{\mathcal{K}}^{\star}}\end{eqnarray}
(respectively,   \begin{eqnarray} \label{3.11}
\mu_{{\mathcal{K}}^d}=
 \frac{\mu_{{\mathcal{K}}^d}^{\star}}{1+\mu_{{\mathcal{K}}^d}^{\star}}, \end{eqnarray}
   \begin{eqnarray} \label {3.12}\mu_{\mathcal{N}}=
 \frac{\mu_{\mathcal{N}}^{\star}}{1+\mu_{\mathcal{N}}^{\star}}).\end{eqnarray}
 }

\vskip 0.2 true cm

 \noindent  {\bf Proof.} \  We only prove the case
 for the $G_{\mathcal{K}}$
 (the arguments are similar for  $G_{{\mathcal{K}}^d}$
 and $G_{\mathcal{N}}$). Since
  $G_{\mathcal{K}}^{(\star)}$ is positive, we can easily conclude that the inverse
 $(1+G_{\mathcal{K}}^{(\star)})^{-1}$ exists and is a bounded self-adjoint
 transformation.
   By virtue of  (\ref{3-1}), (\ref{3-6}) and (\ref{3.0}), we have
\begin{eqnarray} \label{3-9}  \langle G_{\mathcal{K}}^{(\star)} u, v\rangle^{\star}
 &=&[u,v] =  \langle G_{\mathcal{K}}u, v\rangle \\
 &=&  \langle G_{\mathcal{K}}u, v\rangle^{\star}
   + \langle G_{\mathcal {K}}^{(\star)} G_{\mathcal{K}} u, v\rangle^{\star}\nonumber\\
   &=&  \langle G_{\mathcal{K}}u, v\rangle^{\star}
   + \langle G_{\mathcal {K}} G_{\mathcal{K}}^{(\star)} u, v\rangle^{\star},
    \quad \;
 (u,v\in \mathcal{K}).\nonumber \end{eqnarray}
 It follows that
 \begin{eqnarray} \label {3-10} G_{\mathcal{K}}=
  G_{\mathcal{K}}^{(\star)}(1+G_{\mathcal{K}}^{(\star)})^{-1},\end{eqnarray}
   from which the desired result follows immediately.
\ $\square$

\vskip 0.25 true cm

 \noindent  {\bf Proposition 3.5.} \   {\it Let $u$ and $v$ be two eigenfunctions in
 $({\mathcal{K}},\langle \cdot, \cdot\rangle)$ (respectively,
 $({\mathcal{K}}^d,\langle \cdot, \cdot\rangle)$,
   $({\mathcal{N}},\langle \cdot, \cdot\rangle)$)
 of the transformation $G_{\mathcal{K}}$ (respectively,
 $G_{{\mathcal{K}}^d}$,
  $G_{\mathcal{N}}$) at least one of which corresponds to a non-vanishing eigenvalue.
 Then $u$ and $v$ are orthogonal if and only if the
 $u\big|_{\Gamma_\varrho}$
 and $v\big|_{\Gamma_\varrho}$
 are orthogonal in $L_\varrho^2 (\Gamma_\varrho)$, that is,
 \begin{eqnarray} [u,v]=\int_{\Gamma_\varrho} \varrho^3
  u  v\,ds=0.\end{eqnarray}}

 \vskip 0.15 true cm

 \noindent  {\bf Proof.} \ \  Without loss of generality,
 we suppose that $u$ is the eigenfunction corresponding
 to the eigenvalue $\mu \ne 0$. Then
 \begin{eqnarray*} [u,v] = \langle G_{\mathcal{K}}u, v\rangle
 =\mu
 \langle u, v\rangle, \end{eqnarray*}
   which implies the desired result.
    $\;\; \square$

\vskip 0.25 true cm

We can now prove  \vskip 0.20 true cm

 \noindent  {\bf Theorem 3.6.} \ \  {\it  Let  $D\subset (\mathcal{M},g)$
  be a bounded domain with piecewise smooth boundary $\Gamma$, and let $D\in \mathcal{F}$.
   If $u$ is an eigenfunction of
 the transformations $G_{\mathcal{K}}^{(\star)}$
  or $G_{\mathcal{N}}^{(\star)}$ with eigenvalue $\mu^\star \ne 0$,
  then
  $u$ has derivatives of any order in $D$ and is such that
  \begin{eqnarray}  \label{3-15} \left\{\begin{array}
  {ll}\triangle_g^2 u=0 \quad \mbox{in}\;\; D,  \\
   \frac{\partial u}{\partial \nu}=0 \; \; \mbox{on}\;\;
   \Gamma_\varrho\cup \Gamma_0,\quad\;
  u=0 \;\; \mbox{on}\;\; \Gamma-\Gamma_\varrho,\\
  \Delta_g u=0\;\;\mbox{on}\;\; \Gamma-(\Gamma_\varrho\cup
  \Gamma_0),\\
    \frac{\partial (\triangle_g u)}{\partial \nu}- \lambda^3  \varrho^3
    u=0 \;\; \mbox{on}\;\; \Gamma_\varrho,
 \; \quad \; \mbox{with}\;\; \lambda^3=\frac{1}{\mu^\star}. \end{array}\right.\end{eqnarray}}

\vskip 0.2 true cm

 \noindent  {\bf Proof.} \ \
  Since $u$ is an eigenfunction of $G_{\mathcal{K}}^{(\star)}$
  (i.e., $G_{\mathcal{K}}^{(\star)} u=
\mu^{\star} u$), we have that $\frac{\partial u}{\partial
\nu}\big|_{\Gamma_\varrho\cup \Gamma_0}=0$ and
$u\big|_{\Gamma-\Gamma_\varrho}=0$, and that
 \begin{eqnarray*}  \label{03-03}\int_{\Gamma_\varrho} \varrho^3 u
  v \, ds= \mu^{\star} \int_D
 (\triangle_g u)(\triangle_g v) dR \quad \, \mbox{for all} \;\;
 v\in {\mathcal{K}} (D). \end{eqnarray*}
 Applying Green's formula (see, p.$\,$114-120 of \cite{LM}, \cite{Ch1})
 to the right-hand side of the above equation,
 we obtain that
\begin{eqnarray*} \frac{1}{\mu^{\star}} \int_{\Gamma_\varrho} \varrho^3  u
  v \, ds=   \int_D
 (\triangle_g^2 u) v\, dR - \int_{\Gamma} (\triangle_g u)\frac{\partial v}{\partial \nu}\,ds
 +\int_{\Gamma} \frac{\partial (\triangle_g u)}{\partial \nu} \,v \, ds \end{eqnarray*}
 for all $v\in {\mathcal {K}}(D)$,
  where $\frac{\partial (\triangle_g u)}{\partial \nu}\in
  H^{-\frac{3}{2}}(\Gamma)$ (see \cite{LM}).  From $\,\frac{\partial
 v}{\partial \nu}
  \big|_{\Gamma_\varrho\cup \Gamma_0}=0$
   and $v\big|_{\Gamma-\Gamma_\varrho}=0$, we get
 \begin{eqnarray} \quad \quad \quad  \;\int_D  (\triangle_g^2 u)v\, dR -
 \int_{\Gamma-(\Gamma_\varrho\cup\Gamma_0)} (\triangle_g u)
 \frac{\partial v}{\partial \nu} \, ds+
\int_{\Gamma_\varrho}  \left(\frac{\partial (\triangle_g
u)}{\partial \nu} - \frac{1}{\mu^*}\,\varrho^3 u\right)v \, ds
  =0\end{eqnarray}
  for all $v\in {\mathcal K}(D)$.
   By taking all $v\in C_0^\infty (D)$, we have $\Delta^2_g u=0$ in
   $D$. It follows from the interior regularity of elliptic equations that $u\in C^\infty(D)$.
    Noticing  that
$v\big|_{\Gamma_{\varrho}}$ and $\frac{\partial v}{\partial
\nu}\big|_{\Gamma-(\Gamma_\varrho\cup \Gamma_0)}$
 run throughout
   space $L^2(\Gamma_\varrho)$  and
   $L^2(\Gamma-(\Gamma_\varrho\cup\Gamma_0))$, respectively,
    when $v$ runs throughout space $K(D)$,
   we see
 that \begin{eqnarray*}  \Delta_g u=0\,\, \mbox{on}\;\; \Gamma-(\Gamma_\varrho\cup \Gamma_0),\quad \,
 \mbox{and}\;\; \frac{\partial (\triangle_g u)}{\partial \nu}  -\frac{1}{\mu^{\star}}
 \varrho^3
   u  =0 \;\; \mbox{on}\;\; \Gamma_\varrho.\end{eqnarray*}
 Therefore, (\ref{3-15}) holds.  In a similar way, we can prove
 the desired result for $G_{\mathcal{N}}$.
 \ \ $\square$

\vskip 0.24 true cm

  \noindent  {\bf Theorem 3.7.} \ \  {\it  Let $(\mathcal{M},g)$ be a real
 analytic  Riemannian manifold, and let $D\subset (\mathcal{M},g)$
  be a bounded domain with piecewise smooth boundary $\Gamma$. Let
 $D\in \mathcal{F}$.
 Assume that $\Gamma_0$ is a $C^{2,\varepsilon}$-smooth $(n-1)$-dimensional
 surface in $\Gamma-\bar \Gamma_\varrho$.
   If $u$ is an eigenfunction of
 the transformations $G_{{\mathcal{K}}^d}^{(\star)}$
  with eigenvalue $\mu^\star\ne 0$,
  then $u$ has derivatives of any order in $D$ and is such that
  \begin{eqnarray}  \label{3-30} \left\{\begin{array} {ll} \triangle_g^2 u=0 \quad \mbox{in}\;\; D,\\
             \frac{\partial u}{\partial \nu}=0 \; \; \mbox{on}\;\; \Gamma_\varrho,\\
      u=\frac{\partial u}{\partial \nu}=0 \quad \mbox{on}\;\;
      \Gamma_0, \\
     \Delta_g u =0, \;\; \frac{\partial (\Delta_g u)}{\partial \nu}=0
       \;\;\mbox{on}\;\; \Gamma-(\Gamma_\varrho\cup \Gamma_0),\\
   \frac{\partial (\triangle_g u)}{\partial \nu}- \lambda^3 \, \varrho^3 \,
   u=0 \;\; \mbox{on}\;\; \Gamma_{\varrho}, \quad \; \mbox{with}\;\;
   \lambda^3=\frac{1}{\mu^\star}. \end{array}\right.\end{eqnarray}}

\vskip 0.2 true cm

 \noindent  {\bf Proof.} \
 If $G_{{\mathcal{K}}^d}^{(\star)} u= \mu^\star
u$, then we have that $\frac{\partial u}{\partial \nu}=0$ on
$\Gamma_\varrho$ and $u=\frac{\partial u}{\partial \nu}=0$ on
$\Gamma_0$, and that
 \begin{eqnarray} \label{3/25}
 \int_{\Gamma_\varrho} \varrho^3 u v\, ds =
   \mu^\star \int_D (\triangle_g u)(\triangle_g v) dR \quad \; \mbox{for all}\;\;
 v\in {\mathcal{K}}^d (D),\end{eqnarray}
  By using Green's formula and noticing that
   $\frac{\partial v}{\partial \nu} \big|_{\Gamma_\varrho}=0$ and $v\big|_{\Gamma_0}=
  \frac{\partial v}{\partial \nu}\big|_{\Gamma_0}=0$, we get
   that \begin{eqnarray} \label{3..27}
  &&  \int_D  (\triangle_g^2 u) v\, dR+
 \int_{\Gamma_\varrho} \left(\frac{\partial (\triangle_g u)}{\partial \nu}- \frac{1}{\mu^\star}
 \varrho^3
   u \right)v\, ds +\int_{\Gamma-(\Gamma_\varrho\cup \Gamma_0)} \frac{\partial (\Delta_g u)}{\partial \nu}\,
   v\, ds\\
  && \quad \;\; -
\int_{\Gamma-(\Gamma_\varrho\cup\Gamma_0)} (\triangle_g u)
\frac{\partial v}{\partial \nu} \, ds =0 \quad \; \mbox{for all}\;\;
v\in {K}^d (D),\nonumber
  \end{eqnarray}
where $\frac{\partial (\triangle_g u)}{\partial \nu}\in
H^{-\frac{3}{2}}(\Gamma-(\Gamma_\varepsilon\cup \Gamma_0))$.
  By taking all $v\in C^\infty_0(D)$, we obtain that $\triangle_g^2 u=0$ in
  $D$.
  Note that
   $\frac{\partial v}{\partial \nu}\big|_{\Gamma -(\Gamma_\varrho\cup \Gamma_0)}$
  and $v\big|_{\Gamma -\Gamma_0}$
   run throughout
  the spaces  $L^2(\Gamma -(\Gamma_{\varrho}\cup \Gamma_0))$
  and $L^2(\Gamma -\Gamma_0)$, respectively,
  when $v$ runs throughout the space $K^d(D)$.
    Thus we have
\begin{eqnarray*} \label{3..28}
    && \triangle_g u=0 \;\; \mbox{and}\;\; \frac{\partial (\Delta_g u)}{\partial \nu}=0\;\;
     \mbox{on}\;\;
  \Gamma -(\Gamma_\varrho\cup \Gamma_0),\\
   &&  \quad  \quad \,\mbox{and}\;\;
   \frac{\partial (\triangle_g u)}{\partial \nu}
   -\frac{1}{\mu^\star} \varrho^3
    u =0 \quad \,\mbox{on}\;\; \Gamma_\varrho. \quad \qquad \qquad \qquad \qquad
    \square
\end{eqnarray*}

\vskip 0.25 true cm

  \noindent  {\bf Theorem 3.8.} \ \  {\it  Let $(\mathcal{M},g)$, $D$, $\Gamma_\varrho$  and $\Gamma_0$ be as
  in Theorem 3.7.
  Assume that $\gamma_k^3$ and $\kappa_k^3$ are the $k$-th Steklov eigenvalues of the following problems:
  \begin{eqnarray}  \label{3-30-0} \left\{\begin{array} {ll} \triangle_g^2 u=0 \quad \mbox{in}\;\; D,\\
      \frac{\partial u}{\partial \nu}=0 \; \; \mbox{on}\;\; \Gamma_\varrho,\\
      u=\frac{\partial u}{\partial \nu}=0 \quad \mbox{on}\;\;
      \Gamma_0, \\
   \frac{\partial u}{\partial \nu}=0 \;\; \mbox{and}\;\;
   \frac{\partial (\triangle_g u)}{\partial \nu}=0 \;\; \mbox{on}\;\;
  \Gamma -(\Gamma_\varrho \cup \Gamma_0),\\
   \frac{\partial (\triangle_g u)}{\partial \nu} - \varsigma^3  \varrho^3
  u=0 \;\; \mbox{on}\;\; \Gamma_{\varrho} \end{array}\right.\end{eqnarray}
  and
  \begin{eqnarray*}  \label{3-31} \left\{\begin{array} {ll} \triangle_g^2 u=0 \quad \mbox{in}\;\; D,\\
      \frac{\partial u}{\partial \nu}=0 \; \; \mbox{on}\;\; \Gamma_\varrho,\\
      u=\frac{\partial u}{\partial \nu}=0 \quad \mbox{on}\;\;
      \Gamma_0, \\
  \Delta_g u=0\;\; \mbox{and}\;\; \frac{\partial (\triangle_g u)}{\partial \nu}=0
   \;\;  \mbox{on}\;\;
  \Gamma -(\Gamma_\varrho \cup \Gamma_0),\\
   \frac{\partial (\triangle_g u)}{\partial \nu}- \kappa^3 \varrho^3
  u=0 \;\; \mbox{on}\;\; \Gamma_{\varrho}, \end{array}\right.\end{eqnarray*}
   respectively.
  Then $\varsigma_k^3\le \kappa_k^3$ for all $k\ge 1$.}

\vskip 0.26 true cm

 \noindent  {\bf Proof.} \  \ For $0\le \alpha\le 1$, let $u_k=u_k(\alpha, x)$  be
 the normalized eigenfunction corresponding to the $k$-th Steklov eigenvalue $\lambda_k$ for the
 following problem:
   \begin{eqnarray} \label{3---31}  \left\{\begin{array} {ll} \triangle_g^2 u_k=0 \quad \mbox{in}\;\; D,\\
      \frac{\partial u_k}{\partial \nu}=0 \; \; \mbox{on}\;\; \Gamma_\varrho,\\
      u_k=\frac{\partial u_k}{\partial \nu}=0 \quad \mbox{on}\;\;
      \Gamma_{0}, \\
  \alpha \Delta_g u_k+(1-\alpha) \frac{\partial u_k}{\partial \nu}=0\;\; \mbox{and}\;\;
  \frac{\partial (\Delta_g u_k)}{\partial \nu}=0
   \;\;  \mbox{on}\;\;
  \Gamma -(\Gamma_\varrho \cup \Gamma_{0}),\\
   \frac{\partial (\triangle_g u_k)}{\partial \nu}- \lambda \, \varrho^3 \,
  u_k=0 \;\; \mbox{on}\;\; \Gamma_{\varrho}. \end{array}\right.\end{eqnarray}
   It is easy to verify (cf, p.$\;$410 or Theorem 9 of p.$\;$419 in \cite{CH})
   that the $k$-th Steklov eigenvalue
    $\lambda_k=\lambda_k(\alpha)$ is continuous on the closed interval $[0,1]$
   and differentiable in the open interval $(0,1)$, and that
   $u_k(\alpha, x)$ is also differentiable with respect to $\alpha$ in
   $(0,1)$ (cf. \cite{Fri}).
   \  We will denote by $'$ the derivative with respect to $\alpha$.
    Then
 \begin{eqnarray} \label{32} \, \qquad \;\;\quad \left\{\begin{array} {ll}
  \triangle_g^2 u'_k=0 \quad \mbox{in}\;\; D,\\
      \frac{\partial u'_k}{\partial \nu}=0 \; \; \mbox{on}\;\; \Gamma_\varrho,\\
      u'_k=\frac{\partial u'_k}{\partial \nu}=0 \quad \mbox{on}\;\;
      \Gamma_{0}, \\
   \Delta_g u_k+ \alpha \triangle_g u'_k -\frac{\partial u_k}{\partial
   \nu} +  (1-\alpha) \frac{\partial u'_k}{\partial \nu}=0\,\, \mbox{and}\,\,
  \frac{\partial (\Delta_g u'_k)}{\partial \nu}=0
   \,\,  \mbox{on}\,\,
  \Gamma -(\Gamma_\varrho \cup \Gamma_{0})\\
   \frac{\partial(\triangle_g u'_k)}{\partial \nu}- \lambda' \, \varrho^3 \,
 \frac{\partial u_k}{\partial \nu}-\lambda \varrho^3 \frac{\partial u'_k}{\partial \nu}=0
 \,\, \mbox{on}\,\, \Gamma_{\varrho}. \end{array}\right. \end{eqnarray}
  Multiplying  (\ref{32}) by $u_k$,  integrating the product over
  $D$, and then applying Green's formula,
 we get that for $0<\alpha<1$
\begin{eqnarray*} 0&=&\int_D (\Delta^2_g u'_k)u_k \,dR =
\int_D (\Delta_g^2 u_k)u'_k\, dR  - \int_{\partial D} (\Delta_g
  u_k)\frac{\partial u'_k}{\partial \nu}\, ds \\
  && +\int_{\partial D} u'_k
  \frac{\partial (\Delta_g u_k)}{\partial \nu}\, ds -\int_{\partial D}
  u_k\frac{\partial (\Delta_g u'_k)}{\partial \nu}\, ds +\int_{\partial
  D} (\Delta_g u'_k) \frac{\partial u_k}{\partial \nu}\,ds \nonumber \\
  &=& \left [ \int_{\Gamma_\varrho} u'_k \frac{\partial
  (\Delta u_k)}{\partial \nu} \, ds - \int_{\Gamma -(\Gamma_\varrho \cup
  \Gamma_{0})} (\Delta_g u_k)\frac{\partial
  u'_k}{\partial \nu} \, ds \right]\nonumber \\
  && +  \left [- \int_{\Gamma_\varrho} u_k\frac{\partial
  (\Delta u_k')}{\partial \nu} \, ds + \int_{\Gamma -(\Gamma_\varrho \cup
  \Gamma_{0})} (\Delta_g u'_k)\frac{\partial
  u_k}{\partial \nu} \, ds \right]\nonumber \\&=&
  \left [ \int_{\Gamma_\varrho}
  \lambda \varrho^3 u_k\,u'_k \, ds
   + \int_{\Gamma -(\Gamma_\varrho \cup
  \Gamma_{0})} \left(\frac{1-\alpha}{\alpha} \,\frac{\partial u_k}{\partial \nu}\right)\frac{\partial
  u'_k}{\partial \nu} \, ds \right]\nonumber \\
  && +   \int_{\Gamma_\varrho} \left(-\lambda' \varrho^3  u_k-\lambda \varrho^3  u'_k\,
\right)  u_k \, ds \nonumber \\
  && + \int_{\Gamma -(\Gamma_\varrho \cup
  \Gamma_{0})} \left( -\frac{1}{\alpha}\, \Delta_g u_k
  +\frac{1}{\alpha}\,
  \frac{\partial
  u_k}{\partial \nu} -\frac{1-\alpha}{\alpha} \, \frac{\partial u'_k}{\partial \nu} \right)
  \frac{\partial u_k}{\partial \nu} \, ds \nonumber \\
&=& -\lambda' \int_{\Gamma_\varrho} \varrho^3  u_k^2\, ds +
 \int_{\Gamma -(\Gamma_\varrho \cup
  \Gamma_{0})} \left[\left(\frac{1-\alpha}{\alpha^2}\right) \frac{\partial u_k}{\partial \nu}
  +\frac{1}{\alpha} \,\frac{\partial u_k}{\partial \nu}\right] \frac{\partial u_k}{\partial \nu}\,
  ds \nonumber \\
  &=& -\lambda' \int_{\Gamma_\varrho} \varrho^2 u_k^2 ds +
 \int_{\Gamma -(\Gamma_\varrho \cup
  \Gamma_{0})} \bigg(\frac{1}{\alpha} \, \frac{\partial u_k}{\partial
  \nu} \bigg)^2 ds,\nonumber \end{eqnarray*}
 i.e., $$\lambda'_k(\alpha)= \frac{\int_{\Gamma -(\Gamma_\varrho
\cup
  \Gamma_{0})} \left(\frac{1}{\alpha} \, \frac{\partial u_k}{\partial
  \nu} \right)^2 ds}{\int_{\Gamma_\varrho} \varrho^3 \,u_k^2 ds} >0 \quad \; \mbox{for all}\;\;
0<\alpha<1.$$
 This implies that $\lambda_k$ is increasing with respect to $\alpha$ in
 $(0,1)$.  Note that
 if we change the $\alpha$ from $0$ to $1$, each individual Steklov
 eigenvalue $\lambda_k$ increase monotonically form the value
 $\varsigma_k$ which is the $k$-th Steklov eigenvalue of (\ref{3-30}) to
the value $\kappa_k$ which is the $k$-th Steklov eigenvalue
 (\ref{3-31}).
 Thus, we have that $ \varsigma_k \le \kappa_k$  for all $k$.   \ \ $\square$

\vskip 0.26 true cm

 Conversely, we can show that a sufficiently smooth function satisfying
(\ref{3-15}) (respectively, (\ref{3-30})) is an eigenfunction of
$G_{\mathcal{K}}^{(\star)}$ or $G_{\mathcal{N}}^{(\star)}$
(respectively, $G_{{\mathcal{K}}^d}^{(\star)}$).

 \vskip 0.20 true cm

 \noindent  {\bf Proposition 3.9.} \ \  {\it  Let $\bar D$ be bounded domain with  piecewise smooth
  boundary, and let $D\in \mathcal{F}$. Assume that $u$ belongs to $C^4(\bar D)$ and let $\lambda> 0$.

a) \ \  If $\Gamma_{\varrho}\ne \Gamma$ and $u$ satisfies
(\ref{3-15}), then $u\in \mathcal{K}$ and $u$ is an  eigenfunction
of $G_{\mathcal{K}}^{(\star)}$ with the eigenvalue
$\mu^\star=\lambda^{-3}$,
\begin{eqnarray}  \label{3-27} G_{\mathcal{K}}^{(\star)}u=\lambda^{-3} u.\end{eqnarray}

b) \ \  If $\Gamma_{\varrho}\ne \Gamma$ and $u$ satisfies
(\ref{3-30}), then $u\in {\mathcal{K}}^d$ and $u$ is an
eigenfunction of $G_{{\mathcal{K}}^d}^{(\star)}$ with the eigenvalue
$\mu^\star=\lambda^{-3}$,
\begin{eqnarray}  \label{3-28} G_{{\mathcal{K}}^d}^{(\star)}u=\lambda^{-3} u.\end{eqnarray}

c) \ \  If $\Gamma_{\varrho}= \Gamma$ and $u$ satisfies
(\ref{3-15}), then $u\in \mathcal{N}$ and $u$ is an  eigenfunction
of $G_{\mathcal{N}}^{(\star)}$ with the eigenvalue
$\mu^\star=\lambda^{-3}$,
\begin{eqnarray} \label {3-29} G_{\mathcal{N}}^{(\star)}u=\lambda^{-3} u.\end{eqnarray}}

 \noindent  {\bf Proof.} \ \  i) \ \  $\Gamma_\varrho \ne \Gamma$. We
 claim that there is no eigenvalue $\lambda^3=0$. Suppose by contradiction that there is a  function
 $u$ in $C^4(\bar D)$ satisfying
 \begin{eqnarray} \label {3-36}\left\{\begin{array} {ll} \triangle_g^2 u=0 \;\;\mbox{in}\;\; D, \\
  \frac{\partial u}{\partial \nu}=0\;\;\mbox{on}\;\; \Gamma_\varrho\cup
  \Gamma_0,\quad \;
   u=0 \;\;\mbox{on}\;\;
  \Gamma-\Gamma_\varrho\\
  \Delta_g u=0\;\;\mbox{on}\;\; \Gamma-(\Gamma_\varrho\cup \Gamma_0),
    \quad \;\; \mbox{and}\;\; \frac{\partial (\triangle_g
  u)}{\partial\nu}=0 \;\; \mbox{on} \;\; \Gamma_\varrho.\end{array} \right.\end{eqnarray}
  Multiplying the above equation by $u$, integrating the result over
  $D$, and then using Green's formula, we derive
  \begin{eqnarray*}
  && 0=\int_D u(\triangle_g^2 u) dR= \int_D |\triangle_g u|^2 dR -  \int_\Gamma u
  \frac{\partial (\triangle_g u)}{\partial \nu} \,ds \\
  && \;\; \quad \; \; + \int_{\Gamma} (\triangle_g u)
 \frac{\partial u}{\partial \nu} ds=
 \int_D |\triangle_g u|^2 dR.  \end{eqnarray*}
This implies that $\triangle_g u=0$ in $D$, so that $$ 0=\int_D
u(\Delta_g u) dR =-\int_D |\nabla u|^2 dR -\int_{\partial D}
u\frac{\partial u}{\partial \nu} ds=-\int_D |\nabla u|^2 dR.$$ That
is, $u\equiv constant$ in $D$.  Since $u=0$ on
$\Gamma-\Gamma_\varrho$, we get that $u=0$ in $D$. The claim is
proved.

  In view of assumptions, we see that $u\in \mathcal{K}$.
  By (\ref{3-15}) and Green's formula, it follows that
  for an arbitrary $v\in K(D)$
  \begin{eqnarray*}  \langle G_{\mathcal{K}}^{(\star)}
   u, v\rangle^\star &=&[u,v] = \int_{\Gamma_\varrho}
  \varrho^3  u v\, ds\\
    &=& \lambda^{-3}\int_{\Gamma_\varrho}
    \frac{\partial (\triangle_g u)}{\partial \nu}\, v \, ds
  =   \lambda^{-3}
  \int_{\Gamma} \frac{\partial (\triangle_g u)}{\partial \nu} \,v\, ds \\
  &=&  \lambda^{-3}  \left[\int_{\Gamma}    (\triangle_g u)\,\frac{\partial  v}{\partial \nu}\, ds
   +\int_D (\triangle_g u) (\triangle_g v) dR - \int_D
  v(\triangle_g^2 u)dR\right]\\
  &=& \lambda^{-3} \int_D (\triangle_g
  u)(\triangle_g v) dR = \lambda^{-3}\langle u, v \rangle^\star, \end{eqnarray*}
  Therefore,  \begin{eqnarray*} \langle G_{\mathcal{K}}^{(\star)}u
     - \lambda^{-3}u, v \rangle^\star=0 \quad\quad \mbox{for all } v\in K(D), \end{eqnarray*}
 which  implies (\ref{3-27}).
 By a similar way, we can prove b).

   ii) \ \  $\Gamma_\varrho = \Gamma$. In this case, for the eigenvalue $\lambda^3=0$,
   the problem (\ref{3-15}) has the
   solutions $u_{01}=constant$ and $u_{02}(x)=\int_D F(x,y)dR_y$ in
   $D$, here $F(x,y)$ is Green's function with Neumann boundary condition (see, Section 1).
   These solutions do not belong to
   $N(D)$. If, however, $u$ is a solution with eigenvalue
   $\lambda^3>0$ then $u\in \mathcal{N}$.
  Indeed, by Green's formula we get
\begin{eqnarray*} \int_{\partial D} \frac{\partial (\triangle_g
u)}{\partial \nu} \,ds =\int_D \triangle_g^2 u \, dR=0
\end{eqnarray*}
and hence form (\ref{3-15}) we obtain
 \begin{eqnarray*} \int_{\partial D} \varrho^3 u\,u_{01} ds=\int_{\partial D} \varrho^3 u\, ds =0.\end{eqnarray*}
In addition, from (\ref{3-15}) we get
\begin{eqnarray*} \left\{ \begin{array}{ll} \triangle_g (\triangle_g u)=0\quad \mbox{in}\;\;
 D,\\
 \frac{\partial (\Delta_g u)}{\partial \nu}=\lambda^3 \varrho^3 u \quad \mbox{on}\;\;
 \partial D.\end{array}\right.\end{eqnarray*}
so that
 \begin{eqnarray*} \Delta_g u(x)= \int_{\partial D} F(x,y)
 \frac{\partial (\Delta_g u)}{\partial \nu_y} ds_y
 =\lambda^3\int_{\partial D} F(x,y) \varrho^3(y) u(y)\,
 ds_y.\end{eqnarray*}
  Combining this and Green's formula, we have
  \begin{eqnarray*} 0 &=& -\int_{\partial D} \frac{\partial
  u}{\partial \nu} ds =\int_D \Delta_g u\, dR =\lambda^3 \int_D
  \left( \int_{\partial D} F(x,y)\varrho^3(y)  u(y) ds_y\right) dR_x
  \\ &=& \lambda^3 \int_{\partial D} \varrho^3 (y) u(y) \left(\int_D
  F(x,y)dR_x\right) ds_y\\
&=& \lambda^3 \int_{\partial D} \varrho^3 (x) u(x) \left(\int_D
  F(x,y)dR_y\right) ds_x
 =\lambda^3 \int_{\partial D} \varrho^3 (x) u(x)\,u_{02}(x)\, ds_x, \end{eqnarray*}
  i.e.,
$\int_{\partial D} \varrho^3  u_{02} u \, ds =0,$ so that $u\in
\mathcal{N}$.
 Proceeding as in a), we can prove that (\ref{3-29}) holds. \ \ \ $\square$

\vskip 0.23 true cm

\noindent {\bf Remark 3.10.} \  Each of transformations
$G_{\mathcal{K}}^\star$, $G_{{\mathcal{K}}^d}^\star$ and
$G_{\mathcal{N}}^\star$ corresponds to a biharmonic Steklov problem
given by
  the quadratic forms
  \begin{eqnarray*}  \langle u, u\rangle^\star =\int_D
  |\triangle_g u|^2 dR \end{eqnarray*}
  and \begin{eqnarray*} [u,u] =\int_{\Gamma_\varrho}  \varrho^3
  u^2    ds \end{eqnarray*}
  and the function classes of ${\mathcal{K}}^\star$, ${{\mathcal{K}}^d}^\star$
  and ${\mathcal{N}}^\star$, respectively.
  The eigenvalues $\lambda_k^3$  of these  biharmonic Steklov problems are given by
  \begin{eqnarray} \label {3-37} \lambda_k^3 =1/\mu^\star_k, \quad k=1,2,3, \cdots.\end{eqnarray}
Since $0$ is the only limit point of $\mu_k^\star$, the only
possible limit points of $\lambda_k^3$ are $+\infty$.

 \vskip 1.39 true cm

\section{Biharmonic Steklov eigenvalues on a rectangular parallelepiped}

\vskip 0.45 true cm

  Let  $D=\{x\in {\Bbb R}^n \big|
 0\le x_i\le l_i, \, i=1, \cdots, n\}$ with boundary $\Gamma$, and let
 $\Gamma_\varrho =\{x\in {\Bbb R}^n\big| 0\le x_i \le l_i \,\,
\mbox{when}\,\, i<n, \, x_n =0\}$. Let
 $\Gamma^{l_n} =\{x\in {\Bbb R}^n\big| 0\le x_i \le l_i \,\,
\mbox{when}\,\, i<n, \, x_n =l_n\}$.
 Our first purpose, in this section, is to discuss
 the biharmonic Steklov eigenvalue problem on
 $n$-dimensional rectangular parallelepiped $D$:
 \begin{eqnarray} \label {4-1}  \left\{ \begin{array}{ll} \triangle^2 u = 0\;
 \;&
   \mbox{in}\;\; D,\\
  \frac{\partial u}{\partial \nu}=0 \; \; &\mbox{on}\;\; \Gamma_\varrho,\\
   u=\frac{\partial u}{\partial \nu}=0 \;\; &\mbox{on}\;\;
   \Gamma^{l_n},\\
 u= \Delta u=0\;\; &\mbox{on}\;\; \Gamma -(\Gamma_{\varrho} \cup \Gamma^{l_n}),\\
   \frac{\partial(\Delta u)}{\partial \nu}-\lambda^3 \varrho^3  u=0 \;\;
    &\mbox{on}\;\; \Gamma_{\varrho}, \quad \,  \varrho =constant>0\;\;\mbox{on}\;\;
     \Gamma_\varrho.
   \end{array}\right.\end{eqnarray}

   We consider nonzero product solution of (\ref{4-1}) of the form:
\begin{eqnarray*} u=X(x_1,\cdots, x_{n-1})\, Y(x_n),\end{eqnarray*}
where $X(x_1,\cdots, x_{n-1})$ is a function of variables $x_1,
\cdots, x_{n-1}$ and $Y(x_n)$ is a function of $x_n$ alone.
 Since \begin{eqnarray*} \Delta u&=&
\big(\Delta_{n-1} X(x_1, \cdots, x_{n-1})\big)Y(x_n)+2\nabla X (x_1,
\cdots, x_{n-1}) \cdot \nabla Y(x_n) \\ && +\big(X(x_1, \cdots,
 x_{n-1})\big)Y''(x_n)=
  \big(\Delta_{n-1} X(x_1, \cdots, x_{n-1})\big)Y(x_n)\\&&
  +\big(X(x_1, \cdots,
 x_{n-1})\big)Y''(x_n) \end{eqnarray*}
 and
    \begin{eqnarray*}  \Delta^2 u &=&(\Delta^2_{n-1} X(x_1, \cdots, x_{n-1}))Y(x_n)
    + 2 \big(\Delta_{n-1} X(x_1, \cdots, x_{n-1})\big) Y'' (x_n)\\
    &&+ \big(X(x_1,\cdots, x_{n-1})\big)Y''''(x_n),\end{eqnarray*}
where $$\Delta_{n-1} X(x_1, \cdots, x_{n-1}) =\sum_{i=1}^{n-1}
\frac{\partial^2
 X}{\partial x_i^2},$$
 we find by  $\Delta^2 u=0$ that
  \begin{eqnarray*} & \big(\Delta^2 X(x_1, \cdots, x_{n-1})\big) Y(x_n)
 +2 (\Delta X(x_1, \cdots, x_{n-1})\big) Y''(x_n)\\
 &+\big(X(x_1, \cdots, x_{n-1})\big)
 Y''''(x_n)=0,\end{eqnarray*}
  so that  \begin{eqnarray} \label {4..2}\quad\; \; \frac{\Delta^2 X(x_1, \cdots,
 x_{n-1})}{X(x_1, \cdots, x_{n-1})} +2\,\frac{\Delta X(x_1, \cdots,
 x_{n-1})}{X(x_1, \cdots, x_{n-1})}\, \frac{Y''(x_n)}{Y(x_n)} +
 \frac{Y''''(x_n)}{Y(x_n)}=0.\end{eqnarray}
Differentiating (\ref{4..2}) with respect to $x_n$, we obtain that
\begin{eqnarray*} \label{4..3}  2 \frac{\Delta X(x_1, \cdots, x_{n-1})}{ X(x_1, \cdots,
x_{n-1})}\left[\frac{Y''(x_n)}{Y(x_n)}\right]'
+\left[\frac{Y''''(x_n)}{Y(x_n)}\right]'=0.\end{eqnarray*} The above
equation holds if and only if
\begin{eqnarray} \label{4..4}   \frac{\Delta X(x_1, \cdots, x_{n-1})}{ X(x_1, \cdots,
x_{n-1})} =-\frac{\left[\frac{Y''''(x_n)}{Y(x_n)}\right]'}{ 2
\left[\frac{Y''(x_n)}{Y(x_n)}\right]'}= -\eta^2,\end{eqnarray} where
$\eta^2$ is a constant.  Therefore, we have that
\begin{eqnarray}\label{4..5} \Delta X(x_1, \cdots, x_{n-1}) +\eta^2 X(x_1,
\cdots, x_{n-1}) =0 \end{eqnarray} and
$$ \left[\frac{Y''''(x_n)}{Y(x_n)}\right]' -2 \eta^2
\left[\frac{Y''(x_n)}{Y(x_n)}\right]'=0.$$
 From (\ref{4..5}), we get \begin{eqnarray} \Delta^2 X
=-\eta^2 \Delta X =\eta^4 X. \end{eqnarray}
  Substituting this into  (\ref{4..2}), we obtain the following equation
\begin{eqnarray} \label{4..6} Y''''(x_n)- 2\eta^2 Y''(x_n) +\eta^4
Y(x_n)=0.\end{eqnarray}
   It is easy to verify  that the general solutions of (\ref{4..6}) have
   the form:
  \begin{eqnarray} \label{4.6'} Y(x_n) = A \cosh \, \eta x_n +B\sinh\, \eta
  x_n  +  C x_n \cosh\, \eta x_n +D x_n \sinh\, \eta x_n.\end{eqnarray}
 By setting   $Y(0)=1,\; \; Y(l_n)=0, \;\; Y'(0)=0, \; \; Y'(l_n)=0$, we get
 \begin{eqnarray}  \label{4--7} \quad \;  Y(x_n)&=&
 \cosh \eta x_n -
 \left[\frac{(\sinh \eta l_n)(\cosh \eta l_n)+\eta l_n}{\sinh^2 \eta l_n
  - \eta^2 l_n^2}\right] \sinh \eta x_n \\
  && +\left[
  \frac{\eta (\sinh \eta l_n) (\cosh \eta l_n) +\eta^2 l_n}{
 \sinh^2 \eta l_n -\eta^2 l_n^2} \right] x_n \cosh \eta x_n\nonumber
 \\ && - \left[\frac{\eta \sinh^2 \eta l_n}{
 \sinh^2 \eta l_n -\eta^2 l_n^2}\right] x_n\sinh \eta x_n.\nonumber
   \end{eqnarray}
     It is well-known that for the Dirichlet eigenvalue problem
     \begin{eqnarray}  \left\{ \begin{array} {ll} \Delta X (x_1,\cdots, x_{n-1})+ \eta^2
      X(x_1, \cdots, x_{n-1})=0
     \quad \; \mbox{in}\;\; D,\\
     u=0 \quad \; \mbox{on}\; \; \partial \{(x_1,\cdots, x_{n-1}) \big|
  0\le x_i\le l_i, \,\, i=1, \cdots, n-1\},
     \end{array} \right. \end{eqnarray}
  there exist the eigenfunctions
  \begin{eqnarray} \label{4..10} X(x_1, \cdots, x_{n-1})= c \left(\sin \frac{m_1\pi}{l_1}x_1\right)\cdots
  \left(\sin
  \frac{m_{n-1} \pi } {l_{n-1}}x_{n-1}\right), \end{eqnarray}
  which correspond to the eigenvalues
   \begin{eqnarray*}  \eta^2=\sum_{i=1}^{n-1} \left(\frac{m_i \pi}{l_i}\right)^2,\quad \;\;
 \mbox{where}\;\; m_i=1, 2,3, \cdots. \end{eqnarray*}
    Therefore,
    \begin{eqnarray} \label{4.11}  \quad \;\,\; u&=& \big(X(x_1, \cdots, x_{n-1})\big)Y(x_n)
    \\
    &= & c \left(\sin \frac{m_1\pi}{l_1}x_1\right)\cdots
  \left(\sin \frac{m_{n-1} \pi } {l_{n-1}}x_{n-1}\right)
\left[\cosh \eta x_n\nonumber \right.\\&&\left. -
 \left(\frac{(\sinh \eta l_n)(\cosh \eta l_n)+\eta l_n}{\sinh^2 \eta l_n
  - \eta^2 l_n^2}\right) \sinh \eta x_n \right.\nonumber\\
  && \left.+\left(
  \frac{\eta (\sinh \eta l_n) (\cosh \eta l_n) +\eta^2 l_n}{
 \sinh^2 \eta l_n -\eta^2 l_n^2} \right) x_n \cosh \eta x_n\right.\nonumber
 \\ && \left.- \left(\frac{\eta \sinh^2 \eta l_n}{
 \sinh^2 \eta l_n -\eta^2 l_n^2}\right) x_n\sinh \eta x_n\right].\nonumber
  \end{eqnarray}
 Since $$Y'''(0)= 2 \eta^3 \left(
\frac{(\sinh \eta l_n)(\cosh \eta l_n)+\eta l_n}{\sinh^2 \eta l_n
-\eta^2 l_n^2}\right),$$
 we obtain
 \begin{eqnarray*}\frac{\partial (\triangle u)}{\partial \nu}\big|_{x_n=0}&=&
  \big(\Delta_{n-1} X(x_1,\cdots, x_{n-1})\big)Y'(0) +\big(X
  (x_1, \cdots, x_{n-1})\big)Y'''(0)\\
    &=& 2\eta^3 \left( \frac{(\sinh \eta l_n)(\cosh \eta l_n)+\eta l_n}{\sinh^2
\eta l_n -\eta^2 l_n^2}\right)X(x_1, \cdots,
x_{n-1}),\end{eqnarray*} so that
$$\frac{\partial (\triangle u)}{\partial \nu} -\lambda^3 \varrho^3  u
=0\quad \; \mbox{on}\;\; \Gamma_\varrho$$ with
  $$\lambda^3 = \frac{2\eta^3 l_n^3}{\rho^3 l_n^3} \left(
\frac{(\sinh \eta l_n)(\cosh \eta l_n) +l_n\eta}{\sinh^2 \eta l_n
-\eta^2 l_n^2}\right).$$

\vskip 0.29 true cm

  Our second purpose is to discuss the biharmonic
 Steklov eigenvalue problem on the $n$-dimensional rectangular parallelepiped $D$:
\begin{eqnarray} \label {4-13}  \left\{ \begin{array}{ll} \triangle^2 u = 0\quad \;
   \mbox{in}\;\; D,\\
   \frac{\partial u}{\partial \nu}=0 \; \; \mbox{on}\;\; \Gamma_\varrho, \quad \,
   u=\frac{\partial u}{\partial \nu}=0 \;\;\mbox{on}\;\;
   \Gamma^{l_n},\\
      \frac{\partial (\Delta u)}{\partial \nu}= 0
      \;\; \mbox{on}\;\; \Gamma-(\Gamma_\varrho\cup \Gamma^{l_n}),\\
    \frac{\partial (\triangle u)}{\partial \nu}-\lambda^3 \varrho^3 u=0 \;\;
    \mbox{on}\;\; \Gamma_\varrho, \quad \,  \varrho =constant>0\;\;\mbox{on}\;\;
     \Gamma_\varrho.
   \end{array}\right.\end{eqnarray}

\vskip 0.22 true cm

  Similarly, (\ref{4-13}) has the special solution
  $u=\big(X(x_1, \cdots, x_{n-1})\big)Z(x_n)$ with $Z(x_n)$ having form
 (\ref{4.6'}). According to the boundary conditions of (\ref{4-13}),
 we get that the problem (\ref{4-13}) has the solutions
 \begin{eqnarray*}  \label{4-19} u(x)&=&  u(x_1, \cdots, x_n)\\
   &=& c\left(\cos \frac{m_1 \pi}{l_1}\, x_1\right)\cdots
   \left(\cos \frac{m_{n-1} \pi}{l_{n-1}}\, x_{n-1}\right)
  Z(x_n),\nonumber \end{eqnarray*}
where $m_1, \cdots, m_{n-1}$ are whole numbers, and $Z (x_n)$ is
given by
\begin{eqnarray}  \label{4-9} \quad \;  Z(x_n)&=&
 \cosh \beta x_n -
 \left[\frac{(\sinh \beta l_n)(\cosh \beta l_n)+\beta l_n}{\sinh^2 \beta l_n
  - \beta^2 l_n^2}\right] \sinh \beta x_n \\
  && +\left[
  \frac{\beta (\sinh \beta l_n) (\cosh \beta l_n) +\beta^2 l_n}{
 \sinh^2 \beta l_n -\beta^2 l_n^2} \right] x_n \cosh \beta x_n\nonumber
 \\ && - \left[\frac{\beta \sinh^2 \beta l_n}{
 \sinh^2 \beta l_n -\beta^2 l_n^2}\right] x_n\sinh \beta x_n,\nonumber
   \end{eqnarray}
 $\beta=\big[\sum_{i=1}^{n-1}
(m_i\pi/l_i)^2\big]^{1/2}$ with $\sum_{i=1}^{n-1} m_i\ne 0$. Since
  $Z'''(0)=2\beta^3
   \left(\frac{(\sinh \beta l_n) (\cosh \beta l_n)+\beta l_n}{
 \sinh^2 \beta l_n -\beta^2 l_n^2}\right)$
 and  $\frac{\partial (\Delta u)}{\partial \nu}\big|_{x_n=0}=\big(X(x_1,\cdots, x_{n-1})\big)
 Z'''(0)$, we get
$\frac{\partial (\Delta u)}{\partial \nu}-\lambda^3 \varrho^3
  u=0$ on $\Gamma_\varrho$, where
$$\lambda^3 = \frac{2\beta^3 l_n^3}{\rho^3 l_n^3} \left( \frac{(\sinh \beta
l_n)(\cosh \beta l_n) +\beta l_n}{\sinh^2 \beta l_n -\beta^2
l_n^2}\right).$$

\vskip 1.39 true cm

\section{Asymptotic distribution of eigenvalues on special domains}

\vskip 0.45 true cm

\noindent {\bf 5.1. Counting function $A(\tau)$.}

\vskip 0.2 true cm

In order to obtain our asymptotic formula, it is an effective way
 to investigate the
distribution of the eigenvalues of the transformation
$G_{\mathcal{K}}$ (respectively, $G_{{\mathcal{K}}^d}$,
$\,G_{\mathcal{N}}$) instead of the transformations
$G_{\mathcal{K}}^{(\star)}$ (respectively,
  $ G_{{\mathcal {K}}^d}^{(\star)}$, $\,G_{\mathcal{N}}^{(\star)}$).
 It follows from (\ref{3.10}), (\ref{3.11}), (\ref{3.12}) and (\ref{3-37}) we obtain
 \begin{eqnarray} \mu_k= (1+\lambda_k^3)^{-1}, \quad \; k=1, 2, 3,
 \cdots, \end{eqnarray}
 where $\mu_k$ denote the $k$-th eigenvalue of
 $G_{\mathcal{K}}$ or
 $G_{{\mathcal{K}}^d}$ or
 $G_{\mathcal{N}}$, and $\frac{1}{\lambda_k}$ is the
 $k$-th eigenvalue of $G_{\mathcal{K}}^{(\star)}$ or $G_{{\mathcal{K}}^d}^{(\star)}$
  or $G_{\mathcal{N}}^{(\star)}$.
 Since $A(\tau)= \sum_{\lambda_k\le \tau} 1$, we have
 \begin{eqnarray} A(\tau)= \sum_{\mu_k \ge (1+\tau^3)^{-1}} 1.
 \end{eqnarray}

\vskip 0.48 true cm

\noindent {\bf 5.2. $\,D$ is an $n$-dimensional  rectangular
parallelepiped and
 $g^{ik}=\delta^{ik}$}.

\vskip 0.3 true cm

 Let $D$ be an $n$-dimensional rectangular parallelepiped,
 $g^{ik}=\delta^{ik}$ in the whole of $\bar D$,
  $\varrho=constant>0$ on one face $\Gamma_{\varrho}^+$ of the
  rectangular parallelepiped and $\varrho=0$ on $\Gamma_\varrho-\Gamma_\varrho^+$,
  i.e., $D=\{x\in {\Bbb R}^n\big|
 0\le x_i\le l_i, \, i=1, \cdots, n\}$,
 $\Gamma_\varrho^+ =\{x\in {\Bbb R}^n\big| 0\le x_i \le l_i \,\,
\mbox{when}\,\, i<n, \, x_n =0\}$) and  $\Gamma_0=\Gamma^{l_n}
=\{x\in {\Bbb R}^n\big| 0\le x_i \le l_i \,\, \mbox{when}\,\, i<n,
\, x_n =l_n\}$.
 Without loss of generality, we assume $l_i < l_n$ for all $i<n$.

For the above domain $D$, except for the $K(D)$ and $K^d (D)$ in
Section 3, we introduce the linear space of functions
     \begin{eqnarray*}  K^0(D)= \{u\big| u\in Lip(\bar D) \cap
     H^2 (D), \;  \frac{\partial u}{\partial \nu}=0\;\;
  \mbox{on}\;\; \Gamma_\varrho \cup \Gamma_{0}, \;\; u=0 \;\;
  \mbox{on}\;\; \Gamma-\Gamma_\varrho^+\},\end{eqnarray*}
 Clearly,
   \begin{eqnarray}  \label {5-3} K^0(D) \subset K(D) \subset K^d
  (D),\end{eqnarray}
     Closing $K^0$, $K$ and $K^d$ respect to the norm $\|u\|=
   \sqrt{\langle u, u\rangle}$, we obtain the Hilbert spaces
   ${\mathcal{K}}^0$, $\mathcal{K}$ and ${\mathcal{K}}^d$, and
  \begin{eqnarray} {\mathcal{K}}^0 \subset \mathcal{K} \subset
  {\mathcal{K}}^d.\end{eqnarray}
   According to Theorem 3.3, we see that the bilinear functional
   \begin{eqnarray} \label{5-4} [u,v] =\int_{\Gamma_\varrho^+}
   \varrho^3
   uv   \, ds\end{eqnarray}
   defines self-adjoint, completely continuous transformations $G^0$,
   $G$ and $G^d$ on ${\mathcal{K}}^0$, $\mathcal{K}$ and
   ${\mathcal{K}}^d$, respectively (cf. Section 3).
 Obviously,   \begin{eqnarray*} \label{5-5}  \langle G^0 u, v\rangle
= \langle Gu,
      v\rangle
   \quad \; \mbox{for all}\;\; u, v  \;\; \mbox{in}\;\;
   {\mathcal{K}}^0,\end{eqnarray*}
     \begin{eqnarray*}  \label{5-6} \langle G u, v\rangle = \langle G^d u,
     v\rangle
   \quad \; \mbox{for all}\;\; u, v  \;\; \mbox{in}\;\;
   {\mathcal{K}},\end{eqnarray*}
  from which and applying Theorem 1.4 of \cite{Sa} we immediately get
  \begin{eqnarray} \mu_k^0 \le \mu_k \le \mu_k^d, \quad k=1,2,3,
  \cdots, \end{eqnarray}
  where $\{\mu_k^0\}$ and $\{\mu_k^d\}$ are the eigenvalues of $G^0$
  and $G^d$, respectively. Hence
  \begin{eqnarray} \label{5-66} A^0(\tau) \le A(\tau) \le A^d(\tau) \quad \;\;
  \mbox{for all}\;\; \tau,\end{eqnarray}
  where \begin{eqnarray} \label {5-7} A^0(\tau)=\sum_{\mu_k^0\ge (1+\tau^3)^{-1}}
1\end{eqnarray}
 and \begin{eqnarray} \label {5-8} A^d(\tau)=\sum_{\mu_k^d\ge (1+\tau^3)^{-1}}
1.\end{eqnarray}

 We shall estimate the asymptotic behavior
  of $A^0(\tau)$ and $A^d(\tau)$.
 It is easy to verify (cf. Theorems 3.6, 3.7) that
 the eigenfunctions of the transformations $G^0$ and $G^d$,
  respectively, satisfy
 \begin{eqnarray} \label {5-13} \quad \;\; \left\{ \begin{array}{ll} \triangle^2 u = 0\quad \;
   \mbox{in}\;\; D,\\
   \frac{\partial u}{\partial \nu}=0 \; \; \mbox{on}\;\; \Gamma_\varrho^+,\quad \\
   u=\frac{\partial u}{\partial \nu} =0 \;\; \mbox{on}\;\; \Gamma^{l_n},
\,\, \mbox{and}\;\;
u=\Delta u=0\;\;\mbox{on}\;\; \Gamma-(\Gamma_\varrho^+ \cup \Gamma^{l_n}),\\
  \frac{\partial (\triangle u)}{\partial \nu}-\lambda^3 \varrho^3
   u=0 \;\;
    \mbox{on}\;\; \Gamma_{\varrho}^+, \quad \,  \varrho =constant>0\;\;\mbox{on}\;\;
     \Gamma_\varrho^+.
   \end{array}\right.\end{eqnarray}
and
\begin{eqnarray} \label {5-14} \quad  \left\{ \begin{array}{ll} \triangle^2 u = 0\quad \;
   \mbox{in}\;\; D,\\
   \frac{\partial u}{\partial \nu}=0 \; \; \mbox{on}\;\; \Gamma_\varrho^+,\quad \;\; u=
   \frac{\partial u}{\partial \nu}=0 \;\; \mbox{on}\;\;
   \Gamma^{l_n},\\
    \frac{\partial (\Delta u)}{\partial \nu}=
    0 \;\; \mbox{and}\;\; \Delta u=0\;\; \mbox{on}\;\;
    \Gamma-(\Gamma_\varrho^+\cup \Gamma^{l_n}),\\
    \frac{\partial (\triangle u)}{\partial \nu}-\kappa^3 \varrho^3 u=0 \;\;
    \mbox{on}\;\; \Gamma_\varrho^+, \quad \,  \varrho =constant>0\;\;\mbox{on}\;\;
     \Gamma_\varrho^+.
   \end{array}\right.\end{eqnarray}

As being verified in Section 4, the functions of form
\begin{eqnarray} \label{-5-15} u(x)=
c \left(\sin \frac{m_1\pi}{l_1}x_l\right)\cdots
  \left(\sin
  \frac{m_{n-1} \pi } {l_{n-1}}x_{n-1}\right)
 Y (x_n)\end{eqnarray}
 are the solutions of the problem
(\ref{5-13}), where $Y(x_n)$ is given by (\ref{4--7}). Since the
functions in (\ref{-5-15})
 have derivatives of any order in $D$, it follows from Proposition 3.9 that they
  are eigenfunctions of the transformation $G^0$ with eigenvalues $(1+\lambda^3)^{-1}$,
 where
\begin{eqnarray} \label{-5-16} &&  \; \lambda^3  =\frac{2\eta^3 l_n^3}{\varrho^3 l_n^3}
  \left(\frac{(\sinh \eta l_n)(\cosh \eta l_n) + \eta l_n}{\sinh^2 \eta \l_n -\eta^2
  l_n^2}\right),\\
  && \quad \;\; \mbox{and}\;\;  \eta= \left[ \sum_{i=1}^{n-1}
  \big(\frac{m_i\pi}{l_i}\big)^2 \right]^{1/2},\;\; i=1,2,3, \cdots.\nonumber\end{eqnarray}
    Note that the restriction of $u$ on $\Gamma_\varrho$
       \begin{eqnarray} \label{5;17}  u\big|_{\Gamma_\varrho} =
   c\left(\sin \frac{m_1 \pi}{l_1}  x_1\right)
   \cdots \left(\sin \frac{m_{n-1} \pi}{l_{n-1}} x_{n-1}\right), \end{eqnarray}
 when $m_1, \cdots, m_{n-1}$ run through all positive integers (see, Section 4), form
  a complete system of orthogonal functions in $L^2_\varrho(\Gamma_\varrho)$.
It follows from Proposition 3.5 that if $m_1,
  \cdots, m_{n-1}$ run through all positive integers,
  then the functions (\ref{-5-15})
 form an orthogonal basis of the subspace of ${\mathcal{K}}^0$, spanned by the
 eigenfunctions of $G^0$, corresponding to
  positive eigenvalues. That is, when
   $m_1, \cdots,  m_{n-1}$ run through all
  positive integers, then $(1+\lambda^3)^{-1}$,
   where $\lambda^3$ is given by (\ref{-5-16}), runs through all positive eigenvalues of $G^0$.

 Similarly, for the problem (\ref{5-14}), the eigenfunctions $\{u_k\}$ of the operator
 $G^d$ on ${\mathcal{K}}^d$, corresponding to non-zero eigenvalues,
  form  an orthogonal basis
 of the subspace of ${\mathcal{K}}^d$.
 The non-zero eigenvalues  of $G^d$ are $\mu^d_k=(1+\kappa_k^3)^{-1}$, where $\kappa_k^3$
 is the $k$-th Steklov eigenvalue of (\ref{5-14}).

In order to give the upper bound estimate of $A^d(\tau)$, we further
introduce the following Steklov
 eigenvalue problem
\begin{eqnarray} \label {5.*.15}  \left\{ \begin{array}{ll} \triangle^2 u = 0\quad \;
   \mbox{in}\;\; D,\\
   \frac{\partial u}{\partial \nu}=0 \; \; \mbox{on}\;\; \Gamma_\varrho^+, \quad \,
   u=\frac{\partial u}{\partial \nu}=0 \;\;\mbox{on}\;\;
   \Gamma^{l_n},\\
      \frac{\partial u}{\partial \nu}=\frac{\partial (\Delta u)}{\partial \nu}= 0
      \;\; \mbox{on}\;\; \Gamma-(\Gamma_\varrho^+\cup \Gamma^{l_n}),\\
    \frac{\partial (\triangle u)}{\partial \nu}-\gamma^3 \varrho^3  u=0 \;\;
    \mbox{on}\;\; \Gamma_\varrho^+, \quad \,  \varrho =constant>0\;\;\mbox{on}\;\;
     \Gamma_\varrho^+.
   \end{array}\right.\end{eqnarray}
 Let $\gamma_k^3$ be the $k$-th eigenvalue of (\ref{5.*.15}). By Theorem 3.8, we have
 \begin{eqnarray}\label{5.-.15}  \gamma_k^3\le \kappa_k^3\quad \;\; \mbox{for all }\;\;
 k\ge 1.\end{eqnarray}
   We define \begin{eqnarray} \label{5..0}\mu_k^f=\frac{1}{1+\gamma_k^3}, \quad \;\; A^f(\tau)
   = \sum_{\mu_k^f\ge (1+\tau^3)^{-1}} 1.\end{eqnarray}
It follows from (\ref{5.-.15}) and (\ref{5..0}) that
\begin{eqnarray} \label{5*1} A^d (\tau)\le A^f(\tau)\quad \; \mbox{for all}\;\; \tau.\end{eqnarray}

  We know (cf. Section 4) that the problem (\ref{5.*.15}) has the solutions of form
 \begin{eqnarray}  \label{5-21} u(x)=c \left(\cos \frac{m_1 \pi}{l_1} x_1\right)\cdots
 \left(\cos \frac{m_{n-1} \pi}{l_{n-1}} x_{n-1}\right)
  Z (x_n), \end{eqnarray}
where $m_1, \cdots, m_{n-1}$ are non-negative integers with
$\sum_{i=1}^{n-1}m_{i}\ne 0$, and $Z (x_n)$ is given by (\ref{4-9}).
 This implies that if $m_1, \cdots, m_{n-1}$  run through all non-negative integers with
  $\sum_{i=1}^{n-1}m_{i} \ne 0$, then
\begin{eqnarray} \quad \quad \; \label{5--16}   \gamma^3  =\frac{2\beta^3 l_n^3}{\rho^3 l_n^3}
  \left(\frac{(\sinh \beta l_n)(\cosh \beta l_n) + \beta l_n}{\sinh^2 \beta \l_n -\beta^2 l_n^2}\right),
  \quad \; \beta =  \left[ \sum_{i=1}^{n-1}
  \big(\frac{m_i\pi}{l_i}\big)^2 \right]^{1/2}\end{eqnarray}
 runs throughout
 all eigenvalues of problem (\ref{5.*.15}).

 We first compute the asymptotic behavior of $A^f (\tau)$.
By (\ref{5..0}), (\ref{5--16}) and the argument as in p.$\,$44 of
\cite{We4} or p.$\,$373 of \cite{CH} or p.$\,$51-53 of \cite{Sa},
$A^f(\tau)=$the number of $(n-1)$-tuples $(m_1, \cdots, m_{n-1})$
satisfying the inequality
\begin{eqnarray} \quad \quad \; \label{5-22'}  \frac{2\beta^3 l_n^3}{\varrho^3 l_n^3}
  \left(\frac{(\sinh \beta l_n)(\cosh \beta l_n) +
  \beta l_n}{\sinh^2 \beta \l_n -\beta^2 l_n^2}\right)\le \tau^3,
  \end{eqnarray}
  where $m_1, \cdots,
 m_{n-1}$ are non-negative integers with $\sum_{i=1}^{n-1} m_i \ne
 0$.
   By setting \begin{eqnarray} \label {005.10} t(s)=2s^3 \left(
\frac{(\sinh s)(\cosh s) +s}{\sinh^2 s -s^2}\right),\end{eqnarray}
  we see that $$\lim_{s\to +\infty} t(s)/s^3= 2.$$
 We claim that
for all $s\ge 1$, \begin{eqnarray*} t'(s) =\frac{2s^2 \left[ -s^3
+3s\, \sinh^2 s +3 (\sinh^3 s)(\cosh s) -3 s^2 (\sinh s)(\cosh
s)-2s^3 \cosh^2 s\right]}{(\sinh^2 s -s^2)^2}
>0.\end{eqnarray*}
  In fact, let \begin{eqnarray*}
\theta(s)=- s^3 +3s \,\sinh^2 s +3 (\sinh^3 s) (\cosh s) -3 s^2
(\sinh s)(\cosh s) -2s^3 \cosh^2 s.\end{eqnarray*} Then
\begin{eqnarray*} &&\theta (1)>0, \,\, \mbox{and}\\
  &&\theta'(s)= 4(\cosh^2 s) \left[ 3\sinh^2 s -3 s^2 -s^3 \frac{\sinh s}{\cosh s}\right]\\
&& \quad \quad \; > 4(\cosh^2 s) \left[ 3\sinh^2 s -3 s^2 -s^3
\right]\\
&& \quad \quad \; >  4(\cosh^2 s) \left[
\frac{3}{4}(e^{2s}+e^{-2s})-\frac{3}{2}-3 s^2
 -s^3\right]\\
  && \quad \quad \; >4(\cosh^2 s) \left[ \frac{3}{4} \left(2+4s^2 +\frac{4}{3}s^4\right)-\frac{3}{2} -3s^2
-s^3 \right]\\
 &&\quad \quad \; = 4(\cosh^2 s) (s^4-s^3) \ge 0 \quad \; \mbox{for}\;\; s\ge
 1.
 \end{eqnarray*}
 This implies that $\theta(s)\ge 0$ for $s\ge 1$.
  Thus, the function ${t}(s)$ is increasing in $[1, +\infty)$.
  Denote by $s={h}(t)$ the inverse of function
 ${t}(s)$ for $s\ge 1$. Then $$\lim_{t\to +\infty} \frac{({h}(t))^3}{t}=\frac{1}{2}.$$
  Note that, for $s\ge 1$, the inequalities $t(s)\le t$ is equivalent to
  $s\le h(t)$.
 Hence (\ref{5-22'}) is equivalent to
\begin{eqnarray*}  \beta l_n  \le h(\varrho^3 l_n^3 \tau^3), \end{eqnarray*}
which can be written as
  \begin{eqnarray} \label{5-23} \quad \quad\; \sum_{i=1}^{n-1}
  (m_i/l_i)^2 \le \left[\frac{1}{\pi l_n} h(\varrho^3 l_n^3 \tau^3)\right]^2,
  \quad \, m_i=0, 1,2,\cdots\;\;\mbox{with}\;\; \sum_{i=1}^{n-1}m_i\ne 0. \end{eqnarray}
 We consider the $(n-1)$-dimensional ellipsoid
\begin{eqnarray*} \label{5-24}  \sum_{i=1}^{n-1}
  (z_i/l_i)^2 \le \left[\frac{1}{\pi l_n} h(\varrho^3 l_n^3 \tau^3)\right]^2. \end{eqnarray*}
  Since $A^f(\tau)+1\,$ just is the number of those $(n-1)$-dimensional unit cubes of the
  $z$-space that have corners whose
  coordinates are non-negative integers in the ellipsoid (see, VI. \S4 of \cite{CH}).
 Hence $A^f(\tau)+1$ is the sum of the volumes of these cubes.
  Let $V(\tau)$ denote the volume and $T(\tau)$ the
  area of the part of the ellipsoid situated in
  the positive octant $z_i\ge 0,\, i=1,\cdots, n-1$. Then
$$ V(\tau) \le A^f(\tau)+1 \le V(\tau)+ (n-1)^{\frac{1}{2}} T(\tau),$$
where $(n-1)^{\frac{1}{2}}$ is the diagonal length of the unit cube
(see, \cite{CH} or \cite{Sa}).  Since
  \begin{eqnarray*} \label{5-25}  V(\tau)= D_{n-1}
  2^{-(n-1)} l_1\cdots l_{n-1}\left[\frac{h(\varrho^3 l_n^3
  \tau^3)}{\pi l_n}\right]^{(n-1)},
 \end{eqnarray*}
  by $h(t)\sim \left(\frac{t}{2}\right)^{1/3} \;\;\mbox{as}\;\; t\to +\infty$, we get that
  \begin{eqnarray*}  \label{5-26}  V(\tau)\sim D_{n-1}
 (\sqrt[3]{16}\,\pi)^{-(n-1)} l_1\cdots l_{n-1} \varrho^{n-1}
 \tau^{n-1}, \quad \; \mbox{as}\;\; \tau\to +\infty.\end{eqnarray*}
   Note that
\begin{eqnarray*}
  T(\tau)\sim \mbox{constant}\cdot \tau^{n-2}.\end{eqnarray*}
    It follows that
      \begin{eqnarray*} \lim_{\tau\to +\infty} \frac{A^f (\tau)}{\tau^{n-1}}= \omega_{n-1}
  (\sqrt[3]{16}\,\pi)^{-(n-1)} l_1 \cdots l_{n-1} \varrho^{n-1}, \end{eqnarray*}
  i.e.,
  \begin{eqnarray} \label{5.2.8}  A^f(\tau) \sim \frac{\omega_{n-1}}{(\sqrt[3]{16}\,\pi)^{(n-1)}}
  |\Gamma_\varrho^+|\varrho^{n-1} \tau^{n-1},\;\; \mbox{as}\;\; \tau \to +\infty,\end{eqnarray}
 where $|\Gamma_\varrho^+|$ denotes the area of the face $\Gamma_\varrho^+$.

Next, we consider $A^0(\tau)$.
    Similarly,
\begin{eqnarray} \quad \quad \; \label{5-22}  \frac{2\eta^3 l_n^3}{\varrho^3 l_n^3}
  \left(\frac{(\sinh \eta l_n)(\cosh \eta l_n) + \eta l_n}{\sinh^2 \eta \l_n -\eta^2 l_n^2}\right)\le \tau^3,
  \end{eqnarray}
   is equivalent to
  \begin{eqnarray*} \eta l_n \le {h}(\varrho^3 l_n^3 \tau^3),\end{eqnarray*}
  i.e.,
  \begin{eqnarray*} \sum_{i=1}^{n-1} \big[m_i /l_i\big]^2 \le \left( \frac{h
  (\varrho^3 l_n^3 \tau^3)}{\pi l_n}\right)^2, \;\; m_i=1,2,3,\cdots.\end{eqnarray*}
    Similar to the argument for $A^f(\tau)$, we find  (see also, \S4 of \cite{CH}) that
    \begin{eqnarray*}  \# \{ (m_1, \cdots, m_{n-1})\big|
  \sum_{i=1}^{n-1} \big(\frac{m_i}{l_i}\big)^2 \le \big( \frac{h
  (\varrho^3 l_n^3 \tau^3)}{\pi l_n}\big)^2,\;\, m_i= 1,2,3,\cdots \big\} \\
     \sim
   \frac{\omega_{n-1}}{(\sqrt[3]{16}\,\pi)^{(n-1)}}
  |\Gamma_\varrho^+|\varrho^{n-1} \tau^{n-1} \;\; \mbox{as}\;\; \tau\to +\infty.\qquad \quad\end{eqnarray*}
  i.e.,
    \begin{eqnarray} \label {5-30}   \lim_{\tau \to +\infty}
    \frac{A^0(\tau)}{\tau^{n-1}}= \frac{\omega_{n-1}}{(\sqrt[3]{16}\,\pi)^{(n-1)}}
  |\Gamma_\varrho^+|\varrho^{n-1}.\end{eqnarray}
 Not that $\varrho=0$ on $\partial D -\Gamma_\varrho^+$, by
  (\ref{5-66}), (\ref{5*1}), (\ref{5.2.8}) and (\ref{5-30}), we have
  \begin{eqnarray} \label{5-31} A(\tau) \sim \frac{\omega_{n-1}\tau^{n-1}}{(\sqrt[3]{16}\,\pi)^{(n-1)}  }
  \int_{\Gamma_{\varrho}} \varrho^{n-1}ds \quad \;
   \mbox{as}\;\; \tau\to +\infty. \end{eqnarray}

\vskip 0.63 true cm

\noindent {\bf 5.3.  $\,$A cylinder $D$ whose base is an
 $n$-polyhedron of ${\Bbb R}^{n-1}$ having $n-1$ orthogonal plane surfaces and
  $g^{ik}=\delta^{ik}$}.

\vskip 0.38 true cm

 \noindent  {\bf Lemma 5.1.}  {\it Let $D^{(j)}
 =\Gamma^{(j)}_\varrho \times [0, l_n]$, $\,j=1,2$,
    where
    $\Gamma^{(1)}_\varrho=\{(x_1,\cdots, x_{n-1})\in {\Bbb R}^{n-1}\big|
    \, x_i\ge 0\;\;\mbox{for}\;\; 1\le i\le n-1,\,\,\mbox{and}\,\, \sum_{i=1}^{n-1} \frac{x_i}{l_i}\le 1\}$,
  and $\Gamma^{(2)}_\varrho$ is an $(n-1)$-dimensional
  cube with side length $l=\max_{1\le i\le n-1} l_i$.
     Assume that
    $\Gamma^{(j)}_{0}=\Gamma^{(j)}_\varrho\times \{l_n\}$, $\; j=1,2$.
     Assume also that  $\varrho$ is a positive constant
     on $\Gamma^{(j)}_\varrho$, $\;j=1,2$. If $l<l_n$, then
      \begin{eqnarray} \label {7-01}\gamma_k^f(D^{(1)}) \ge  \gamma_k^f(D^{(2)}) \quad \; \mbox{for}
    \;\;  k=1,2,3,\cdots,\end{eqnarray}
 where $(\gamma_k^f (D^{(j)}))^3$ (similar to $\gamma^3$ of (\ref{3-30-0}) in Theorem 3.8)
  is the $k$-th Steklov eigenvalue for
 the domain $D^{(j)}$ .}

\vskip 0.14 true cm

 \noindent  {\bf Proof.} \ \  Let $v^{(j)}_k$  be the $k$-th Neumann
 eigenfunction
 corresponding to $\alpha^{(j)}_k$  for the $(n-1)$-dimensional
 domain
 $\Gamma^{(j)}_\varrho$, $(j=1,2)$, \ i.e.,
 \begin{eqnarray} \label {7-02} \left\{ \begin{array} {ll}
   \triangle v^{(j)}_k +
 \alpha^{(j)}_k v^{(j)}_k =0 \;\; &\mbox{in}\;\; \Gamma^{(j)}_\varrho,\\
  \frac{\partial v^{(j)}_k}{\partial \nu} =0 \,\, &\mbox{on}\;\; \partial
  \Gamma^{(j)}_\varrho. \end{array} \right.\end{eqnarray}
    Put \begin{eqnarray*}u^{(j)}_k(x) = \big(v^{(j)}_k(x_1,
    \cdots, x_{n-1})\big)(Z^{(j)}(x_n)) \quad \, \mbox{in}\;\;
   D^{(j)},\end{eqnarray*}
   where $Z^{(j)}(x_n)$ is as in (\ref{4-9}) with $\beta$
   being replaced by
  $\sqrt{\alpha^{(j)}_k}$. It is easy to verify that $u^{(j)}_k(x)$ satisfies
  \begin{eqnarray}  \label{7-08} \left\{\begin{array}
   {ll} \triangle^2 u^{(j)}_k=0 \;\; &\mbox{in}\;\; D^{(j)},\\
      \frac{\partial u^{(j)}_k}{\partial \nu}=0 \; \; &\mbox{on}\;\; \Gamma^{(j)}_\varrho,\\
      u^{(j)}_k=\frac{\partial u^{(j)}_k}{\partial \nu}=0 \;\; & \mbox{on}\;\;
 \Gamma_0^{(j)}\times \{l_n\},\\
  \frac{\partial u^{(j)}_k}{\partial \nu}= \frac{\partial
   (\triangle u^{(j)}_k)}{\partial \nu}=0 \;\; &\mbox{on}\;\;
 (\partial \Gamma_\varrho^{(j)})\times [0,l_n], \\
   \frac{\partial (\triangle u^{(j)}_k)}{\partial \nu}- (\gamma^f_k (D^{(j)}))^3 \varrho^3
  u_k^{(j)}=0 \;\; &\mbox{on}\;\; \Gamma^{(j)}_{\varrho}.
  \end{array}\right.\end{eqnarray}
 with
\begin{eqnarray*}   \big(\gamma_k^f (D^{(j)})\big)^3 =\frac{2\sqrt{\alpha_k^{(j)}}}{\varrho^3}
  \left(\frac{(\sinh \sqrt{\alpha_k^{(j)}} \,l_n)\cosh \sqrt{\alpha_k^{(j)}}\, l_n +
  \sqrt{\alpha_k^{(j)}}\, l_n}{\sinh^2 \sqrt{\alpha_k^{(j)}} \,\l_n -\alpha_k^{(j)} l_n^2}\right).
 \end{eqnarray*}
 It follows from p.$\,$437-438 of \cite{CH} that
 the $k$-th Neumann eigenvalue $\alpha^{(1)}_k$ for the domain $\Gamma^{(1)}_\varrho$ is
  at least as large as the $k$-th Neumann eigenvalue $\alpha_k^{(2)}$
  for the domain $\Gamma_\varrho^{(2)}$.
  Recalling that
  $2s \left( \frac{(\sinh s)\cosh s +s}{\sinh^2 s
-s^2}\right)$ is increasing when $s\ge 1$,
 we get  $$\left(\gamma_k^f (D^{(1)})\right)^3 \ge \left(\gamma_k^f(D^{(2)})\right)^3, \;\;
 k=1,2,3,\cdots$$
 if $l<l_n$. Here we have used the fact that $ \sqrt{\alpha_k^{(2)}}\,l_n
 \ge 1$ since any Neumann eigenvalue for $\Gamma_\varrho^{(2)}$ has the form
 $\sum_{i=1}^{n-1} \big(\frac{m_i\pi}{l}\big)^2$.
 In other words, if $l<l_n$, then the number
 $A^f(\tau)$ of the eigenvalues less than or equal
 to a given bound $\tau^3$ for the domain $D^{(1)}$
 is at most equal to the corresponding number of the eigenvalues for the domain $D^{(2)}$.
 \ \ $\square$

 Similarly, we can easily verify that
 the number $A^f(\tau)$
 of the eigenvalues less than or equal
 to a given bound $\tau^3$
  for the $n$-dimensional rectangular parallelepiped $D$
 is never larger than the corresponding number for
 an $n$-dimensional
 rectangular parallelepiped of the same height whose base is an
 $(n-1)$-dimensional cube and contains the base of $D$.

\vskip 0.52 true cm

\noindent {\bf 5.4. $\,D$ is a cylinder and $g^{ik}=\delta^{ik}$}.

\vskip 0.3 true cm

Let $D$ be an open $n$-dimensional cylinder in ${\Bbb R}^n$, whose
boundary consists of an $(n-1)$-dimensional cylindrical surface and
two parallel plane surfaces perpendicular to the cylindrical
surface.  Assume that
  $g^{ik}=\delta^{ik}$ in the whole of $\bar D$, that $\Gamma_\varrho$
 includes at least one of the plane surfaces, which we call
 $\Gamma_\varrho^+$, and that $\varrho$ is positive constant on $\Gamma_\varrho^+$
  and vanishes on $\Gamma_\varrho-\Gamma_\varrho^+$.
  We let the plane surface $\Gamma_\varrho^+$ be situated in the
  plane $x_n=0$ and let another parallel surface $\Gamma^{l_n}$ be situated in the plane
  $\{x\in {\mathbb{R}}^n \big|x_n =l_n\}$. We now divide the
  plane $x_n=0$ into a net of $(n-1)$-dimensional
  cubes, whose faces are parallel to the coordinate-planes in $x_n=0$. Let $\Gamma_1, \cdots,
  \Gamma_p$ be those open cubes in the net, closure of which
  are entirely contained in $\Gamma_\varrho^+$,
   and let $Q_{p+1}, \cdots, Q_q$ be the remaining open cubes, whose closure
     intersect $\Gamma_\varrho^+$. We may let the subdivision  into
   cubes  be so fine that, for every piece of the boundary of $\Gamma_\varrho^+$
   which is contained  in one of the closure cubes, the direction of the normal varies by less than a
   given angle $\vartheta$, whose size will be determined  later.
 (This can be accomplished by repeated halving of the side of cube.)
         We can make the side length $l$ of each
   cube be less than $l_n$.
   Furthermore, let $D_j$, $(j=1, \cdots, p)$, be the open $n$-dimensional
   rectangular parallelepiped with the cube $\Gamma_j$
   as a base and otherwise bounded by the ``upper''
  plane surface $\Gamma^{l_n}$ of the cylinder $\bar D$ and planes parallel to the
  coordinate-planes $x_1=0, \cdots, x_{n-1}=0$ (cf.$\,$\cite{Sa}).

  We define the linear spaces of functions
\begin{eqnarray*} K=\{u \big| u \in Lip (\bar D)\cap  H^2(D),
\;  \frac{\partial u}{\partial \nu}=0 \;\; \mbox{on}\;\;
\Gamma_{\varrho}\cup \Gamma^{l_n},\;\; u=0\,\,\mbox{on}\;\; \Gamma-\Gamma_\varrho\},\qquad \quad \quad \\
 K_j^0 =\{u_j \big| u_j \in Lip(\bar D_j)\cap H^2(D_j),\;
 \frac{\partial u}{\partial \nu}=0 \;\; \mbox{on}\;\;
 \Gamma_j\cup \Gamma_j^{l_n},\, \; \mbox{and}\;\; u=0\;\;\mbox{on}\;\;
 \partial D_j-\Gamma_j
 \},\\(j=1,\cdots, p)\qquad \qquad \qquad \qquad \qquad \qquad \qquad \qquad \quad \end{eqnarray*}
 with the inner products
\begin{eqnarray*} & \langle u, v\rangle =
  \int_{D} (\triangle u)(\triangle v) dR  +\int_{\Gamma_\varrho} \varrho^3 u v\,ds \quad \; \mbox{for}\;\; u,
  v\in K,\\
&\langle u_j, v_j\rangle_j =\langle u_j, v_j \rangle_j^{\star} +[u_j
, v_j]_j  =
  \int_{D_j} (\triangle u_j)(\triangle v_j)dR +\int_{\Gamma_j} \varrho^3 u_j v_j\,ds
  \quad \; \mbox{for} \;\; u_j, v_j\in K_j^0,\end{eqnarray*}
   respectively.  Closing $K$ and $K_j^0$
   with respect to the norms $\|u\| =
  \sqrt{ \langle u, u\rangle}$ and $\|u_j\|_j =
  \sqrt{ \langle u_j, u_j\rangle_j}$, we obtain the
   Hilbert spaces ${\mathcal{K}}$ and  ${\mathcal{K}}_j^0\,$ ($j=1, \cdots, p$), respectively.
 Clearly, the bilinear functional
\begin{eqnarray*}  &&[u,v] = \int_{\Gamma_\varrho} \varrho^3 \,
  u v\, ds\\
 &&  [u_j, v_j]_j = \int_{\Gamma_j} \varrho^3 \,u_jv_j\, ds, \quad
\;\; (j=1, \cdots, p), \end{eqnarray*}
  define  self-adjoint, completely continuous transformations $G$ and $G^0_j$ on
  ${\mathcal{K}}$ and ${\mathcal{K}}_j^0$
  by \begin{eqnarray} \label{5-2-1}
   \langle Gu, v\rangle&=&[u,v] \quad \;\; \mbox{for}\;\; u,v
\;\;\mbox{in}\;\; {\mathcal{K}},\\
  \label {5-2/1} {\langle G_j^0 u_j, v_j\rangle}_j &=&[u_j,v_j]_j
  \quad \;\; \mbox{for}\;\;
  u_j \,\,\mbox{and}\;\; v_j\,\, \mbox{in}\;\; {\mathcal{K}}_j^0,\end{eqnarray}
 respectively.
 By defining a space \begin{eqnarray*} {\mathcal{K}}^0
 =\sum_{j=1}^p \oplus {\mathcal{K}}_j^0 =\{u^0 \big| u^0=
u_1+\cdots + u_p, \,\, u_j\in {\mathcal{K}}_j^0\}\end{eqnarray*}
    with its inner product
   \begin{eqnarray}  \label{5-2-2}\langle u^0, v^0\rangle =\sum_{j=1}^p
   \langle u_j, v_j\rangle_j, \;\;
   \end{eqnarray}
 we find that the space ${\mathcal{K}}^0$ becomes a Hilbert space.
 If  we define the transformation $G^0$ on ${\mathcal{K}}^0$ by
 \begin{eqnarray} \label{5-2/2} G^0 u^0=G_1^0 u_1+ \cdots+ G^0_p u_p\quad \;\; \mbox{for}\;\;
 u^0 =u_1+\cdots+ u_p\,\, \mbox{in}\,\,
  {\mathcal{K}}^0,\end{eqnarray}
  we see  that $G^0$ is a self-adjoint, completely continuous
  transformation on ${\mathcal{K}}^0$. If we put
  \begin{eqnarray}\label {5-2-3} [u^0, v^0]= \sum_{j=1}^p [u_j,v_j]_j, \end{eqnarray}
  we find by (\ref{5-2/1})---(\ref{5-2-3}) that
  \begin{eqnarray} \label{5-2-4} \langle G^0 u^0, v^0\rangle
  = [ u^0, v^0] \quad \; \,\mbox{for all}\;\;
  u^0\;\;\mbox{and}\;\; v^0 \;\; \mbox{in}\;\; {\mathcal{K}}^0.\end{eqnarray}
   We now define a mapping of ${\mathcal{K}}^0$ into $\mathcal{K}$.
    Let $u^0=u_1+ \cdots+ u_p, \, u\in H_j^0$, be an element of ${\mathcal{K}}^0$
     and define
     \begin{eqnarray} \label{5-2-5} u=\Pi^0 u^0,\end{eqnarray}
     where $u(x)= u_j(x)$, when $x\in \bar D_j$, and $u(x)=0$,
     when $x\in \bar D -\cup_{j=1}^p \bar D_j$. Then $u\in \mathcal{K}$
     and thus (\ref{5-2-5}) defines a transformation $\Pi^0$ of ${\mathcal{K}}_1^0 \oplus \cdots \oplus
     {\mathcal{K}}_p^0$ into
     $\mathcal{K}$.
  It is readily seen that
\begin{eqnarray}\label{5-2-6}
      [\Pi^0 u^0, \Pi^0 v^0] =[u^0, v^0] \quad \; \mbox{for all}\;\;
      u^0 \;\; \mbox{and}\;\; v^0 \;\; \mbox{in}\;\; {\mathcal{K}}^0.
      \end{eqnarray}
and
\begin{eqnarray}
      \label{5-2-7}\langle G^0 u^0,  v^0\rangle
      =\langle G\Pi^0 u^0, \Pi^0 v^0\rangle  \quad \; \mbox{for all}\;\;
      u^0 \;\; \mbox{and}\;\; v^0 \;\; \mbox{in}\;\; {\mathcal{K}}^0.
      \end{eqnarray}
 From (\ref{5-2-6}) and (\ref{5-2-7}), we find by applying Corollary 1.4.1 of \cite{Sa} that
 \begin{eqnarray*} \mu_k^0 \le \mu_k \quad \; \mbox{for}\;\; k=1,2,3, \cdots\end{eqnarray*}
 Therefore \begin{eqnarray} \label {5-38} A^0(\tau)\le A(\tau).\end{eqnarray}
 The definition of $G^0$ implies that
 \begin{eqnarray} \label{5-2-8}  G^0 {\mathcal{K}}_j^0
 \subset {\mathcal{K}}_j^0, \quad \; (j=1, \cdots, p),\end{eqnarray}
  and
\begin{eqnarray} \label{5-2-9} G^0 u^0 = G_j^0 u^0,
 \quad \, \mbox{when}\;\; u^0\in  {\mathcal{K}}_j^0.\end{eqnarray}
 From (\ref{5-2-3}), (\ref{5-2-4}), (\ref{5-2-8}), (\ref{5-2-9}) and Theorem 1.6 of \cite{Sa}, we obtain
\begin{eqnarray} \label {5-40} A^0(\tau) =\sum_{j=1}^p A^0_j (\tau),\end{eqnarray}
  where $A_j^0 (\tau)$ is the number of eigenvalues
  of the transformation
   $G_j^0$ on ${\mathcal{K}}_j^0$
 which are greater or equal to $(1+\tau^3)^{-1}$.
 Because  $\bar D_j, (j=1, \cdots, p)$, is an $n$-dimensional
    rectangular parallelepiped
   we find  by (\ref{5-30}) that
\begin{eqnarray} \label{5-2-10} A^0_j (\tau)  \sim \omega_{n-1}(\sqrt[3]{16}\,\pi)^{-(n-1)}
 |\Gamma_j| \varrho^{n-1}  \tau^{n-1}\quad \;\; \mbox{as}\;\; \tau\to +\infty,\end{eqnarray}
where $|\Gamma_j|$ denotes the area of the face $\Gamma_j$ of $D_j$.
By (\ref{5-40}) and (\ref{5-2-10}) we infer that
\begin{eqnarray} \label{5-2-11} A^0 (\tau) \sim \omega_{n-1}(\sqrt[3]{16}\,\pi)^{-(n-1)}
  \sum_{j=1}^p |\Gamma_j| \varrho^{n-1}  \tau^{n-1}\quad \, \mbox{as}\;\; \tau\to +\infty.\end{eqnarray}

\vskip 0.12 true cm

 Next, we shall calculate the upper estimate of $A(\tau)$.
  Let $\bar P_j, (j=p+1, \cdots,q)$, be the $n$-dimensional rectangular parallelepiped
   with the cube $\bar Q_j$ as a base and otherwise
 bounded by the ``upper'' plane surface $\Gamma^{l_n}$ of the cylinder $\bar D$ and planes parallel to the
 coordinate-planes $x_1=0, \cdots, x_{n-1}=0$.
  The intersection $\bar P_j \cap \bar D$ is a cylinder $\bar D_j, (j=p+1, \cdots, q)$,
  with $\bar \Gamma_j:=\bar Q_j \cap \bar \Gamma_\varrho$  as a basis. Then
  \begin{eqnarray} \bar D=\sum_{j=1}^q \bar D_j.\end{eqnarray}

  We first define the linear spaces of functions
  \begin{eqnarray*}& K^d =\{u\big| u\in  Lip (\bar D)\cap H^2(D),
  \; \frac{\partial u}{\partial \nu}=0\,\, \mbox{on} \;\; \Gamma_\varrho, \; u=\frac{\partial u}{\partial \nu}
  =0 \,\, \mbox{on}\,\,
   \Gamma^{l_n})\},\\
 & K_j^d =\{u\big| u\in  Lip (\bar D)\cap H^2(D_j),
  \; \frac{\partial u}{\partial \nu}=0\,\, \mbox{on} \;\; \Gamma_j, \,\,
   u=\frac{\partial u}{\partial \nu}
  =0 \,\, \mbox{on}\,\,
   \Gamma^{l_n}_j\},\; \; (j=1, \cdots, q)\end{eqnarray*}
  with the inner products
    \begin{eqnarray} \label{5-3-1}  & \langle u, v\rangle = \int_{D} (\triangle u)(\triangle v) dx
  +\int_{\Gamma_\varrho}  \varrho^3  u v \, ds, \end{eqnarray}
  and
  \begin{eqnarray} \label {5-3.2}  \langle u_j, v_j\rangle_j = \int_{D_j}  (\triangle u_j)(\triangle v_j) dx
  +\int_{\Gamma_j}  \varrho^3u_j v_j  \, ds,\end{eqnarray}
  respectively.
    Closing $K^d$ and $K_j^d$ with respect to the norms $\|u\|=\sqrt{\langle u, u\rangle}$
  and $\|u_j\|_j
  =\sqrt{\langle u_j, u_j\rangle_j}$,
 we get  Hilbert spaces
 ${\mathcal{K}}^d$ and ${\mathcal{K}}_j^d, (j=1, \cdots, q)$, and then we define the Hilbert
 space
\begin{eqnarray} \label{5-3-2} {\mathcal{K}}^d = \sum_{j=1}^q \oplus {\mathcal{K}}^d_j=\{
  u^d\big| u^d =u_1+\cdots+ u_q, \; \; u_j \in
{\mathcal{K}}_j^d\}\end{eqnarray}
 with its inner product
  \begin{eqnarray} \label{5-3-3}\langle u^d,v^d \rangle
   = \sum_{j=1}^q \langle u_j, v_j\rangle_j.    \end{eqnarray}
  The bilinear functional
\begin{eqnarray} \label{5-3-4} [u_j, v_j]_j = \int_{\Gamma_j}  \varrho^3  u_jv_j\, ds, \quad \;\; (j=1, \cdots,
q),\end{eqnarray}
 define a self-adjoint, completely continuous transformation $G_j^d$
 on ${\mathcal{K}}_j^d$ given by
\begin{eqnarray} \label{5-3-5}
      \langle G^d_j u_j,  v_j\rangle_j =[u_j, v_j]_j \quad \; \mbox{for all}\;\;
      u_j \;\; \mbox{and}\;\; v_j \;\; \mbox{in}\;\; {\mathcal{K}}_j^d.
      \end{eqnarray}
       The self-adjoint, completely continuous transformation $G^d$ on ${\mathcal{K}}^d$ is defined by
 \begin{eqnarray} \label{5-3-6}G^d u^d=G_1^d u_1+\cdots+ G^d_q u_q\quad \;\; \mbox{for}\;\;
 u^d =u_1+\cdots+ u_q\,\, \mbox{in}\,\,
  {\mathcal{K}}^d.\end{eqnarray}
  With
    \begin{eqnarray} \label{5-3-7} [u^d, v^d]= \sum_{j=1}^q [u_j,v_j]_j, \end{eqnarray}
  it follows from  (\ref{5-3-3}), (\ref{5-3-5})---(\ref{5-3-7}) that
  \begin{eqnarray} \label{5-3-8} \langle G^d u^d, v^d\rangle = [u^d, v^d] \quad \; \,\mbox{for all}\;\;
  u^d \,\, \mbox{and}\,\, v^d \;\; \mbox{in}\;\; {\mathcal{K}}^d.\end{eqnarray}
 Now we define a mapping  $\Pi$ of $\mathcal{K}$ into ${\mathcal{K}}^d$.
 Let $u\in K(D)$, and put
 \begin{eqnarray*} u^d =\Pi u=u_1+ \cdots + u_q,\end{eqnarray*}
where $u_j(x)= u(x)$, when $x\in \bar D_j$. It can be easily
 verified that
  \begin{eqnarray} \label{5-3-9} \langle \Pi u, \Pi v\rangle =\langle u, v\rangle
      \quad \, \mbox{for all}\;\; u \;\; \mbox{and}\;\; v \;\; \mbox{in}\;\;
      {\mathcal{K}}. \end{eqnarray}
and
 \begin{eqnarray} \label{5-3-10} \langle Gu, v\rangle =\langle G^d \Pi u, \Pi v\rangle
        \quad \, \mbox{for all}\;\; u \;\; \mbox{and}\;\; v \;\; \mbox{in}\;\;
            \mathcal{K}.\end{eqnarray}
 Combining (\ref{5-3-9}), (\ref{5-3-10}) and using Corollary 1.4.1 of \cite{Sa}, we obtain
  \begin{eqnarray*} \mu_k \le \mu_k^d \quad \; \mbox{for}\;\; k=1,2,3, \cdots,\end{eqnarray*}
 and hence
  \begin{eqnarray} \label {5..56} A(\tau)\le A^d(\tau).\end{eqnarray}
 From $G^d {\mathcal{K}}_j^d\subset {\mathcal{K}}_j^d$, $\;(j=1, \cdots,
 q)$, and $G^d u^d =G^d_j u^d$ when $u^d\in {\mathcal{K}}_j^d$,  we
 get
\begin{eqnarray} A^d(\tau) =\sum_{j=1}^q A^d_j (\tau),\end{eqnarray}
  where $A_j^d (\tau)$ is the number of eigenvalues
   of the transformation
   $G_j^d$ on ${\mathcal{K}}_j^d$
   which are
   greater than or equal to $(1+\tau^3)^{-1}$.
Also, we define $A_j^f (\tau)$ similar to (\ref{5.*.15}) and
(\ref{5..0}), i.e.,
 \begin{eqnarray*}  A^f_j(\tau)
   = \sum_{\mu_k^f\ge (1+\tau^3)^{-1}} 1\; \quad \;\mbox{with}\;\;
   \mu_k^f= \frac{1}{1+\gamma_k^3},\end{eqnarray*}
 where  $\gamma_k^3$ is  the $k$-th Steklov eigenvalue of the following problem
  \begin{eqnarray*} \label {5..15}  \left\{ \begin{array}{ll} \triangle^2 u_j = 0\quad \;
   \mbox{in}\;\; D_j,\\
   \frac{\partial u_j}{\partial \nu}=0 \; \; \mbox{on}\;\; \Gamma_j, \quad \,
   u_j=\frac{\partial u_j}{\partial \nu}=0 \;\;\mbox{on}\;\;
   \Gamma^{l_n}_j,\\
      \frac{\partial u}{\partial \nu}=\frac{\partial (\Delta u_j)}{\partial \nu}= 0
      \;\; \mbox{on}\;\; \partial D_j-(\Gamma_j\cup \Gamma^{l_n}_j),\\
    \frac{\partial (\triangle u_j)}{\partial \nu} -\gamma^3 \varrho^3  u_j=0 \;\;
    \mbox{on}\;\; \Gamma_j, \quad \,  \varrho =constant>0\;\;\mbox{on}\;\;
     \Gamma_\varrho^+.
   \end{array}\right.\end{eqnarray*}
    From Theorem 3.8, it follows that
 \begin{eqnarray*} \gamma_k^3\le \kappa_k^3\quad \;\; \mbox{for all }\;\;
 k\ge 1,\end{eqnarray*}
 and hence
 \begin{eqnarray}  \label{5---1} A^d_j(\tau)\le A^f_j(\tau)\quad \;\mbox{for all}\;\;
   \tau\;\;\mbox{and}\;\;j=1,\cdots, p,\end{eqnarray}
  where $\kappa_k^3$ is the $k$-th eigenvalue of the transformation
 $G^d_j$.
   Since $\bar D_j, (j=1, \cdots, p)$, is an $n$-dimensional rectangular
   parallelepiped,
  we find from
 (\ref{5.2.8}) that
 \begin{eqnarray}\label{5..58}  A^f_j (\tau) \sim \omega_{n-1}(\sqrt[3]{16}\,\pi)^{-(n-1)}
 |\Gamma_j| \varrho^{n-1}  \tau^{n-1}, \quad \;(j=1, \cdots, p). \end{eqnarray}

 It remains to estimate $A^f_j(\tau), (j\ge p+1)$.
 According to the argument in p.$\,$438-440 of \cite{CH},
  each of the $(n-1)$-dimensional domains $\Gamma_j$ is
 bounded either by $n-1$ orthogonal plane surfaces of the partition
 (the diameter of the intersection  of any two plane surfaces lies
   between $l$ and $3l$), and an $(n-2)$-dimensional
  surface of the boundary $\partial \Gamma_\varrho$
  (see, in two dimensional case, Figure 5 of p.$\,$439 of \cite{CH}),
 or by $2n-3$ orthogonal plane surfaces of the partition (the diameter of the intersection of any
 two plane surfaces lies between $l$ and $3l$),
 and a surface of the boundary $\partial \Gamma_\varrho$ (see, in two dimensional case,
   Figure 6 of p.$\,$439 of \cite{CH}).
The number $q-p$ is evidently smaller than a constant $C/l^{n-2}$,
where $C$ is independent of $l$ and depends essentially on the area
of
 the boundary $\partial \Gamma_\varrho$.
  Now, we take any point on the  boundary surface of $\Gamma_j$
  and take the tangent plane through it.
   This tangent plane together with the plane parts of $\partial \Gamma_j$
   bounds an $n$-polyhedron of ${\Bbb R}^{n-1}$ with a vertex at which $n-1$
   orthogonal plane surfaces meet  (see,
    Figure 5 of p.$\,$439 of \cite{CH} in two dimensions),
   e.g., if $\vartheta$ is sufficiently small it forms an $(n-1)$-dimensional
   $n$-polyhedron of ${\Bbb R}^n$
   with a vertex having $n-1$ orthogonal plane surfaces (the diameter of the intersection  of any two
   plane surfaces is also smaller than $4l$), or else an  $(n-1)$-dimensional $2(n-1)$-polyhedron
   of ${\Bbb R}^{n-1}$ (see, Figure 6 of p.$\,$439 of \cite{CH} in two dimensional case),
      the diameter of the intersection of any two plane surfaces (except for the top inclined plane surface)
      of the $2(n-1)$-polyhedron is also smaller than $4l$;
     The shape of the result domain depends on the type to which  $\bar \Gamma_j$ belongs.
   We shall denote the result domains by $S'_j$.  The domain $\Gamma_j$
   can always be deformed into the domain $S'_j$ by a transformation
   of the form (\ref{02-001}), as defined in Section 2.  In the case of
   domains of the first  type,  let the intersection point of $n-1$ orthogonal plane surfaces
    be the pole of a system  of pole coordinates
   $r$, $\theta_1$,  $\theta_2, \cdots, \theta_{n-2}$, and let $r=f(\theta_1,
   \theta_2, \cdots, \theta_{n-2})$ be the equation of the boundary
   surface of $\Gamma_\varrho$, $r=h(\theta_1,
   \theta_2, \cdots, \theta_{n-2})$
  the equation of the inclined plane surface of the $n$-polyhedron of
  ${\Bbb R}^{n-1}$ having a vertex
  of $n-1$ orthogonal plane surfaces. Then the
  equations
   $$\theta_1'=\theta_1, \quad  \theta'_2=\theta_2, \;\cdots,\;\,
      \theta'_{n-2}=\theta_ {n-2},  \quad r'= r\,\frac{h(\theta_1,
   \theta_2, \cdots, \theta_{n-2})}{f(\theta_1,
   \theta_2, \cdots, \theta_{n-2})}$$
  represents a transformation of the domain  $\Gamma_j$ into
  the $n$-polyhedron  $S'_j$ of ${\Bbb R}^{n-1}$.
     For a domain of the second type, let $x_{n-1}=
  h (x_1, \cdots, x_{n-2})$  be the equation of top plane surface of the
  $2(n-1)$-polyhedron  and let $x_{n-1}=f(x_1, \cdots, x_{n-2})$
  be the  equation of the boundary  surface of $\Gamma_\varrho$.
  We then consider the  transformation
  $$ x'_1=x_1, \;\; \cdots, \quad x'_{n-2}=x_{n-2}, \quad
   x'_{n-1}= x_{n-1}\,\frac{h(x_1, \cdots, x_{n-2})}
   {f(x_1, \cdots, x_{n-2})}. $$
If we assume that the side length $l$ of cube in the partition is
sufficiently small, and therefore the rotation of the normal on the
 boundary surface  is taken sufficiently small, then the
transformations considered here evidently have precise the form
(\ref{02-001}),  and the quantity denoted by $\epsilon$ in
(\ref{02-001}) is arbitrarily small.
  From Corollary to Theorem 10 of p.$\,$423 of \cite{CH}, we know that
  there exists a number $\delta>0$ depending on $\epsilon$
 and approaching zero with $\epsilon$, such that
  \begin{eqnarray*} \bigg| \frac{\alpha_k(S'_j)}{\alpha_k(\Gamma_j)}
  -1\bigg|<\delta \quad \, \mbox{uniformly for all}\;\; k,
  \end{eqnarray*}
  where $\alpha_k(\Gamma_j)$ and $\alpha_k (S'_j)$  are the $k$-th Neumann
  eigenvalues of $\Gamma_j$ and $S'_j$, respectively.
  According to the argument as in the proof of Lemma 5.1, we see
  that
  \begin{eqnarray*}
 (\gamma_k^f (E_j))^3=\frac{1}{\varrho^3 l_n^3}\, t(
  l_n \, \alpha_k(\Gamma_j)), \quad \;  \big(\gamma_k^f (E'_j)\big)^3=\frac{1}{\varrho^3 l_n^3}\, t(
 l_n \, \alpha_k(S'_j)),\end{eqnarray*}
 where  $t(s)$ is given by (\ref{005.10}), and $(\gamma_k^f (E_j))^3$
  and $\big(\gamma_k^f(E'_j)\big)^3$ (similar to $\gamma^3$ of (\ref{3-30-0}))
  are the $k$-th Steklov eigenvalue
   for the $n$-dimensional  domains
  $E_j=\Gamma_j\times [0, l_n]$
  and $E'_j= S'_j \times [0, l_n]$, respectively.  Recalling
  that the function
  $t=t(s)$ is continuous and increasing for $s\ge 1$,
   and $\lim_{s\to +\infty} \frac{t(s)}{s}=\frac{1}{2}$,
      we get
 that there exists a constant $\delta'>0$ depending on $\epsilon$
 approaching
 zero with $\epsilon$, such that
\begin{eqnarray*} \bigg|\frac{\big(\gamma_k^f (E_j)\big)^3}{\big(\gamma_k^f
(E'_j)\big)^3}-1\bigg|<\delta'.\end{eqnarray*}
 In other words, the corresponding
$k$-th eigenvalues for the $n$-dimensional
 domains $E_j=\Gamma_j\times [0, l_n]$ and $E'_j=S'_j\times
[0,l_n]$ differ only by a  factor which itself differs by a small
amount from $1$, uniformly  for all $k$. Therefore, the same is true
also for the corresponding
 numbers $A^f_{E_j} (\tau)$ and $A^f_{E'_j} (\tau)$ of
  the eigenvalues less or equal to the bound
 $\tau^3$.

The  domain $E'_j$ is either a cylinder whose base is an
$n$-polyhedron of ${\Bbb R}^{n-1}$ having $(n-1)$ orthogonal plane
surfaces with its largest side length small than $4l$ or a cylinder
whose base is a
 combination of such an $n$-polyhedron of ${\Bbb R}^{n-1}$
 and an $(n-1)$-dimensional  cube with
side-length smaller than $3l$; it follows that if $l$ is taken
sufficiently small,  the number $A^f_{E_j}(\tau)$ from some $\tau$
on satisfies the inequality
  \begin{eqnarray*} \label{07-00}  A^f_{E_j}(\tau) < C_1 l^{n-1}
  \tau^{n-1} +C_2 l^{n-2} \tau^{n-2} \end{eqnarray*}
where $C_1$, $C_2$ are constants, to be chosen suitably.
  Thus, $A_{E_j}^f (\tau)$ can be  written as
  $A_{E_j}^f (\tau)=\theta (C_3 l^{n-1} \tau^{n-1}+C_4
  l^{n-2}\tau^{n-2})$, where $\theta$ denotes a number between $-1$
  and $+1$ and $C_3, C_4$ are constants independent of $l,j$ and
  $\tau$. It follows that $$\sum_{j=p+1}^q A_{E_j}^f(\tau)=\tau^{n-1}
  \big[ \theta C_3 (q-p)l^{n-1} +\theta C_4 (q-p)l^{n-2}
  \,\frac{1}{\tau}\big].$$
  As pointed out before, $(q-p)l^{n-2}<C$; therefore, for sufficiently small
  $l$, $\,(q-p)l^{n-1}$ is arbitrarily
 small and we have  the asymptotic relation
 \begin{eqnarray} \label{5.-01} \lim_{\tau\to +\infty}
 \sum_{j=p+1}^q \frac{A_{E_j}^f (\tau)}{\tau^{n-1}}
 =\varsigma(l),  \end{eqnarray}
 where $\varsigma(l)\to 0$ as $l\to 0$.
 For, we may choose the quantity $l$ arbitrarily, and  by taking
   a sufficiently small fixed $l$, make the factor of $\tau^{n-1}$ in
  the above equalities arbitrarily close to zero for sufficiently
  large $\tau$.
 Since \begin{eqnarray}   \label{5.-02}  A^d_{E_j}(\tau)\le A^f_{E_j}(\tau)
 \quad \, \mbox{for}\;\; j=p+1, \cdots, q,\end{eqnarray}
 we get
 \begin{eqnarray} \label{5.-03} \lim_{\tau\to +\infty}
 \frac{ \sum_{j=p+1}^q A_{E_j}^d (\tau)}{\tau^{n-1}}\le
 \lim_{\tau\to +\infty}
 \frac{ \sum_{j=p+1}^q A_{E_j}^f (\tau)}{\tau^{n-1}}=\varsigma(l).  \end{eqnarray}
  Combining (\ref{5-38}), (\ref{5-2-11}), (\ref{5..56}), (\ref{5---1}),
  (\ref{5..58}), (\ref{5.-01}),  (\ref{5.-02}) and
  (\ref{5.-03}) we obtain
  \begin{eqnarray} \label {5.63} &&
 \omega_{n-1}(\sqrt[3]{16}\,\pi)^{-(n-1)} \varrho^{n-1}
  \sum_{j=1}^p |\Gamma_j| \le
 \underset {\tau\to \infty} {\underline{\lim}}\, \frac{A(\tau)}{\tau^{n-1}}
\le \overline{\lim_{\tau\to \infty}} \,\frac{A(\tau)}{\tau^{n-1}}
\\
 && \;\;\quad \qquad \quad\;\;  \le
 \bigg(\omega_{n-1}(\sqrt[3]{16}\,\pi)^{-(n-1)}
   \varrho^{n-1}\sum_{j=1}^p |\Gamma_j|\bigg)+\varsigma(l). \nonumber
   \end{eqnarray}
 Letting $l\to 0$, we immediately see that
$\sum_{j=1}^p |\Gamma_j|$ tends to the area $|\Gamma_\varrho|$ of
 $\Gamma_\varrho$ and $\lim_{l\to 0} \varsigma(l)=0$. Therefore,
 (\ref{5.63}) gives
\begin{eqnarray} A (\tau) \sim
\frac{\omega_{n-1}}{(\sqrt[3]{16}\,\pi)^{(n-1)}}
 |\Gamma_\varrho| \varrho^{n-1}  \tau^{n-1} \quad \; \mbox{as}\;\; \tau\to +\infty,\end{eqnarray}
or
\begin{eqnarray}\label{5-90} A (\tau) \sim \frac{\omega_{n-1}\tau^{n-1}}{(\sqrt[3]{16}\,\pi)^{(n-1)}}
   \int_{\Gamma_\varrho}  \varrho^{n-1} ds \quad \mbox{as} \;\; \tau\to +\infty.\end{eqnarray}

In the above argument, we first made the assumption that the
boundary $\partial \Gamma_\varrho$ of $\Gamma_\varrho$  is smooth.
However,  the corresponding discussion and result remain essentially
valid if $\partial \Gamma_\varrho$ is composed of a finite number of
$(n-2)$ dimensional smooth surfaces.

 \vskip 1.39 true cm

\section{Proofs of main results}

\vskip 0.45 true cm

 \noindent  {\bf Lemma 6.1.} \ \  {\it  Let $g^{il}$ and $g'^{il}$
  be two metric tensors on manifold $\mathcal{M}$ such that
    \begin{eqnarray} \label {006-1}
 \big|g^{il} -g'^{il}\big|< \epsilon, \quad \, i,l=1, \cdots, n
    \end{eqnarray}
and
 \begin{eqnarray} \label {006-2} \bigg| \frac{1}{\sqrt{|g|}} \,\frac{\partial}{\partial x_i}
   \big(\sqrt{|g|} g^{il}\big)\,- \, \frac{1}{\sqrt{|g'|}} \,\frac{\partial}{\partial x_i}
   \big(\sqrt{|g'|} g'^{il}\big)\bigg|\le \epsilon, \quad \,
   i,l=1,\cdots,n   \end{eqnarray}
   for all points in $\bar D$, where $D$ is a bounded domain in $\mathcal{M}$ (see, Section 3).
    Let \begin{eqnarray*} \mu_1\ge \mu_2\ge \cdots \ge \mu_n \ge
  \cdots >0 \;\; \mbox{and}\;\;  \mu'_1\ge \mu'_2\ge \cdots \ge \mu'_n \ge
  \cdots >0\end{eqnarray*}
  be positive eigenvalues of $G$ and $G'$, respectively,
  where $G$ and $G'$ are given by \begin{eqnarray*} \langle Gu,
v\rangle =\int_{\Gamma_\varrho} \varrho^3 \, uv\, ds, \quad\;
\mbox{for}\;\; u \,\, \mbox{and}\;\; v
\,\, \mbox{in}\;\; \mathcal{K}, \\
\langle G'u, v\rangle' =\int_{\Gamma_\varrho} \varrho^3 \, uv\, ds',
\quad\; \mbox{for}\;\; u \,\, \mbox{and}\;\; v \,\, \mbox{in}\;\;
{\mathcal{K}}'.\end{eqnarray*}
  Then, for $k=1,2,3, \cdots$,
  \begin{eqnarray}\label {006-3} &&
\frac{(1+\epsilon{\tilde M})^{-n-\frac{1}{2}}}{
 \big(\max \{ (1+\epsilon M), (1+\epsilon{\tilde M})^{1/2}\}\big) }
   \mu_k \le \mu'_k
 \\
  && \quad \quad \quad\;\le
\frac{(1+\epsilon{\tilde M})^{n+\frac{1}{2}}}
 {\big(\min \{(1-\epsilon M), (1+\epsilon{\tilde M})^{-1/2}\}\big) }
 \mu_k,  \nonumber \end{eqnarray}
where $\tilde M$ and $M$ are constants depending only on $g$, $g'$,
$\partial g$, $\partial g'$  and $\bar D$.}

 \vskip 0.24 true cm

\noindent  {\bf Proof.} \
 It follows from (\ref{006-1})  that
 there exists a positive constant $\tilde M$ independent of
 $\epsilon$ and depending only on $g^{ij}$, $g'^{ij}$ and $\bar D$ such that
 \begin{eqnarray*} (1+\epsilon \tilde M)^{-1} \sum_{i,l=1}^{n}
 g^{il} t_it_l \le \sum_{i,l=1}^{n}
 g'^{il} t_it_l \le (1+\epsilon \tilde M) \sum_{i,l=1}^{n}
 g^{il} t_it_l  \end{eqnarray*}
 for all points in $\bar D$ and all real numbers $t_1, \cdots, t_n$.
 Thus we have
 \begin{eqnarray*} (1+\epsilon{\tilde M})^{-n/2} \sqrt{|g|}\le \sqrt{|g'|}
\le (1+\epsilon{\tilde M})^{n/2} \sqrt{|g|}, \end{eqnarray*}
   which implies (see,
  p.$\;$64-65 of \cite{Sa}) that \begin{eqnarray} \label{-0-6.1} (1+\epsilon{\tilde M})^{-n/2} dR \le dR'
  \le (1+\epsilon{\tilde M})^{n/2} dR
\end{eqnarray}
 and \begin{eqnarray} \label{-0-6.2}   (1+\epsilon\tilde{M})^{-(n+1)/2} ds \le ds' \le
 (1+\epsilon{\tilde M})^{(n+1)/2} ds.\end{eqnarray}
 Thus \begin{eqnarray} \label {006-4} (1+\epsilon{\tilde M})^{-(n+1)/2} [u,u]\le [u, u]' \le
 (1+\epsilon{\tilde M})^{(n+1)/2} [u, u].\end{eqnarray}
  Putting   $$\omega_{il}=g'^{il}- g^{il}, \quad \; \theta_{il}
  =\frac{1}{\sqrt{|g'|}} \,\frac{\partial}{\partial x_i}
   \big(\sqrt{|g'|} g'^{il}\big)\,- \, \frac{1}{\sqrt{|g|}} \,\frac{\partial}{\partial x_i}
   \big(\sqrt{|g|} g^{il}\big),$$
  we  immediately see that \begin{eqnarray*}
   \max_{x\in \bar D}|\omega_{il}|\le \epsilon \quad \;
    \mbox{and}\;\; \max_{x\in \bar D} |\theta_{il}|\le
   \epsilon.\end{eqnarray*}
 Thus, for any $u\in {\mathcal{K}} (D)$ or $u\in {\mathcal{K}}^d
 (D)$,
we have \begin{eqnarray*} \quad \triangle_{g'} u
   = \sum_{i,l=1}^n (\omega_{il} +g^{il} ) \frac{\partial^2 u}{\partial
x_i\partial x_l} + \sum_{i,l=1}^n \left[\theta_{il}+
  \frac{1}{\sqrt{|g|}} \,\frac{\partial}{\partial x_i}
   \big(\sqrt{|g|} g^{il}\big)\right] \frac{\partial
   u}{\partial x_l}, \end{eqnarray*}
so that \begin{eqnarray*}\label{006-5}  \triangle_{g'} u -
 \triangle_g u = \sum_{i,l=1}^n \left[ \omega_{il}  \frac{\partial^2 u}{\partial
  x_i\,\partial x_l} +\theta_{il} \,
   \frac{\partial u}{\partial
   x_l}\right]\nonumber\end{eqnarray*}
  It follows that
\begin{eqnarray*} | \triangle_{g'} u -
 \triangle_g u| \le \epsilon \big(M_1|\nabla_g^2 u|+M_2|\nabla_g
 u|\big),\end{eqnarray*}
 where $|\nabla^2_g u|^2$ is defined in an invariant ways as
  \begin{eqnarray*} |\nabla_g^2 u|^2 =\nabla^l \nabla^k u \,
  \nabla_l \nabla_k u = g^{pl} g^{kq} \left(\frac{\partial^2
  u}{\partial x^k \partial x^l} -\Gamma_{kl}^m \frac{\partial
  u}{\partial x^m}\right)\left(\frac{\partial^2
  u}{\partial x^p \partial x^q} -\Gamma_{pq}^r \frac{\partial
  u}{\partial x^r}\right),\end{eqnarray*}
 and $M_1$ and $M_2$ are constants depending only on $g$, $g'$,
 $\partial g$, $\partial g'$ and
 $\bar D$.
 Thus, \begin{eqnarray} \label {006-7}
 \int_{D} | \triangle_{g'} u -
 \triangle_g u|^2 dR  \le 2 \epsilon^2 \left(M_1^2\int_D |\nabla_g^2 u|^2dR+M_2^2\int_D |\nabla_g
 u|^2dR   \right).\end{eqnarray}

Set  \begin{eqnarray} \label {006-10}
  {\tilde \Lambda}_1^0 (D) =\inf_{v\in K(D), \;\int_D |\nabla_g v|^2 dR =1}
 \,\frac{\int_D|\triangle_g v|^2 dR }{\int_D |\nabla_g v|^2
dR},\end{eqnarray}
  \begin{eqnarray} \label {006-11}
  {\tilde\Lambda}_1^d (D) =\inf_{v\in K^d(D), \;\int_D |\nabla_g v|^2 dR =1}
\, \frac{\int_D|\triangle_g v|^2 dR }{\int_D |\nabla_g v|^2
dR},\end{eqnarray}
 and the spaces  $K(D)$ and  $K^d(D)$ are  as in  Section 3.
  Furthermore, set \begin{eqnarray} \label {006-8}
  \Theta_1^0 (D) =\inf_{v\in K(D), \;\int_D |\nabla_g^2 v|^2 dR =1}
 \,\frac{\int_D|\triangle_g v|^2 dR }{\int_D |\nabla_g^2 v|^2
dR},\end{eqnarray}
  \begin{eqnarray} \label {006-9}
  \Theta_1^d (D) =\inf_{v\in K^d(D), \;\int_D |\nabla_g^2 v|^2 dR =1}
 \,\frac{\int_D|\triangle_g v|^2 dR }{\int_D |\nabla_g^2 v|^2
dR}.\end{eqnarray}
  Clearly, ${\tilde \Lambda}_1^0 (D)\ge {\tilde \Lambda}_1^d (D)$, and
   $\Theta_1^0(D)\ge \Theta_1^d
  (D)$.
  Similar to the proofs of Lemmas 2.1, 2.2, it is easy to prove that
  the existence of the minimizers to (\ref{006-11}) and (\ref{006-9}),
  respectively. Therefore, we have that ${\tilde \Lambda}_1^d(D)>0$ and $\Theta_1^d
  (D)>0$ (Suppose by contradiction that ${\tilde \Lambda}_1^d(D)=0$ and $\Theta_1^d
  (D)=0$.
  Then $\triangle_g u=0$ in $D$ for the corresponding minimizer $u\in K^d(D)$ in two cases.
 By applying
   Holmgren's uniqueness theorem for the minimizer $u\in K^d(D)$
   in each case, we immediately see that
   $u\equiv 0$ in $D$.
 This contradicts the assumption $\int_D |\nabla_g u|^2
 dR=1$ or $\int_D |\nabla_g^2 u|^2 dR =1$ for the minimizer $u\in K^d(D)$ in the corresponding cases).
  Combining these inequalities, we obtain
\begin{eqnarray*} \label {006-12}
 \int_{D} | \triangle_{g'} u -
 \triangle_g u|^2 dR  \le 2\epsilon^2 \left(\frac{M^2_1}{\Theta_1^0(D)}  +\frac{M^2_2}{{\tilde \Lambda}_1^0
 (D)}\right)
  \int_D |\triangle_g
 u|^2, \quad \mbox{for}\;\; u\in K(D)\end{eqnarray*}
and \begin{eqnarray*} \label {006-13}
 \int_{D} | \triangle_{g'} u -
 \triangle_g u|^2 dR  \le 2 \epsilon^2 \left(\frac{M^2_1}{\Theta_1^d(D)}  +\frac{M^2_2}{{\tilde \Lambda}_1^d (D)}
  \right) \int_D |\triangle_g
 u|^2 \quad \mbox{for}\;\; u\in K^d(D).\end{eqnarray*}
  Thus we have  that, for all $u\in K(D)$ or $u\in
 K^d(D)$,
 \begin{eqnarray*} \label {006-14}
    (1-\epsilon M)\int_D | \triangle_g u|^2 dR
  \le    \int_{D} |\triangle_{g'} u|^2 dR
   \le (1+\epsilon M) \int_D |\triangle_{g} u|^2 dR,\end{eqnarray*}
where $M$ is a constant depending only $g$, $g'$, $\partial g$,
 $\partial g'$ and $\bar D$.
 Combining this and (\ref{-0-6.1}) we get
   \begin{eqnarray*}
(1-\epsilon M) (1+\epsilon \tilde
  M)^{-\frac{n}{2}}\int_D |\Delta_g u|^2 dR\le
  \int_D |\Delta_{g'} u|^2 dR' \le (1+\epsilon M) (1+\epsilon \tilde
  M)^{\frac{n}{2}}\int_D |\Delta_g u|^2 dR.
   \end{eqnarray*}
 That is,
\begin{eqnarray} \label {006-15} (1-\epsilon M)(1+\epsilon \tilde
  M)^{-\frac{n}{2}} \langle u,u\rangle^\star \le \langle u, u\rangle'^\star \le
 (1+\epsilon M) (1+\epsilon \tilde
  M)^{\frac{n}{2}} \langle u, u\rangle^\star.\end{eqnarray}
By (\ref{006-4}) and (\ref{006-15}) we obtain that, for all $u\in K
(D)$ or $u\in K^d(D)$,
\begin{eqnarray*} \label {006-16} &&\frac{(1+\epsilon{\tilde M})^{-(n+1)/2}[u,u]}{
 \big(\max \{ (1+\epsilon M)(1+\epsilon \tilde
  M)^{\frac{n}{2}}, (1+\epsilon{\tilde M})^{(n+1)/2}\}\big) \big(
  \langle u, u\rangle^* +[u,u]\big)} \le
   \frac{[u,u]'}{\langle u,u\rangle'^\star
   +[u,u]'}\\
&&\quad \;\;\le
 \frac{(1+\epsilon{\tilde M})^{(n+1)/2} [u,u]}
 {\big(\min \{(1-\epsilon M)(1+\epsilon \tilde
  M)^{-\frac{n}{2}}, (1+\epsilon{\tilde M})^{-(n+1)/2}\}\big) \big( \langle u,
 u\rangle^\star  +[u,u]\big)},
 \end{eqnarray*}
 which implies (\ref{006-3}). \ \  $\square$

\vskip 0.23 true cm
 \noindent  {\bf Remark 6.2.} \  Let $\tilde \Gamma$ and $\Gamma$ be two bounded domains
 in ${\Bbb R}^{n-1}$, and let $\tilde \Gamma$ is similar to $\Gamma$ (in
 the  elementary sense of the term; the length of any line in $\tilde \Gamma$
  is to the corresponding length in $\Gamma$ as $h$  to $1$), and
  let $\Gamma_0 =\Gamma\times \{\sigma\}$ and $\tilde \Gamma_0 =\tilde \Gamma\times \{h\sigma \}$.
   It is easy to verify that \begin{eqnarray*} {\tilde\Lambda}_1^d(\tilde D) = h^{-2} {\tilde \Lambda}_1^d(D), \quad
   \;\, \Theta_1^d(\tilde D)=\Theta_1^d (D), \end{eqnarray*}
where $D=\Gamma\times [0, \sigma]$, $\tilde D=\tilde \Gamma \times
[0, h\sigma]$, and ${\tilde \Lambda}_1^d (D)$ and $\Theta_1^d(D)$
are defined as in (\ref{006-11}) and
 (\ref{006-9}), respectively.

\vskip 0.26 true cm

 \noindent  {\bf Lemma 6.3.} \ \  {\it
  Let $G$ and  $G'$ be the continuous linear transformations defined by
\begin{eqnarray*} \langle Gu, v\rangle = \int_{\Gamma_\varrho}
\varrho^3 uv\, ds \quad \; \; \mbox{for} \;\; u\;\; \mbox{and}\;\;
 v \;\; \mbox{in}\;\; K(D)\;\;\mbox{or}\;\; K^d (D)\end{eqnarray*}
 and \begin{eqnarray*} \langle G'u, v\rangle' = \int_{\Gamma_\varrho}
 \varrho'^3
  u v\, ds
\quad \; \; \mbox{for} \;\; u\;\; \mbox{and}\;\;
 v \;\; \mbox{in}\;\; K(D)\;\;\mbox{or}\;\; K^d
 (D),\end{eqnarray*}
 respectively.
 Let \begin{eqnarray*}
 \mu_1\ge \mu_2\ge\cdots\ge \mu_k \ge \cdots >0\quad \,
 \mbox{and}\;\;
\mu'_1\ge \mu'_2\ge\cdots\ge \mu'_k \ge \cdots >0\end{eqnarray*} be
the positive eigenvalues of $G$ and $G'$, respectively.
 If $\varrho\le \varrho'$, then
 \begin{eqnarray} \label{a-1} \mu_k\le \mu'_k \quad \;\mbox{for}\;\;
 k=1,2,3,\cdots.\end{eqnarray}}

\vskip 0.22 true cm

\noindent  {\bf Proof.} \ Since $\varrho\le \varrho'$, we see that
for any $u\in K(D)$ or $K^d(D)$,
\begin{eqnarray*} \frac{\langle Gu, u\rangle}{\langle u, u\rangle}
=\frac{\int_{\Gamma_\varrho} \varrho^3 u^2 \, ds }{\langle u,
u\rangle^* +\int_{\Gamma_\varrho} \varrho^3 u^2\, ds}\le
\frac{\int_{\Gamma_\varrho} \varrho'^3 u^2 \, ds }{\langle u,
u\rangle^* +\int_{\Gamma_\varrho} \varrho'^3  u^2\, ds}
=\frac{\langle G'u, u\rangle'}{\langle u, u\rangle'},\end{eqnarray*}
which
 implies (\ref{a-1}).  \ \ $\square$

 \vskip
0.586 true cm

\noindent  {\bf Proof  of Theorem 1.1.} \ \  a)  \ \  First, let
$(\mathcal{M}, g)$ be a real analytic
 Riemannian manifold, and let the boundary $\partial D$
 of $D$ be $C^{2,\varepsilon}$-smooth.
    As in \cite{Sa}, we divide the domain $\bar D$
 into subdomains in the following manner.
  It is clear that  $\partial D$
  is the union of the
 portions $\bar \Gamma_1, \cdots, \bar \Gamma_p$ (without common
 inner point on the surface).
    Let $U$ be a coordinate neighborhood which contains
 $\bar \Gamma_j$, let $x_i=x_i (Q)$ and $a_i= a_i(\nu_Q)$ be the coordinates
 of a point $Q$ in $\bar \Gamma_j$
  and the interior Riemannian normal $\nu_Q$ at $Q$, respectively.
   We define the subdomain  $D_j$ and surface $\Gamma_j^\sigma$ by
  \begin{eqnarray*} \label{6-1}   D_j = \{P\big| x(P)=x(Q)+\xi_n a(\nu_Q),
  \quad  Q\in \Gamma_j, \,\, 0<\xi_n< \sigma\}\end{eqnarray*}
  and \begin{eqnarray*} \Gamma_j^\sigma =\{ P\big| x(P)=x(Q)+\sigma \,
    a(\nu_Q),\; Q\in \Gamma_j\},\nonumber
  \end{eqnarray*}
 where $\sigma$ is a positive constant. The closure of  $D_j$ is
 \begin{eqnarray} \bar D_j = \{P\big| x(P) =x(Q)+\xi_n a(\nu_Q),
 \; Q\in \bar \Gamma_j, \; 0\le \xi_n \le \sigma\}.\end{eqnarray}
By the assumption,
 each $\bar \Gamma_j$, which is contained in a
 coordinate neighborhood, can be represented by equations
  \begin{eqnarray} \label {0006-1} x_i =\psi_i (\xi_1, \cdots,
  \xi_{n-1})\end{eqnarray}
  with $C^{2,\varepsilon}$-smooth functions $\psi_i$, i.e., it is the imagine of the closure
  $\bar \Upsilon_j$ of an open domain $\Upsilon_j$ of ${\Bbb R}^{n-1}$.
Hence, if $\sigma$ is sufficiently small, the definitions have a
sense  and the formula \begin{eqnarray}  x(P)=x(Q) +\xi_n a(\nu_Q),
\quad \,
  Q\in \bar \Gamma_j, \;\; 0\le \xi_n \le \sigma\end{eqnarray}
  defines a $C^{2,\varepsilon}$-smooth homeomorphism of a neighborhood of the
  image of $\bar D_j$ in ${\Bbb R}^n$  given by the coordinates $x$
  and a neighborhood $U_j$ of the  closed cylinder
  $\bar F_j$ in ${\Bbb R}^n$ defined by $\bar F_j= \{\xi \big| (\xi_1,
  \cdots, \xi_{n-1})\in \bar \Upsilon_j, \, 0\le \xi_n \le \sigma\}$.
   (Therefore all Lemmas of Section 2 are
 true for every $\bar D_j$). Moreover, the domains $\bar D_1, \cdots, \bar D_p$
 have no common inner points and the remainder
 $D_0=D -\cup_{j=1}^p \bar D_j$ of $D$ has a finite number of connected parts.
 Note that  the boundary of $\bar D_0$ contains
 no part of $\Gamma_\varrho$ with measure $>0$.

Define the space $N=N(D)$, $\mathcal{N}$ and the transformation $G$
on $\mathcal N$ as in Section 3.
  We shall investigate the asymptotic behavior of $A(\tau)$ with
  regard to transformation $G$ on space $\mathcal{N}$.
 Moreover, we define the function spaces
   \begin{eqnarray*} & K_j^0 =\{u_j\big|u_j \in Lip(\bar D_j)\cap
   H^2(D_j),\,
\frac{\partial u_j}{\partial \nu}=0 \;\; \mbox{on}\;\; \Gamma_j \cup
\Gamma_j^{\sigma}, \, \mbox{and}\;\; u=0\;\;\mbox{on}\;\;
\partial D_j -\Gamma_j \},\\
   & H_0^0 = \{u_0 \big|Lip(\bar D_0)\cap H^2(D_0), \;
    u=\frac{\partial u}{\partial \nu} =0 \;\; \mbox{on}\;\;
   \partial D_0\},\\
   & K_j^d =
  \{u_j\big|u_j\in   Lip (\bar D_j)\cap H^2(D_j), \, \frac{\partial u_j}{\partial \nu}=0
  \;\; \mbox{on}\;\;
  \Gamma_j, \, u_j=\frac{\partial u_j}{\partial \nu}=0\;\;\mbox{on}\;\;
  \Gamma_j^\sigma\},\\
 & \, \quad (j=0, 1, \cdots, p), \end{eqnarray*}
     and the bilinear functionals
     \begin{eqnarray}  \langle u_j, u_j \rangle_j^{\star} =\int_{D_j} |\triangle_g u_j |^2 dR,
     \quad \, (j=0, 1, \cdots, p), \end{eqnarray}
      \begin{eqnarray} [u_j, v_j]_j = \int_{\Gamma_j} \varrho^3  u_j v_j \, ds,
      \quad (j=1, \cdots, p),  \quad \, [u_0, v_0]=0, \end{eqnarray}
   and
  \begin{eqnarray} \langle u_j, v_j\rangle_j=  \langle u_j,
  v_j\rangle_j^{\star} +[u_j, v_j]_j, \quad (j=0, 1, \cdots, p),\end{eqnarray}
  where $u_j, v_j \in K_j^0$ or $K_j^d$. Closing $K_j^0$ and $K^d_j$ with respect to
 the norm $|u_j|_j=\sqrt{\langle u_j, u_j\rangle_j}$,
 we get the Hilbert spaces ${\mathcal{K}}_j^0$ and ${\mathcal{K}}_j^d$,
  $\,(j=0, 1, \cdots, p)$.
 Then, in the same manner as in Section 5 we can define the  Hilbert
 ${\mathcal{K}}^0$ and  ${\mathcal{K}}^d$, and define the
 positive, completely continuous transformations
  $G^0$, $G^d$, $G_j^0$
 and $G_j^d$ on  ${\mathcal{K}}^0$,
    ${\mathcal{K}}^d$, ${\mathcal{K}}_j^0$
    and ${\mathcal{K}}_j^d$, respectively. Furthermore, we can prove
 \begin{eqnarray} \label{6-100} A^0(\tau) \le A(\tau)\le A^d(\tau)
 \quad \, \mbox{for all}\;\; \tau,\end{eqnarray}
and \begin{eqnarray} \label{6-101} A^0(\tau) =\sum_{j=0}^p
A_j^0(\tau),\quad\;
  A^d(\tau)=\sum_{j=0}^p A_j^d(\tau),\end{eqnarray}
 where $A^0(\tau)$, $A^d(\tau)$, $A_j^0(\tau)$ and $A_j^d(\tau)$
 are the numbers of eigenvalues
 of the transformations  $G^0$, $G^d$, $G_j^0$ and $G^d_j$
 on ${\mathcal{K}}^0$, ${\mathcal{K}}^d$, ${\mathcal{K}}_j^0$ and
 ${\mathcal{K}}_j^d$
  which are greater than or equal to $(1+\tau^3)^{-1}$,
  respectively.

  Since  $[u_0, u_0]_0=0$ for all $u_0, v_0\in K_0^0$ or $K_0^d$  and  $\langle
  G_0^0 u_0, u_0\rangle_0 =\langle G_0^d u_0, u_0\rangle_0 =[u_0, u_0]_0$, we find fairly
  easily that $G_0^0=G_0^d=0$,
   so that  $A_0^0(\tau) =A_0^d(\tau) =0$, ($\tau\ge 0$).
  Thus, we have to  estimate $A_j^0(\tau)$ and $A_j^d(\tau)$ for those domains $D_j$, where
  $\int_{\Gamma_j} \varrho^3 ds>0$.

 We can choose a finer subdivision of $\partial D$ by
subdividing the domains $\bar \Upsilon_j$ into smaller  ones, e.g.
 by means of a cubical net in the coordinates $\xi$.
 According to p.$\,$71 of \cite{Sa}, by a linear transformation of the coordinates
 we can choose a new coordinate system $\eta$ such that
 \begin{eqnarray*}
 g^{il} (\bar \eta)=\delta^{il},\;\quad (i,l=1, \cdots, n),
\end{eqnarray*}
 for one point $\bar \eta$  in the mapping $T_j$ of $\Upsilon_j$,
and \begin{eqnarray}
 \label {0.6-0.1}   & x_k(P)= \psi_k (\eta_1, \cdots, \eta_{n-1})+\eta_n
 \,a(\nu (\eta_1, \cdots, \eta_{n-1})),\\
   & \mbox{for}\;\; (\eta_1, \cdots, \eta_{n-1})\in \bar T_j,\;\,
    0\le \eta_n \le \sigma  \nonumber \end{eqnarray}
 defines a $C^{2,\varepsilon}$-smooth homeomorphism from  $\bar E_j$ to $\bar
 D_j$, where $\bar E_j=\{\eta=(\eta_1$,$\cdots$, $\eta_n)\big| (\eta_1$,
 $\cdots$, $\eta_{n-1})\in \bar T_j, \,
  0\le \eta_n\le \sigma\}$ is a  cylinder
in ${\Bbb R}^n$. (This can also be realized by choosing a
(Riemannian) normal
 coordinates system at the point $\bar \eta\in T_j$ for the manifold
 $(\mathcal{M},g)$ (see, for example, p.$\,$77 of \cite{Le})
  such that $a(\nu(\eta))=(0, \cdots, 0, 1)$ and by using the mapping
  (\ref{0.6-0.1})).
 If we denote the new
 subdomains of
 $\partial D$ by $\bar \Gamma_j$ as before, it is clear that we can always choose
  them and $\sigma$ (see p.$\,$71 of \cite{Sa}), so that,
\begin{eqnarray} \label{0.6-0.2} |g'^{il} (\eta')-g^{il}(\bar
\eta)|<\epsilon, \quad \, i,l=1,\cdots, n,\end{eqnarray}
 \begin{eqnarray} \label {0006-2}  &&\bigg| \frac{1}{\sqrt{|g(\eta')|}}
  \,\frac{\partial}{\partial x_i}
   \big(\sqrt{|g(\eta')|} g^{il}(\eta')\big)\,- \,
   \frac{1}{\sqrt{|g(\bar \eta)|}} \,\frac{\partial}{\partial x_i}
   \big(\sqrt{|g(\bar \eta)|} g^{il}(\bar \eta)\big)\bigg|<\epsilon, \\
   &&\,\qquad \;\;\qquad  \; \qquad \, i,l=1,\cdots,n,\nonumber
   \end{eqnarray}
for any given $\epsilon>0$, and all points $\eta' \in {\bar E}_j$.
 The inequalities (\ref{0.6-0.2}) imply that
\begin{eqnarray} \label {0006-1} (1+{\tilde M_j}\epsilon)^{-1}\sum_{i=1}^n
t_i^2  \le \sum_{i,l=1}^n g^{il}
    (\eta') t_i t_l \le (1+{\tilde M_j}\epsilon)\sum_{i=1}^n
   t_i^2\end{eqnarray}
  for all points $\eta'\in {\bar E}_j$ and all real numbers $t_1, \cdots,
  t_n$, where ${{\tilde M}}_j$ is a positive constant depending only on
  $g^{il}$ and $\bar E_j$ (cf. Lemma 6.1).  This and formula (128) of \cite {Sa} say that
  \begin{eqnarray} \label{6>>1} (1+\epsilon{\tilde M}_j)^{-n/2} |T_j|\le |\Gamma_j|
  \le (1+\epsilon{\tilde M}_j)^{n/2} |T_j|,\end{eqnarray}
  where
  \begin{eqnarray*}  |\Gamma_j|=
 \int_{T_j} \sqrt{g(\eta)}\, d\eta_1 \cdots d\eta_{n-1},
 \quad \, |T_j|=\int_{T_j} d\eta \cdots d\eta_{n-1}\end{eqnarray*}
 are the Riemannian and Euclidean areas of $\Gamma_j$ and
  $|T_j|$, respectively,
  \vskip 0.22 true cm

  Now, we consider the Hilbert spaces ${\mathcal{K}}_j^0$ and
  ${\mathcal{K}}_j^d$. When transported to $\bar E_j$, the underlying incomplete function spaces
  $K_j^0$ and $K_j^d$ are
  \begin{eqnarray*} K_j^0 = \{ u\big| u\in Lip (\bar E_j)\cap H^2(E_j),  \,
    \, \frac{\partial u_j}{\partial \nu}=0 \;\;\mbox{on}\;\; T_j\cup T_j^\sigma,
    \,\mbox{and}\;\,
     u=0\;\;\mbox{on}\;\; \partial E_j -T_j\}
  \end{eqnarray*}
   and \begin{eqnarray*}  K_j^d =\{u_j\big| u_j\in
    Lip (\bar E_j)\cap H^2(E_j),\; \frac{\partial u_j}{\partial \nu}=0
    \;\;\mbox{on}\;\; T_j,\; u_j=\frac{\partial u_j}{\partial \nu}=0\;\;
  \mbox{on}\;\; T_j^{\sigma}\}, \end{eqnarray*}
   respectively.
    The inner product, which is similar to Section 5, is defined  by
    \begin{eqnarray*} \langle u, v\rangle_j = \int_{E_j} (\triangle_g u)(\triangle_g v)
     \sqrt{g(\eta)} d\eta_1 \cdots d\eta_n +\int_{T_j}
  \varrho^3 uv \sqrt{g(\eta)}\, d\eta_1\cdots \eta_{n-1} \end{eqnarray*}
  and the transformations $G_j^0$ and $G_j^d$ are defined by
  \begin{eqnarray*} \langle G_j^0u, v\rangle_j = \int_{T_j}
  \varrho^3
    uv \sqrt{g(\eta)}\, d\eta_1\cdots \eta_{n-1}, \quad \; \mbox{for}\;\; u,v \;\;
  \mbox{in}\;\; {\mathcal{K}}_j^0, \end{eqnarray*}
and
 \begin{eqnarray*} \langle G_j^d u, v\rangle_j = \int_{T_j}
 \varrho^3
    uv \sqrt{g(\eta)}\, d\eta_1\cdots \eta_{n-1}, \quad \; \mbox{for}\;\; u,v \;\;
  \mbox{in}\;\; {\mathcal{K}}_j^d, \end{eqnarray*}
respectively.

 Let
\begin{eqnarray}  {\underline{\varrho}}_j =\inf_{\bar \Gamma_j} \varrho \quad\; \mbox{and}\;\;
 \bar \varrho_j =\sup_{\bar \Gamma_j} \varrho,\end{eqnarray}
 and let us introduce the inner products
  \begin{eqnarray*} {\underline{\langle u, v\rangle}}_j =\int_{E_j}
  (\triangle u)(\triangle v) d\eta_1 \cdots d \eta_n +
  \int_{T_j} {\underline{\varrho}}_j^3 uv\, d\eta_1\cdots d\eta_{n-1} \end{eqnarray*}
 and
 \begin{eqnarray*} {\overline{\langle u, v\rangle}}_j =\int_{E_j}
  (\triangle u)(\triangle v) d\eta_1 \cdots d \eta_n +
  \int_{T_j} {\overline{\varrho}}_j^3  u v\, d\eta_1\cdots d\eta_{n-1} \end{eqnarray*}
in the spaces $K^0_j$ and $K_j^d$, respectively. By closing these
spaces in the corresponding norms, we get Hilbert spaces
${\underline{\mathcal{K}}}_j^0$
 and ${\overline{\mathcal{K}}}_j^d$. Furthermore, we  obtain the
  positive, completely continuous transformations ${\underline{G}}_j^0$ and
  ${\overline{G}}_j^d$ on
 ${\underline{\mathcal{K}}}_j^0$ and
${\overline{\mathcal{K}}}_j^d$, which are given by
\begin{eqnarray} \underline{\langle {\underline{G}}_j^0 u, v\rangle_j}=\int_{T_j}
{\underline{\varrho}}_j^3  uv \, d\eta_1 \cdots d\eta_{n-1}, \quad
\; \mbox{for}\;\; u\;\; \mbox{and}\;\; v\;\; \mbox{in}\;\;
 {\underline{\mathcal{K}}}_j^0\end{eqnarray}
 and
\begin{eqnarray} \overline{\langle {\overline{G}}_j^d u, v\rangle_j}=\int_{T_j}
 {\overline{\varrho}}_j^3 uv
 \, d\eta_1 \cdots d\eta_{n-1}, \quad \; \mbox{for}\;\; u\;\; \mbox{and}\;\; v\;\; \mbox{in}\;\;
 {\overline{\mathcal{K}}}_j^d,\end{eqnarray}  respectively.

  Let $\mu_k(G_j^0)$ be the $k$-th positive eigenvalue of
  $G_j^0$ and so on. According to Lemma 6.1 and Remark 6.2,
  $\Lambda_1^d(D_j)$ and $\Theta_1^d(D_j)$ have uniformly positive lower bound when
  repeated taking finer division of $D$
  (In fact, by repeated halving the side length of
  every rectangular parallelepiped
  in the partition net of the coordinates $\eta$ for each cylinder $E_j$,
  we see that ${\tilde \Lambda}_1^d(D_j)$ will tend to $+\infty$, and
  that  $\Theta_1^d(D_j)$ will have a positive lower bound).
 This implies that the corresponding positive constants ${\tilde M}_j$ and
 $M_j$ have uniformly upper bound when
 we further divide the domain $D$ into finer a division, where
$\tilde M_j$ is defined as before, and $M_j$ is a constant
independent of $\epsilon$ and depending only on $g$, $\partial g$
   and $\bar E_j$ as in Lemma 6.1.
  Denote by $c_j(\epsilon)$ the maximum value of
 $(1+\epsilon {\tilde M}_j)^{n+\frac{1}{2}}\big(
 \min\{(1-\epsilon M), (1+\epsilon \tilde M)^{-\frac{1}{2}}\}\big)^{-1}$ and $
   (1+\epsilon{\tilde M}_j)^{n+\frac{1}{2}} \big(\max\{(1+\epsilon M),
   (1+\epsilon \tilde M)^{\frac{1}{2}}\}\big)$.
Obviously, $c_j(\epsilon)\to 1$ as $\epsilon\to 0$.
  By virtue of (\ref{0.6-0.2}) and (\ref{006-2}), it follows from Lemmas 6.1 and 6.3 that
\begin{eqnarray} \mu_k(G^d_j)\le c_j(\epsilon)\, \mu_k(\bar G_j^d)\end{eqnarray}
and
\begin{eqnarray} \mu_k(G^0_j)\ge c_j(\epsilon)^{-1}\mu_k({\underline{G}}_j^0),\end{eqnarray}
 so that
\begin{eqnarray*} A_j^d (\tau)\le {\bar A}_j^d
\big( c_j(\epsilon)\tau+c_j(\epsilon)-1\big)\end{eqnarray*} and
\begin{eqnarray} A_j^0 (\tau)\ge {\underline{A}}_j^0
\big(c_j(\epsilon)^{-1}\tau+
  c_j(\epsilon)^{-1} -1\big)\end{eqnarray}
 where ${\bar
 A}^d_j (\tau)$ and
${\underline{A}}^0_j (\tau)$ are
 the numbers of eigenvalues  of the transformation $\bar G_j^d$
 and ${\underline{G}}_j^0$
  which are greater than or equal to
 $(1+\tau^3)^{-1}$, respectively.
 By (\ref{6-100}) and (\ref{6-101}), we obtain
\begin{eqnarray} \label{6>1} && \sum_j {\underline{A}}_j^0
\big(c_j(\epsilon)^{-1} \tau+ c_j(\epsilon)^{-1}
-1\big)\le A(\tau) \\
 && \; \quad \qquad \le \sum_j {\bar A}_j^d
\big(c_j(\epsilon)\tau +
c_j(\epsilon)-1\big).\nonumber\end{eqnarray}

 Finally,   we shall apply the results of Section 5 to
estimate ${\underline{A}}_j^0(\tau)$
 and ${\overline{A}}_j^d(\tau)$.
Note that \begin{eqnarray} \label{5:0} {\overline{A}}_j^d(\tau)\le
{\overline{A}}_j^f(\tau) \quad\;\,\mbox{for all}\;\;
\tau>0,\end{eqnarray} where ${\bar A}_j^f$ is defined similar to
(\ref{5.*.15})---(\ref{5..0}).
 Formula (\ref{5-2-10}), (\ref{5..58}) and  (\ref{5:0}) imply that
\begin{eqnarray}\label{6>2} \underset {\tau\to
+\infty}{\underline{\lim}}
 \frac{{\underline{A}}_j^0 (\tau)}{\tau^{n-1}}\ge  \omega_{n-1}(\sqrt[3]{16}\,\pi)^{-(n-1)}
 |{T}_j| {\underline{\varrho}}_j^{n-1}    \end{eqnarray}
and \begin{eqnarray} \label{6>3} \underset {\tau\to
+\infty}{\overline{\lim}} \frac{{\bar A}_j^f (\tau)}{\tau^{n-1}} \le
\omega_{n-1}(\sqrt[3]{16}\,\pi)^{-(n-1)}
  |{T}_j| {\bar \varrho}_j^{n-1},\end{eqnarray}
where $|{T_j}|$ is the Euclidean area of $T_j$. Combining
(\ref{6>1}), (\ref{5:0}), (\ref{6>2}), (\ref{6>3}), (\ref{6-100}),
(\ref{6-101}), (\ref{5.-03}) and (\ref{5.63})  we find
\begin{eqnarray*}  \underset{\tau\to \infty} {\overline{\lim}} A(\tau)\, \tau^{-(n-1)} \le
  \omega_{n-1}(\sqrt[3]{16}\,\pi)^{-(n-1)}{\tilde c}_j(\epsilon) \sum_j
   {\overline{\varrho}}_j^{n-1} |\Gamma_j|,\end{eqnarray*}
and \begin{eqnarray*} \underset{\tau\to \infty} {\underline{\lim}}
A(\tau)
 \,\tau^{-(n-1)} \ge
  \omega_{n-1}(\sqrt[3]{16}\,\pi)^{-(n-1)}{\tilde c}_j(\epsilon)^{-1} \sum_j
   {\underline{\varrho}}_j^{n-1} |\Gamma_j|,\end{eqnarray*}
where ${\tilde c}_j(\epsilon) =(1+\epsilon{\tilde M}_j )^{n/2}
 {c}_j(\epsilon)^{n-1}$. Note that $\varrho$ is Riemannian integrable
since
 it is non-negative bounded measurable function on $\partial D$.
   Therefore, letting $\epsilon\to 0$, we obtain the desired result that
\begin{eqnarray}  \label {6..0} A(\tau)\sim \frac{\omega_{n-1}\tau^{n-1}}
{(\sqrt[3]{16}\,\pi)^{(n-1)}}\int_{\partial D}
  \varrho^{n-1} ds \quad \; \mbox{as}\;\; \tau\to +\infty. \end{eqnarray}

 b) \ \  Since  a metric $g^{il}$
  in $C^{2}$ can be approximated by a metric $g'^{il}_\epsilon$
  which is $C^2$-smooth on $M$ and piecewise
  $C^{2,\varepsilon}$-smooth
  in any compact submanifold of $(\mathcal{M}, g)$ such
 that
 \begin{eqnarray*} |g'^{il}_\epsilon - g^{il} |<
 \epsilon, \quad i,l=1,\cdots, n, \end{eqnarray*}
   \begin{eqnarray*} \bigg| \frac{1}{\sqrt{|g'_\epsilon|}} \frac{\partial}{\partial x_i}
 \big(\sqrt{|g'_\epsilon|} \,g'^{il}_\epsilon \big)-\frac{1}{\sqrt{|g|}} \frac{\partial}{\partial x_i}
 \big(\sqrt{|g|}\, g^{il}\big)\bigg|< \epsilon, \quad \, i,l=1,\cdots, n, \end{eqnarray*}
  for all points in  $\bar D$, with any given $\epsilon>0$.
 In addition, any bounded domain $D$ with $C^{1}$-smooth  boundary
 can also be approximated (see, the definition
  in Section 2) by
  domain $D'_\epsilon$ with $C^1$-smooth and  piecewise
  $C^{2,\varepsilon}$-smooth
 boundary. Thus, by Lemma 6.1 and a), we can estimate the eigenvalues
  for $g'^{il}_\epsilon$ in $D'_\epsilon$. But for these  eigenvalues (\ref{6..0}) is true, so
  that letting $\epsilon \to 0$ and noticing that $ds'_\epsilon\to ds$, we get
  that (\ref{6..0}) also holds for such metric $g^{il}$ and $D$.
   \  \  $\quad \;\; \square$

 \vskip 0.48 true cm

\noindent  {\bf Proof of Corollary 1.2.} \ \  By (\ref{1-6}), we
have
\begin{eqnarray} \label{07-1}  A(\lambda_k)\sim
\frac{\omega_{n-1}\lambda_k^{n-1}}{(\sqrt[3]{16}\,\pi)^{(n-1)}}
   \big(\mbox{vol}(\partial D)\big), \quad \, \; \mbox{as}\;\; k\to +\infty.\end{eqnarray}
 Since  $A(\lambda_k)=k+2$, we obtain (\ref{1-7}), which completes the proof.
 \  \ $\; \square$

\vskip 0.5 true cm

\noindent  {\bf Proof  of Theorem 1.3.} \ \
   Let $T_\epsilon: \, H^{1/2}(\partial D)\to H^{-1/2}
 (\partial D)$ be defined as follows:
 For any $\phi\in H^{1/2}(\partial D)$, we put $T_\epsilon\phi
 =\left(\frac{1}{(\varrho+\epsilon)^3}\,\frac{\partial (\triangle_g u)}{\partial \nu} \right)\big|_{\partial D}$,
 where $u$
 satisfies \begin{eqnarray*}\left\{ \begin{array}{ll} \triangle_g^2 u=0 \;\;\mbox{in}\;\; D,\\
  \frac{\partial u}{\partial \nu}=0\;\; \mbox{on}\;\; \partial D,\\
   u=\phi \;\;\mbox{on}\;\; \partial
 D,\end{array}\right.\end{eqnarray*}
 and $\epsilon>0$ is a sufficiently small constant.
  Clearly, $T_\epsilon$ is a self-adjoint, elliptic,  pseudodifferential
  operator of order $1$.
  Also, its principal symbol has the form $c(\varrho +\epsilon)^{-3}|\xi|^3$,
  where $\xi\in {\Bbb R}^{n-1}$ and $c$ is a unknown constant
  (we will determine it later).
  It is easily seen that the operator $T_\epsilon$ has the same eigenvalues $\big(\lambda_k(\epsilon)\big)^3$
  (respectively, the corresponding eigenfunctions $u_k$) as the
 following Steklov eigenvalue problem:
 \begin{eqnarray*} \left\{ \begin{array}{ll}  \triangle_g^2 u_k=0
 \;\;\mbox{in}\;\; D,\\
  \frac{\partial u_k}{\partial \nu}=0 \;\ \mbox{on}\;\; \partial D,\\
   \frac{\partial (\triangle_g u_k)}{\partial \nu} -\big(\lambda_k(\epsilon)\big)^3
   (\varrho+\epsilon)^3 u_k=0\;\;\mbox{on}\;\; \partial
 D.\end{array} \right. \end{eqnarray*}
  Applying the standard method (see, for example, p.$\;$163 of \cite{Sh}, \cite{Ho}, \cite{Ta3} or \cite{AGMT}), we
  have
\begin{eqnarray*} A^{(\epsilon)}(\tau)&=&\frac{1}{(2\pi)^{n-1}}\left(\int_{\partial D} ds\int_{c(\varrho
+\epsilon)^{-3}
|\xi|^3<\tau^3} d\xi\right)+O(\tau^{(n-2)})\\
&=&\frac{1}{(2\pi)^{n-1}}\left(\int_{\partial D}ds
\int_0^{\frac{(\varrho+\epsilon)\tau}{c^{1/3}}}\sigma_{n-1}\,r^{n-2}
dr\right)+O(\tau^{n-2})\\&=&
\frac{\omega_{n-1}\tau^{n-1}}{(2\pi)^{n-1}}
 \int_{\partial D} \left(\frac{\varrho (s)+\epsilon}{c^{1/3}}\right)^{n-1} ds +O(\tau^{(n-2)})
 \quad \;
\mbox{as}\;\; \tau\to +\infty,\end{eqnarray*} where
$A^{(\epsilon)}(\tau)=\#\{\big(\lambda_k(\epsilon)\big)^3 \le
\tau^3\}$,  $\sigma_{n-1}$ and
$\omega_{n-1}=\frac{\sigma_{n-1}}{n-1}$ are respectively the area
and volume of the unit sphere and the unit ball of ${\Bbb R}^{n-1}$.
Letting $\epsilon\to 0$, we obtain
\begin{eqnarray*} A(\tau)=\frac{\omega_{n-1}\tau^{n-1}}{(2\pi)^{n-1}}
 \int_{\partial D} \left(\frac{\varrho(s)}{c^{1/3}}\right)^{n-1}\,ds +O(\tau^{(n-2)})
 \quad \;
\mbox{as}\;\; \tau\to +\infty.\end{eqnarray*}
 It follows from the leading asymptotic formula (\ref{1-6}) of Theorem 1.1 that
 the constant $c$ must be $2$. Therefore,
\begin{eqnarray*} A(\tau)=2^{-(\frac{n-1}{3})}\cdot (2\pi)^{-(n-1)}\omega_{n-1}\tau^{n-1}
 \int_{\partial D} \rho^{n-1}(s) ds +O(\tau^{(n-2)})
 \quad \; \mbox{as}\;\; \tau\to +\infty.\end{eqnarray*}
 \qed

 \vskip 0.35 true cm

\noindent  {\bf Remark 6.4.} \
 It is worth noting that for the harmonic Steklov eigenvalue
  problem
 (\ref{1-3}) with smooth boundary $\partial D$, we can also obtain
   a better asymptotic formula with remainder estimate than that of \cite{Sa} (i.e., (\ref{1-4})).
   In fact,
   we first define the pseudodifferential operator
 $S_\epsilon: \, H^{1/2}(\partial D)\to H^{-1/2}
 (\partial D)$ as follows:
 For any $\phi\in H^{1/2}(\partial D)$, we put
 $S_\epsilon \phi:=\left((\varrho +\epsilon)^{-1} \frac{\partial v}{\partial \nu}
  \right)\big|_{\partial D}$,
 where $v$ satisfies \begin{eqnarray*}  \left\{ \begin{array}{ll} \triangle_g v=0 \;\;\mbox{in}\;\; D,\\
    v=\phi \;\;\mbox{on}\;\; \partial
 D.\end{array}\right.\end{eqnarray*}
  This is just the well-known Dirichlet-to-Neumann map, and it has the same eigenvalues
 as the Steklov eigenvalue problem
 (\ref{1-3}) with $\varrho$ being replaced by $\varrho+\epsilon$.
    One knows (see, for example, p.$\;$103 of \cite{SU}) that the principal symbol
    of the pseudodifferential operator $T$
  is $(\varrho+\epsilon)^{-1} |\xi|$.
 So, we have  the following
 asymptotic formula (see, for example, p$\;$163 of \cite{Sh}, \cite{AGMT})
\begin{eqnarray*} B^{(\epsilon)}(\tau)=(2\pi)^{-(n-1)}\left(\int_{\partial D} ds\int_{(\varrho+\epsilon)^{-1}
|\xi|<\tau} d\xi\right)+O(\tau^{(n-2)})\\
 =\frac{\omega_{n-1}\tau^{n-1}}{(2\pi)^{n-1}}
 \int_{\partial D} (\varrho (s)+\epsilon)^{n-1} ds +O(\tau^{(n-2)})
 \quad \;
\mbox{as}\;\; \tau\to +\infty.\end{eqnarray*}  By letting
$\epsilon\to 0$, we get
\begin{eqnarray}
B(\tau)=\frac{\omega_{n-1}\tau^{n-1}}{(2\pi)^{n-1}}
 \int_{\partial D} \varrho^{n-1}(s)\, ds +O(\tau^{(n-2)})
 \quad \;
\mbox{as}\;\; \tau\to +\infty.\end{eqnarray}

\vskip 1.68  true cm

\vskip 0.32 true cm

\end{document}